\documentclass[11pt]{amsart}
\usepackage{amssymb,amsmath,amsthm,amsfonts,amsrefs} 
\usepackage{enumerate}
\usepackage {latexsym}
\usepackage{graphicx}
\usepackage{epsfig}
\usepackage{bbm}
\usepackage{todonotes}
\usepackage{hyperref}
\usepackage{tikz}
\include{srctex}
\usepackage{xcolor}
\usepackage{marginnote}
\usepackage[letterpaper, top=1.2in, bottom=1.2in, outer=1.1in, inner=1.1in, heightrounded, marginparwidth=2.5cm, marginparsep=2mm]{geometry}
\usepackage{verbatim}
\numberwithin{equation}{section}

\allowdisplaybreaks
\newcommand{\nc}{\newcommand}
\nc{\les}{\lesssim}
\nc{\nit}{\noindent}
\nc{\nn}{\nonumber}
\nc{\D}{\partial}
\nc{\diff}[2]{\frac{d #1}{d #2}}
\nc{\diffn}[3]{\frac{d^{#3} #1}{d {#2}^{#3}}}
\nc{\pdiff}[2]{\frac{\partial #1}{\partial #2}}
\nc{\pdiffn}[3]{\frac{\partial^{#3} #1}{\partial{#2}^{#3}}}
\nc{\abs}[1] {\lvert #1 \rvert}
\nc{\cAc}{{\cal A}_c}
\nc{\cE}{{\cal E}}
\nc{\cF}{{\cal F}}
\nc{\cP}{{\cal P}}
\nc{\cV}{{\cal V}}
\nc{\cQ}{{\cal Q}}
\nc{\cGin}{{\cal G}_{\rm in}}
\nc{\cGout}{{\cal G}_{\rm out}}
\nc{\cO}{{\cal O}}
\nc{\Lav}{{\cal L}_{\rm av}}
\nc{\cL}{{\cal L}}
\nc{\cB}{{\cal B}}
\nc{\cZ}{{\cal Z}}
\nc{\mR}{{\mathcal R}}
\nc{\mG}{{\mathcal G}}
\nc{\cT}{{\cal T}}
\nc{\cY}{{\cal Y}}
\nc{\cX}{{\cal X}}
\nc{\cXT}{{{\cal X}(T)}}
\nc{\cBT}{{{\cal B}(T)}}
\nc{\vD}{{\vec \mathcal{D}}}
\nc{\efield}{\mathcal{E}}
\nc{\vE}{{\vec \efield}}
\nc{\vB}{{\vec \mathcal{B}}}
\nc{\vH}{{\vec \mathcal{H}}}
\nc{\F}{  \mathcal{F} }
\nc{\ty}{{\tilde y}}
\nc{\tu}{{\tilde u}}
\nc{\tV}{{\tilde V}}
\nc{\Pc}{{\bf P_c}}
\nc{\bx}{{\bf x}}
\nc{\bX}{{\bf X}}
\nc{\bXYZ}{{\bf XYZ}}
\nc{\bY}{{\bf Y}}
\nc{\bF}{{\bf F}}
\nc{\bS}{{\bf S}}
\nc{\dV}{{\delta V}}
\nc{\dE}{{\delta E}}
\nc{\TT}{{\Theta}}
\nc{\dPsi}{{\delta\Psi}}
\nc{\order}{{\cal O}}
\nc{\Rout}{R_{\rm out}}
\nc{\eplus}{e_+}
\nc{\eminus}{e_-}
\nc{\epm}{e_\pm}
\nc{\eps}{\varepsilon}
\nc{\vnabla}{{\vec\nabla}}
\nc{\G}{\Gamma}
\nc{\w}{\omega}
\nc{\mh}{h}
\nc{\mg}{g}
\nc{\vphi}{\varphi}
\nc{\tlambda}{\tilde\lambda}
\nc{\be}{\begin{equation}}
\nc{\ee}{\end{equation}}
\nc{\ba}{\begin{eqnarray}}
\nc{\ea}{\end{eqnarray}}

\providecommand{\abs}[1]{\lvert #1 \rvert}
\nc{\g}{\gamma}
\nc{\ol}{\overline}

\newcommand\Span[1]{\left\langle #1 \right\rangle}
\newtheorem{theorem}{Theorem}[section]
\newtheorem{lemma}[theorem]{Lemma}
\newtheorem{lm}[theorem]{Lemma}
\newtheorem{pr}[theorem]{Proposition}
\newtheorem{prop}[theorem]{Proposition}
\newtheorem{cor}[theorem]{Corollary}

\newtheorem{rmk}[theorem]{Remark}

\nc{\pT}{\partial_T}
\nc{\pz}{\partial_z}
\nc{\pt}{\partial_t}
\nc{\la}{\langle}
\nc{\ra}{\rangle}
\nc{\infint}{\int_{-\infty}^{\infty}}
\nc{\halfwidth}{6.5cm}
\nc{\figwidth}{10cm}
\newcommand{\f}{\frac}
\nc{\nlayers}{L} \nc{\nsectors}{M}
\nc{\indicator}{\mathbf{1}}
\nc{\chile}{R_{\rm hole}}
\nc{\Rring}{R_{\rm ring}}
\nc{\neff}{n_{\rm eff}}
\nc{\Frem}{F_{\rm rem}}
\nc{\R}{\mathbb R}
\nc{\C}{\mathbb C}
\nc{\Z}{\mathbb Z}
\nc{\DD}{\Delta}
\nc{\cD}{\mathcal D}
\nc{\logorm}{\left\|}
\nc{\rnorm}{\right\|}
\nc{\rnormp}{\right\|_{\ell^{p,\eps}}}
\nc{\rar}{\rightarrow}
\DeclareMathOperator{\Ai}{Ai}
\DeclareMathOperator{\Bi}{Bi}
\DeclareMathOperator{\sign}{sign}
\DeclareMathOperator{\supp}{supp}

\nc{\bbR}{\mathbb{R}}
\nc{\calD}{\mathcal{D}} 
\nc{\bbC}{\mathbb{C}}
\nc{\calL}{\mathcal{L}}

\begin{document}

\title[Schr\"odinger equation with a repulsive Coulomb potential]{Pointwise decay for radial solutions of the Schr\"odinger equation with a repulsive Coulomb potential}

\author[A.~Black]{Adam~Black}
\address{Department of Mathematics, Yale University, New Haven, CT 06511} 
\email{adam.black@yale.edu}

\author[E.~Toprak]{Ebru~Toprak}
\address{Department of Mathematics, Yale University, New Haven, CT 06511} 
\email{ebru.toprak@yale.edu}

\author[B.~Vergara]{Bruno~Vergara}
\address{Department of Mathematics, Brown University, Providence, RI 02912} 
\email{bruno\_vergara\_biggio@brown.edu}

\author[J.~Zou]{Jiahua~Zou} 
\address{Department of Mathematics, Rutgers University, Piscataway, NJ 08854} 
\email{jiahua.zou@rutgers.edu}






\maketitle

\begin{abstract}
We study the long-time behavior of solutions to the Schr\"odinger equation with a repulsive Coulomb potential on $\bbR^3$ for spherically symmetric initial data. Our approach involves computing the distorted Fourier transform of the action of the associated Hamiltonian $H=-\Delta+\frac{q}{|x|}$ on radial data $f$, which allows us to explicitly write the evolution $e^{itH}f$. A comprehensive analysis of the kernel is then used to establish that, for large times, $\|e^{i t H}f\|_{L^{\infty}} \leq C t^{-\frac{3}{2}}\|f\|_{L^1}$. Our analysis of the distorted Fourier transform is expected to have applications to other long-range repulsive problems. 
\end{abstract}

\section{Introduction}
\subsection{Motivation and main result}
In this article, we study the dispersive properties of the Coulombic Hamiltonian
\begin{align*}
    H:=-\Delta+\frac{q}{|x|},\,\,q\in\bbR.
\end{align*}
The associated Schr\"odinger equation  \begin{align}\label{eq.Coulomb}
\begin{split}
&i \partial_{t}u=Hu\\
&u(0,x) = f(x)
\end{split}
\end{align}
describes the electrical interaction between two charged particles.
In the attractive case, $ q<0 $, this equation has been used since the birth of quantum mechanics to describe the time-evolution of the electron in a hydrogen atom \cite{Schroe}. In this paper, we focus on the repulsive case, $ q>0 $, which was introduced by Yost, Wheeler and Breit in 1936 \cite{YWB} in connection with the interaction of charges of the same sign. Stationary solutions to \eqref{eq.Coulomb} are sometimes called Coulomb wavefunctions \cite{ESII}.

As is well known, the free Schr\"{o}dinger equation in $ \R^3 $ (that is \eqref{eq.Coulomb} without the potential term) enjoys many dispersive estimates describing the spreading of the wave packet, the most fundamental of which is
\begin{align}\label{disperse}
\| e^{i t \Delta} f\|_{L_{}^{\infty}(\R^{3})}\leq \frac{C}{|t|^{3/2}} \|f\|_{L^{1}(\R^{3})}\,,
\end{align}
with constant independent of $ f $ and $ t $. This is a direct consequence of the fact that the  propagator of the free Schr\"odinger equation in $\R^3$ can be computed as $e^{it\Delta} f=k_{t}* f$, where
\begin{align*}
	k_{t}(x)= (4\pi i t)^{-3/2} e^{\frac{|x|^{2}}{4i t}}.
\end{align*}
Our goal in this paper is to extend the estimate \eqref{disperse} to $e^{itH}$ for radial initial data. Namely, we prove the following theorem:
\begin{theorem}\label{thm.main}
Let $H=-\Delta+\frac{q}{|x|}$ be the Coulomb Hamiltonian in $\R^{3}$, with $q>0$.
Then, for any spherically symmetric function $f$ one has the dispersive estimate 
\begin{equation}\label{eq.dispersive}
\| e^{i t H} f\|_{L_{}^{\infty}(\R^{3})} \leq \frac{C}{|t|^{3/2}} \|f\|_{L_{}^{1}(\R^{3})}\,,\qquad \abs{t}  \geq 1,
\end{equation}
with constant $C>0$ independent of $ f $ and $ t $.
\end{theorem}
As explained in the appendix, $ H $ is a self-adjoint operator on $ H^2(\R^3) $ so $ e^{itH} $ is defined via the spectral theorem. In later work, we plan to extend this result to non-radial functions, see the discussion in Section 1.3. \par
Besides the intrinsic interest of the estimate \eqref{eq.dispersive}, we expect that it will play a role in understanding nonlinear long-range scattering phenomenon. Beyond this, the distorted Fourier transform of $H$ developed for the proof may itself prove useful in such analysis due to the importance of the usual Fourier transform for problems with short-range potentials. In particular, we were motivated by the desire to study the modified scattering of the \emph{Hartree} or \emph{Gross-Pitaevskii} equation in $\R^{3}$:
 $$(i\partial_{t}-\Delta) u+\left((-\Delta)^{-1} |u|^{2}\right)u=0\,.$$
This is a nonlinear equation modeling the dynamics of non-relativistic bosonic many-body particle systems in the mean-field limit, which has been extensively studied together with the NLS equation, see e.g., \cites{HN,KP}.
Indeed, taking radial perturbations leads to the equation
\[
	\left(i\partial_{t}- \partial^2_r + \frac{1}{r} \right)v+\varepsilon^{2} v\int\limits_{\R_{+}}\frac{1}{\max\{r,s\}} |v(s,t)|^{2}ds=0\,,
\] 
with $0<\varepsilon\ll 1$, whose linearization is given by the equation \eqref{eq.Coulomb} for radial data. Thus, any investigation of the nonlinear scattering phenomenon of the Hartree equation must begin with a detailed understanding of the long-time asymptotics of solutions to \eqref{eq.Coulomb}.

\subsection{Prior work}
The problem of extending pointwise dispersive estimates to Hamiltonians of the form ${-\Delta +V(x)}$ has attracted considerable attention, given that these estimates serve as crucial tools for the subsequent analysis of both linear and nonlinear problems. In the linear setting, they give rise to Strichartz estimates and in the nonlinear realm they can be used to prove the stability of solitons, see e.g \cites{SW,Weder,survey,survey2}. Nevertheless, most of the existing studies rely on studying the resolvent perturbatively or the use of Duhamel's formula. Consequently, these approaches require that the potential $V$ be bounded or small in some suitable sense or decay faster than $|x|^{-1}$ at infinity, see the notable references \cites{BG, EG, GS, JSS, RSchlag, Yaj} for some of the diverse methods employed in this area. From this perspective, the Coulomb potential $ \abs{x}^{-1} $ is pathological because of its singularity at the origin and slow decay at infinity. For instance, while pointwise estimates have been investigated for inverse square potentials \cites{ FFF,KT }, the same methods cannot be directly applied to $H$ due to the differing scaling behavior between the Coulomb potential and the Laplacian.


To the best of our knowledge, there is only one study quantifying the dispersive properties of the Hamiltonian $H$. In \cite{Miz}, Mizutani considers the more general operator 
\begin{align*}
	H_1:= -\Delta + Z |x|^{-\mu} + \varepsilon V_S(x)
\end{align*}
on $\R^n$ where $\mu\in (0,2)$ and $ |\partial_x^{\alpha}\{ V_S(x) \} | \leq C \la x \ra^{-1-\mu -|\alpha|}$, which for $ \varepsilon=0 $ and $ \mu=1 $ is our operator $ H $. For $ \varepsilon\geq 0 $ sufficiently small depending on $ Z,\mu, $ and $ V_S $, they show Strichartz estimates (including the endpoint) for  $H_1$, i.e., 
 \begin{align*}
 \|e^{itH_1} f \|_{L^p(\R:L^q(\R^n))} \leq C \|f\|_{L^2(\R^n)}, \,\,\,\,\,\ \frac{2}{p}= \frac{n}{2} - \frac{n}{q} , \,\,\,\,\,\ (n,p,q)\neq(2,\infty,\infty).
 \end{align*}
We emphasize that our estimates are instead pointwise in that they control the $ L^\infty $ norm for all large $ t $. Furthermore, we explicitly compute the kernel of the evolution operator $ e^{itH} $ for radial data, which may be of independent interest.

\subsection{Overview of the proof}
Though we will mostly focus in radial waves,
let us consider the spherical decomposition of
$L^2(\R^{3}) = \bigoplus_{{\ell} =0}^\infty L^{2}(\R^{+},\,r^2 dr)\otimes L_{{\ell}}$, 
where $r:=|x|$ and
$L_{{\ell}}$ denotes the eigenspace of spherical harmonics of angular momentum $\ell=0,1,2\dots$.
According to this decomposition, any $f \in L^2(\R^3)$ can be represented as 
\begin{equation*}\label{eq.partialwaves}
f(r,\omega) =  r \sum_{\ell =0}^{\infty}
f_{\ell}(r)Y_{\ell,m} (\omega)\,,
\quad f_{\ell}(r):=\sum_{m=-\ell}^{\ell} \la f(r, \cdot), Y_{\ell,m} \ra\,,
\qquad r\in\R^{+}\,,\omega\in\mathbb{S}^{2}\,,
\end{equation*}
where the inner product is understood in the sense of $L^{2}(\mathbb{S}^{2})$,
$f_{\ell} \in L^{2}(\R^{+},\,r^2 dr) $ 
and $Y_{\ell,m}$ denotes the $(\ell,m)$-spherical harmonic,
i.e. $-\Delta_{\mathbb{S}^2} Y_{\ell,m} = \ell (\ell+1) Y_{\ell,m} $.
Because of its radial symmetry, the Coulomb Hamiltonian $H$ preserves each angular momentum sector $L^2(\bbR^+,\,r^2dr)\otimes L_{\ell}$. On this subspace, it is unitarily equivalent to
the Sturm-Liouville operator $- \partial^2_r +\frac{\ell (\ell+1)}{r^2} + \frac{q}{r} $:
\begin{equation*}
H|_{L_{\ell}}=H_{\ell,q} := r^{-1} \Big(- \partial^2_r +\frac{\ell (\ell+1)}{r^2} + \frac{q}{r} \Big) r\,.
\end{equation*}

In this paper, we treat the radial sector, $ \ell=0 $, leaving the analysis of other angular momenta to a subsequent work. As such, we must understand the half-line Schr\"{o}dinger operator 
\begin{align*}
	-\frac{d^2}{dr^2}+\frac{q}{r}\,.
\end{align*}
This operator is not essentially self-adjoint (indeed, it is limit circle at $ r=0 $), so we must be careful to choose the correct self-adjoint extension whose dynamics coincide with that of $ H_{0,q} $, which we denote $ \mathcal{L}_q $. We refer the reader to the appendix for the construction of $ \mathcal{L}_q $ and its relevant properties.\par
To explicitly describe the time-evolution of $ \mathcal{L}_q $, we derive its \emph{distorted Fourier transform}, which diagonalizes the operator. This consists of a distorted Fourier basis of appropriately selected generalized eigenfunctions of $ \mathcal{L}_q $ and an  associated measure $ \rho(d\sigma) $, which yields the representation
\begin{align} \label{Lprop}
 e^{it\mathcal{L}_q}g (r) =\int\limits_0^{\infty}\int\limits_0^{\infty} e^{it \sigma ^2} \phi_q(\sigma,r) \phi_q(\sigma,s)  g(s)  \:ds  \:\rho(d\sigma).  
\end{align}
 While the existence of such a transform is quite classical for many classes of potentials, the potential $ r^{-1} $ is \emph{strongly singular}, so that additional care is required to develop its spectral theory. The distorted Fourier transform of strongly singular potentials has been studied in \cite{GZ}. We also mention work of Fulton \cite{F}, in which Frobenius solutions from zero are used to derive the distorted Fourier transform for such potentials. However, the Coulomb potential is not singular enough to apply the results of \cite{GZ} verbatim. In the appendix, we adapt these results to treat the Coulomb case, and derive a distorted Fourier basis and measure given by \begin{align} \label{phiandmeasure}
 \phi_q(\sigma,r) = (2i\sigma)^{-1} M_{\frac{iq}{2  \sigma}, \f 12} ( 2 i \sigma r), \,\,\,\,\,\,\,  d\rho(\sigma)= 2\mu^2(\sigma) d \sigma \,,
 \end{align} 
 where $ \mu^2(\sigma)= q\sigma [e^{\frac{q\pi}{\sigma}} -1]^{-1}$, $\mu \geq0$  and $M_{\frac{iq}{2 \sigma}, \f 12} ( 2 i \sigma r)$ is the Whittaker M function, see \cite[(13.14)]{NIST}.

 Substituting \eqref{phiandmeasure} into \eqref{Lprop} and employing the relation $H_{0,q} = r^{-1} \mathcal{L}_q r$, we obtain 
\begin{align}\label{eq.kernel} 
e^{itH_{0,q}}g(r)&= \frac{q}{2r}\int\limits_0^{\infty}\int\limits_0^{\infty} e^{it  q^2\sigma^2} e(q\sigma, r) e(q\sigma, s) s g(s)  ds d \sigma
 \\ 
 &=\int\limits_0^{\infty} K_t(r,s) s^2 g(s) ds \nn\,,
\end{align}
where 
\begin{align}\label{dfb}
e(q\sigma,  r)  = 2 \mu(q \sigma) \phi_1(q\sigma,r) ,\,\,\,\,\, K_t(r,s)= \frac{q}{2 rs}\int\limits_0^{\infty} e^{it  q^2\sigma^2} e(q\sigma, r) e(q\sigma, s) d \sigma. 
\end{align}
 Therefore, Theorem \ref{thm.main} holds provided that we establish
\begin{align}\label{eq.kernelest}
\sup_{r\geq s>0} | K_{t}(r,s)|\les t^{-3/2},\,\,\,\,\,\ t\geq 1\,. 
\end{align}

It should be stressed that, despite the apparent simplicity of $ K_t $ in terms of special functions, we require several delicate approximations in order to prove \eqref{eq.kernelest}. In particular, for use in the oscillatory integral defining $ K_t $, it is important that we obtain $ C^2 $-control jointly in the variables $ r $ and $ \sigma $. While for large $ \sigma $ known integral and series representations of the Whittaker M function suffice, as $\sigma \to 0$ these turn out to be useless due to the singularity in the index. Therefore, one needs to perform a detailed analysis of the eigenvalue problem 
\begin{align} \label{eigenfeq}
-\frac{d^2 f}{dr^2} + \frac{f}{r} = \sigma^2 f
\end{align}
as $ \sigma\rightarrow 0 $, in which we have conveniently normalized the charge constant to $q=1$. By rescaling, this loses no generality so we adopt this convention from now on. A fundamental set of solutions to \eqref{eigenfeq} consists of $M_{\frac{i}{2 \sigma}, \frac{1}{2}} (2 i \sigma r)$ and $W_{\frac{i}{2 \sigma}, \frac{1}{2}} (2 i \sigma r)$, the latter being the Whittaker $ W $ function. With $\phi(\sigma,r):=\phi_1(\sigma,r)$, as $r$ approaches 0, for any fixed $\sigma>0$,
\begin{align}\label{phiat0}
\phi(\sigma,r)= r ( 1+ O(\sigma r))\,,
\end{align} 
whereas for $ r\rightarrow\infty $,
by equations (13.14.32), (13.14.21), and (5.4.3) in \cite{NIST},
we have that
\begin{align}\label{MAsy}
	&\phi(\sigma,r)	\sim C_0\sigma^{-\frac{1}{2}} [ e^{\frac{\pi}{\sigma}}-1]^{\f12} \sin(\Theta(\sigma,r))\\
	&\Theta(\sigma,r)=\sigma r-\frac{\log(2\sigma r)}{2\sigma}+\theta_0(\sigma)
\end{align}
for some absolute constant $ C_0 $ and phase correction $ \theta_0(\sigma) $. \par
  Heuristically, one may see these asymptotics via the WKB method. Let
\begin{align}
\rho_{\sigma}(s,r)=\int\limits_s^{r} \sqrt{ \abs{Q}  } du,\,\,\,Q=u^{-1} - \sigma^2\,,
\end{align}
be the \emph{Agmon distance} \cite[Ch.5]{Ag} from $ s $ to $ r $. From $0$ to the \emph{turning point} $ r_*=\sigma^{-2} $, WKB predicts that $ \phi $ will grow as $ e^{\rho_\sigma(0,r_*)} $ before oscillations set in, which are governed by $ e^{i\rho_\sigma(r_*,r)} $. Since
\begin{align*}
	\rho_{\sigma}(0,r_*)=\frac{\pi}{2\sigma}
\end{align*}
and
\begin{align*}
	\rho_{\sigma}(r_*,r)=\sigma^{-1}\left(r\sigma^2+\frac{1}{2}\log(r\sigma^2)+\frac{1}{2}-\log(2)+O(r\sigma^2)\right)\,,	
\end{align*}
the predictions of WKB exactly match the asymptotics \eqref{MAsy}, where $ \sigma^{-\frac{1}{2}} $ may be regarded as the usual WKB prefactor of $ Q^{-\frac{1}{4}} $.  Within equation \eqref{dfb}, the role of the $\mu$ multiplier is precisely to adjust the distorted Fourier basis in a manner that leads to 
$$\lim_{r\to\infty} |e(\sigma,r) - c_0\sin(\Theta(\sigma,r))|=0\,,$$ 
for all $\sigma>0$ and an absolute constant $c_0 \in \R$.

\par 
To make the approximation precise, we use the \emph{Liouville-Green} (LG) transform \cite[Ch. 11]{Olver}, which is standard practice for semi-classical problems with a simple turning point. Indeed, such semi-classical problems with inverse square or exponentially decaying potentials have been studied in \cites{Wmain, DSS, Wexp, survey2} (see Section 2 for a discussion of the differences). The LG transform is a change of variables $ \zeta $ that transforms a second order ODE with a simple turning point into a \emph{perturbed Airy equation} of the form
\[
-\sigma^{2}\frac{d^2f}{d\zeta^2}=(\zeta+\sigma^2V(\zeta)) f
\]
for a suitable potential $V$. By working in the $\zeta$ variable, we obtain fundamental systems of solutions to \eqref{eigenfeq} with good asymptotics in terms of Airy functions on both sides of the turning point, see Proposition~\ref{smallexp} and Proposition~\ref{oscBasisPr}. 

Unfortunately, the approximation we obtain for $e(\sigma,r)$ on the left side of the turning point becomes ineffective as $r\rightarrow 0$ due to the local singularity in $V(\zeta)$ introduced by the Coulomb potential. In fact, referring to Nakamura's investigation in \cite{Nak}, one expects $e(\sigma,r)$ to decay quite rapidly as $\sigma$ approaches 0 near $r=0$. In \cite{Nak}, the resolvent operator of $-\Delta + W$ where $|W(x)| \les \la x \ra^{-\rho}$ for $\rho \in(0,2)$ is studied. It was proved there that if $W$ is also positive, then there exist $\beta, \gamma >0$, such that 
\begin{align}\label{Nakdec}
\| F(|x| \leq \beta \sigma^{-2} ) F((-\Delta+W) \leq \sigma^2 )\|_{L^2 \to L^2} \les \exp(- \gamma \sigma^{ -\frac{2}{\rho}+1} ),   \,\,\,\,  \sigma^2 \in (0,1]\,,
\end{align}
where $F(A)$ denotes the characteristic function of the set $A$. Because the Coulomb potential takes the form of $W$ with $\rho=1$ as $x$ becomes large, we should therefore compare our estimate on $e(\sigma,r)$ with $e^{-\frac{\gamma}{\sigma}}$ as $ \sigma \to 0$. With this in mind, we take a different perspective on problem \eqref{eigenfeq} around zero, and derive an approximation to $e(\sigma,r)$ via modified Bessel functions \cite[Ch.10]{NIST} that captures the expected behavior. Lemma~\ref{lem:bb} below arises as a consequence of this approximation and may be regarded as complementary to the estimate in \eqref{Nakdec}. 

In view of our analysis, it seems that one cannot avoid two connection problems when considering potentials that decay like $ r^{-1} $. Indeed, in \cite{PSV}, the same problem appeared while performing a turning point analysis of a similar ODE arising in the context of the Klein-Gordon equation on a stationary spherically symmetric black hole. The approximations in terms of Airy and Bessel functions obtained in Section 5 of the aforementioned work are similar to those we obtain in Section 2. Our analysis of the large $\sigma$ regime, as well as the oscillatory integrals we consider are different, however.

We also mention that while approximations of Coulomb eigenfunctions via Bessel and Airy functions have appeared in the literature since the 50's \cites{AS,ESI,ESII}, the form of these approximations are completely unsuitable for use inside the oscillatory integral \eqref{eq.kernel}, as we require explicit estimates of $C^{2}$-errors in both $ r $ and the semi-classical parameter $ \sigma $. To the best of our knowledge, the approximations derived in Sections 2 and 3 are novel. In particular, the results in the aforementioned references cannot be applied directly in the proof of estimate \eqref{eq.kernelest}, for which one needs the control of several derivatives in order to extract the time decay from the phase. See Section 4 for the oscillatory estimates that lead to the proof of \eqref{eq.kernelest}.

\subsection{Notation and conventions}

For the benefit of the reader, we define here the notation and conventions we use throughout this paper:
\begin{itemize}
\item $ \la x \ra = (1+x^2)^{\f12}$.
\item  $ f(x) \les g(x)$ indicates  that there exists some constant $C>0$ so $ f(x) \leq Cg(x)$ for all $x$ in some specified domain. We will also use this notation for functions that depend on several variables. 
\item $ f(x) \sim g(x)$ indicates  that $ cg(x) \leq f(x) \leq Cg(x)$ for some $C>c>0$ independent of $x$.
\item $a(x,\sigma) = O_k(\sigma^m x^p )$ indicates that  $|\partial^j_\sigma \{a(x,\sigma)\}| \les \sigma^{m-j}x^p $ for $ j= 0,1,\dots ,k$. 
\item $\chi_c(x):\bbR_+\rightarrow\bbR$ is a smooth cut-off function supported on $[0,c] $ that is equal to $1$ when $x\leq \frac{2c}{3}$, and  $\tilde{\chi}_c(x)= 1 - \chi_c(x)$.
\item For two functions $f(x)$ and $g(x)$, $W[f,g](x)$ denotes their Wronskian
\begin{align*}
    W[f,g](x)=f(x)g'(x)-f'(x)g(x).
\end{align*}
\end{itemize}

\subsection{Organization of the paper}
The paper is organized as follows. In Sections 2 and 3, we derive approximations to the Whittaker M function $ e(\sigma,r) $ for small and large energies $ \sigma $, respectively. We devote Section 4 to the proof of our main Theorem \ref{thm.main}. For this purpose, we employ the previous eigenfunction approximations to estimate the oscillatory integral \eqref{eq.kernelest}, which leads immediately to the estimate \eqref{eq.dispersive}. Finally, in an Appendix, we detail the construction of the distorted Fourier basis \eqref{dfb} and thus demonstrate the form of the kernel $ K_t $.

\subsection*{Acknowledgement} This research was conducted as part of the Brown-Yale PDE seminar. We thank the other participants for stimulating discussions. We are particularly grateful for the organizers, Beno\^{i}t Pausader and Wilhelm Schlag, who suggested this problem and provided guidance in the writing of this paper. We also thank Wilhelm for notes that formed the basis of the Appendix and Haram Ko for helpful revisions to those notes.

\section{Eigenfunction approximation: Small Energies} 

The aim of this section is to find a good approximation for the distorted Fourier basis introduced before,
$e(\sigma,r) = -i \sigma^{-\f12}[e^{ \frac{\pi}{\sigma}} -1]^{-\f 12} M_{ \frac{i}{2 \sigma}, \f 12} ( 2 i \sigma r )$ when $ \sigma <c$ for some sufficiently small $c$. The main results from this section that we employ for the analysis of the oscillatory integrals are summarized in Proposition~\ref{langerBounds}, Corollary~\ref{oscCor}, and Corollary~\ref{cor:low}.

 We construct approximate solutions to the ODE
\begin{align}\label{rODE}
	-g''(r)+\frac{g(r)}{r}=\sigma^2g(r)
\end{align}
for $ \sigma $ small.
Substituting $f(x) = g(\sigma^2 r)$, equation \eqref{rODE} transforms into:
\begin{align}\label{xODE}
	-\sigma^{2}f''(x)+(x^{-1}-1)f(x)=0\,.
\end{align}

As $ \sigma\rightarrow 0 $, \eqref{xODE} is a semi-classical problem with a simple turning point at $ x=1 $. Thus, we use the Liouville-Green transform to obtain solutions in terms of Airy functions, which is standard practice for such ODEs; see \cite[Ch. 11]{Olver}). The outline of this section is close to \cite{Wmain} where a similar analysis is performed for an ODE with a potential with $ r^{-2} $ asymptotics. Unlike in \cite{Wmain}, we cannot effectively use the Airy function approximation in a neighborhood of $ x=0 $. Indeed, the error in the approximation blows up in $\sigma$ as $x \to 0$ due to the exponential growth of $\Bi(z)$ as $z \to \infty$, see Proposition~\ref{smallexp}. Therefore, in this regime we resort to  an approximation in terms of the modified Bessel function $ I_1 $. Note that Erdélyi and Swanson use $ I_1 $ in \cite{ESII} to approximate the Whittaker function $M_{ \frac{i}{2 \sigma}, \f 12} ( 2 i  \sigma r )$, but the error term derived there is not suitable for our purposes. Therefore, our analysis in the next section may be regarded as a refined version of theirs.

\subsection{Bessel function approximation: $ x\ll 1 $ }
In this section, we approximate $e(\sigma,r)$ in terms of the modified Bessel function $ I_1 $ when $x\in  [0,\delta]$ for some $ \delta<1 $. Let $Q(x) = x^{-1} -1$ and for $ 0\leq x\leq 1 $ define 
\begin{align} \label{etadef}
 \eta(x) =\int\limits_0^x \sqrt{Q(y)}\, dy\,.
 \end{align}
The properties of this function are summarized in the following lemma:
\begin{lm}\label{etaLem}
The map $\eta$ is a smooth transformation from $[0,1]$ to $[0,\frac{\pi}{2}]$,
and with respect to the change of variables 
\begin{align*}
 	\omega(\eta)=p^{\frac{1}{2}}f,\,\,p= \frac{\eta'}{\eta} \,,
\end{align*}
Equation \eqref{xODE} transforms to 
\begin{align}\label{bessel}
\sigma^2 \big( \ddot{\omega} (\eta) + \frac{ \dot{\omega} (\eta)}{ \eta}  \big) - \omega(\eta) = \sigma^2  V_-(\eta) \omega(\eta)\,, 
\end{align} 
where 
\begin{align*}
	V_-(\eta)= \eta^{-1}p^{-\frac{1}{2}}\frac{dp^{\frac{1}{2}}}{d\eta}+p^{-\frac{1}{2}}\frac{d^2p^{\frac{1}{2}}}{d\eta^2} 
\end{align*}
Here, $ \dot{ } $ represents the derivative with respect to $\eta$ and $ ' $ the derivative with respect to $ x $.\par 
Furthermore, we may write
\begin{align*}
	V_-(\eta)=\frac{1}{\eta^2}+\tilde{V}_-(\eta)
\end{align*}
for $ \tilde{V}_- $ smooth on any interval of the form $ [0,\delta] $, $ \delta<\frac{\pi}{2} $.
\end{lm}
\begin{proof}
The smoothness of $ \eta $ is clear. Furthermore, one computes that
\begin{align*}
	&\dot{\omega}=\eta^{-1}p^{-\frac{1}{2}}f'+\frac{dp^{\frac{1}{2}}}{d\eta}f\,,\\
		\begin{split}
		\ddot{\omega}
		=\eta^{-2}p^{-\frac{3}{2}}f''-\eta^{-1}\dot{\omega}+\left( \eta^{-1}p^{-\frac{1}{2}}\frac{dp^{\frac{1}{2}}}{d\eta}+p^{-\frac{1}{2}}\frac{d^2p^{\frac{1}{2}}}{d\eta^2} \right) \omega
	\end{split}
\end{align*}
and thus, using that $ \sigma^2f''=Qf $ and $ (\eta')^2=Q $, the above expression for $ \ddot{\omega} $ may be rewritten as \eqref{bessel}.\par

Furthermore, by the chain rule, one can calculate 
\begin{align*}
V_-(\eta) &= \frac{1}{ 4 \eta^2} - \frac{3}{4} \frac{ (\eta^{\prime \prime})^2}{ (\eta^{\prime})^4} + \frac{1}{2} \frac{\eta^{\prime \prime \prime}}{ (\eta^{\prime})^3}  \\
&= \frac{1}{  \eta^2} +\Big[ - \frac{3}{4} \frac{ (\eta^{\prime \prime})^2}{ (\eta^{\prime})^4} + \frac{1}{2} \frac{\eta^{\prime \prime \prime}}{ (\eta^{\prime})^3} - \frac{3}{4 \eta^2} \Big] = \frac{1}{  \eta^2} + \tilde{V}_-(\eta)\,.
\end{align*}
Moreover, since $\eta(x) = 2x^{\f 12} + O_{\infty}( x^{3/2}) $ for $ x < 1$, one has $|\partial^j_{\eta}\tilde{V}_-(\eta)| \les 1$ for $j=0,1,\dots$ for $ \eta < \pi/2$.
\end{proof} 
With Lemma~\ref{xODE} in hand, we look for the solution to \eqref{bessel} that is relevant to $e(\sigma,r)$. Recall that Equation \eqref{rODE} has a basis of solutions given by $M_{\frac{i}{2\sigma},\f12}(2 i x/\sigma )$, which is a multiple of $e(\sigma, x/\sigma^2 )$, and $W_{\frac{i}{2\sigma},\f12}(2 i x/\sigma )$.  For any fixed $\sigma$, the first one vanishes to first order at $x=0$ whereas the latter is non-vanishing there \cite[((13.14.14),13.14.17)]{NIST}. Transforming these solutions under the change of variables defined in Lemma~\ref{etadef}, the relation $ p^{\f12}\sim x^{-\f12}$ shows that \eqref{bessel} must have two linearly independent solutions $\phi_-$ and $\phi_+$ satisfying the asymptotics 
    \begin{align} \label{phi12}
 \phi_-(\sigma,\eta) \sim  \eta \,\,\,\ \text{and}  \,\,\, \phi_+(\sigma,\eta) \sim  \eta^{-1}
 \end{align}
 as $\eta\rightarrow 0$. Therefore, $\phi_-$ is characterized, up to scaling, as the unique solution to \eqref{bessel} that is vanishing (or even finite) at $\eta=0$.
    In the following proposition we identify $\phi_-$ in terms of $I_1$ and connect it to $e(\sigma, x/\sigma^2 )$ using 
    \begin{align} \label{smap}
  e(\sigma, x/\sigma^2 ) &= -i \sigma^{-\f12} [e^{\frac{\pi}{\sigma}} -1]^{-\f12} M_{\frac{i}{2\sigma},\f12}(2ix/\sigma) \\ 
  &= 2\sigma^{-\f 32} [e^{\frac{\pi}{\sigma}} -1]^{-\f12} x ( 1+ O(x/\sigma)) \,\,\,\, \text{as} \,\,\,\ x \to 0\nn \,,
\end{align}
which follows from \eqref{phiat0} and \cite[(13.14.14)]{NIST}.

\begin{pr} \label{prop:bessel} For any $ \delta\in(0,\frac{\pi}{2})$, there exists $c>0$ such that for all $\sigma \in [0,c)$, on $\eta\in[0,\delta] $, $e(\sigma,x/\sigma^2) $ has the form
\begin{align} \label{exp:bessel}
e(\sigma, x/\sigma^2) = \sqrt{2} \sigma^{-\f 12} [ e^{\frac {\pi}{\sigma}} -1]^{-\f 12} \Big( \frac{\eta}{ \eta^{\prime}} \Big)^{\f 12} I_1(\eta/\sigma) (1+a_-(\sigma,\eta))\,,
 \end{align}
where $I_1$ is the modified Bessel function of the first kind, and $ a_- $ is a smooth function satisfying the bounds 
\begin{align*}
&|a_-( \sigma,\eta)|\les \sigma ,\,\,\ \abs{\dot{a}_-(\sigma,\eta)}\les \sigma, \,\,\,\,\,\, |\ddot{a}_-( \sigma,\eta)|\les 1, \\ 
&|\partial_\sigma\{a_-(\sigma,\eta)\}|\les 1,\,\,\abs{\partial_{\sigma }\{\dot{a}_-(\sigma,\eta)\}}\les \sigma^{-1} ,\,\,|\partial_\sigma^2\{a_-(\sigma,\eta)\}|\les\sigma^{-2},\,\,|\partial_\sigma^2\{\dot{a}_-(\sigma,\eta)\}|\les\sigma^{-3}
\end{align*}
uniformly on $ [0,\delta] $.
\end{pr}
 \begin{rmk} We remark that one can see from the calculations in  Lemma~\ref{etaLem} that ${\tilde{V}_- =O((\eta - \pi/2)^{-3})}$ as $ \eta \to \pi/2$. This is the main reason why we cannot approximate  $e(\sigma,x/\sigma^2)$ by $I_1$ beyond the turning point $\eta=\frac{\pi}{2}$.
 \end{rmk}

\begin{proof}
In the variable $\rho= \sigma^{-1} \eta$, Equation \eqref{bessel} becomes the perturbed Bessel equation
\begin{align*}
	\partial^2_\rho\{ \omega(\rho)\} + \rho^{-1} \partial_\rho\{ \omega(\rho)\} - ( \rho^{-2} +1) \omega(\rho) = \sigma^2 \tilde{V}_-( \sigma \rho)  \omega(\rho)\,.
\end{align*}
A fundamental system for the homogeneous equation
\begin{align*}
	\partial^2_\rho\{ \omega(\rho)\} + \rho^{-1} \partial_\rho\{ \omega(\rho)\} - ( \rho^{-2} +1) \omega(\rho) =0
\end{align*}
is given by modified Bessel functions of first order $ I_1(\rho),K_1(\rho) $ so that by variation of parameters, the function 
\begin{align*}
	&\phi_-(\sigma,\rho)=I_1(\rho)+ \sigma^2 \int\limits_0^\rho  \frac{  [-I_1(\rho) K_1(u) +K_1(\rho) I_1(u)]   \tilde{V}_-( \sigma u )\phi_-(\sigma u)}{ W(K_1(u), I_1(u))} \:du 
	\end{align*}
solves  \eqref{bessel}  (provided that the integral on the right converges) and vanishes at $\rho=0$. Evaluating the Wronskian via \cite[(10.28.2)]{NIST}, plugging in the ansatz $ \phi_-(\sigma,\rho)=I_1(\rho)(1+ a_-(\sigma,\rho)) $, and noting that $I_1 (u)$ has no real zeros we obtain the equation for $ a_-(\sigma,\eta) $:
\begin{align}\label{a-}
a_-(\sigma, \eta) &=  \sigma^2\int\limits_0^\rho  u [K_1(u)I_1(u) - I_1^{-1} (\rho)K_1(\rho) I_1^2(u)]   \tilde{V}_-( \sigma u )( 1+ a_- (\sigma, \sigma u ))\:  du \\ 
& =: \sigma^2\int\limits_0^{\rho} M(\sigma,u) ( 1+ a_- (\sigma, \sigma u )) du \nn\,.
\end{align}
\par
  We will first prove that $a_-(\sigma, \eta)$ is well-defined and bounded by analyzing the leading term
\begin{align}
a_{-,0}(\sigma, \eta) := \sigma^2\int\limits_0^\rho  u [K_1(u) I_1(u)- I_1^{-1} (\rho)K_1(\rho) I^2_1(u)]   \tilde{V}_-( \sigma u ) du \,.
\end{align}
For this, we record the following bounds on $I_1$ and $K_1$:
\begin{align}\begin{split}\label{IKbounds}
	&\abs{\partial_{u}^j\{K_1(u)\}} \les u^{-1-j}\la u \ra ^{j+\frac{1}{2}} e^{-u},\,\,\,\,\,\ j=0,1,2,3... \\
	& \abs{\partial_u^j\{I_1(u)\}}  \sim u^{1-j}\la u \ra ^{j-\frac{3}{2}} e^{u}, \,\,\,\,\,\ j=0,1, \\
	&\abs{\partial_u^2\{I_1(u)\} } \les u^{3-j}\la u \ra ^{j-\frac{7}{2}} e^u \,\,\,\,\,\ j=2,3, \,,
\end{split}
\end{align}
which may easily be deduced from \cite[(10.30.1-2) and (10.40.1-2)]{NIST}. In particular, they imply that
\begin{align*}
	&| I_1(u) K_1(u)| \les  \la u \ra^{-1}\,,\\
	&|I_1^{-1} (u)K_1(u)|  \les u^{-2} \la u \ra^2e^{-2u}\,, \\
	&|I_1^2(u)| \les u^2 \la u \ra^{-3} e^{2u} \,.
\end{align*}
Therefore, since by Lemma \ref{etaLem} $ \tilde{V}_-(\sigma u) $ is bounded for $ u\leq \rho $, which is in turn less than $ \sigma^{-1}\delta $, we may write
\begin{align}\begin{split}\label{ltone}
	|a_{-,0}(\sigma, \eta)| &\les \sigma^2 \int\limits_0^{\rho} u\Span{u}^{-1}\: du +  \sigma^2  \rho^{-2} \Span{\rho}^{2} e^{-2 \rho}\int\limits_0^{\rho } u^3\Span{u}^{-3} e^{2u} du \\
	&\les \sigma^2 \la \rho \ra +\sigma^2\les \sigma\,.
	\end{split}
\end{align} 
Moreover, as $x \to 0$, we have $|a_{-,0}(\sigma, \eta)| \les \sigma^2 \rho^2 \les x$. These bounds may be easily extended to $a_-(\sigma, \eta)$ itself by a contraction argument. In particular, we think of \eqref{a-} as the fixed point equation
\begin{align*}
	a_-=T(1)+ T(a_-)
\end{align*}
for the linear operator $T$ given by  $ Ta=  \sigma^2\int\limits_0^\rho M(\sigma,u) a(u)  du$.  For $ \sigma $ small enough, our computations show that $ T $ is a contraction on  $ L^\infty_\eta ((0,\delta))  $ and moreover that $ T(1) $ lies in this space. This implies that the $L^{\infty}$-norm of the fixed point, given by $ a_- =\sum_{n=0}^\infty T^{n+1}(1)$, is bounded by the norm of the first term, which is $ O\left( \sigma \right)$.\par

Having established the existence and boundedness of $a_-(\sigma,\eta)$, we now turn to the bounds on its derivatives. We first treat the $\eta$-derivative. We have that
\begin{align}
\dot{a}_- (\sigma, \eta) = -\sigma \partial_{\rho} \{I_1^{-1} (\rho)K_1(\rho)\}\int_0^\rho u I^2_1(u)  \tilde{V}_-( \sigma u ) ( 1+ a_- (\sigma, \sigma u ))  du\,.
\end{align}
To estimate this integral, we use that by \eqref{IKbounds}
\begin{align}\label{derIK}
	  \abs{\partial^j_{\rho} \{I_1^{-1} (\rho)K_1(\rho)\}} \les \rho^{-2-j} \la \rho \ra^{2+j} e^{- 2\rho} \,\,\,\, j=1,2,3 \,,
\end{align}
so that
\begin{align*}
	\abs{\dot{a}_{-}(\sigma,\eta)} \les \sigma \rho^{-3}\Span{\rho}^{3} e^{-2\rho}\int_{0}^{\rho}u^3\Span{u}^{-3} e^{2u} \: du \les \sigma\,.
\end{align*}
Differentiating again, we find that
\begin{align}\label{adoubleeta}
	\ddot{a}_- (\sigma, \eta) &= -\partial_{\rho} \{I_1^{-1} (\rho)K_1(\rho)\}\rho I^2_1(\rho)   \tilde{V}_-( \eta ) ( 1+ a_- (\sigma, \eta ))\\
 &-\partial^2_{\rho} \{I_1^{-1} (\rho)K_1(\rho)\}\int\limits_0^\rho uI^2_1(u)  \tilde{V}_-( \sigma u ) ( 1+ a_- (\sigma, \sigma u ))  \:du\,.\nn
\end{align}
By \eqref{IKbounds}, the first term is uniformly bounded. For the second, we use \eqref{derIK} to argue similarly that the second term is uniformly bounded as well.\par
Now we estimate the $\sigma$-derivatives of $a(\sigma, \eta)$. For the first derivative, we have 
\begin{align} 
\partial_\sigma \{a_-(\sigma, \eta )\}&  =\sigma^{-1} (2  a_-(\sigma,\eta) -  \eta  \dot{a}(\sigma,\eta))\nn\\
& + \sigma^2\int\limits_0^\rho  u^2 [K_1(u) - I_1^{-1} (\rho)K_1(\rho) I_1(u)]   \tilde{V}_-'( \sigma u )I _1(u) ( 1+ a_- (\sigma, \sigma u ))\:  du\nn \\
& + \sigma^2\int\limits_0^\rho  u^2 [K_1(u) - I_1^{-1} (\rho)K_1(\rho) I_1(u)]   \tilde{V}_-( \sigma u )I _1(u)  \dot{a}_- (\sigma, \sigma u ))\:  du\nn  \\
&+ \sigma^2\int\limits_0^\rho  u [K_1(u) - I_1^{-1} (\rho)K_1(\rho) I_1(u)]   \tilde{V}_-( \sigma u ) I _1(u) \partial_\sigma\{ a_- (\sigma,x)\}|_{x=\sigma u} \:du  \nn \\
&=:A_1(\sigma,\eta)+A_2(\sigma,\eta)+A_3(\sigma,\eta)+A_4(\sigma,\eta)\nn\,.
\end{align}
By the bounds on $a(\sigma,\eta) $ and $\dot{a} (\sigma, \eta)$, it is clear $|A_1(\sigma,\eta)|\les 1$. Since $ \tilde{V}_- $ is uniformly smooth by Lemma \ref{etaLem}, it is easy to argue from the bound $ \abs{a_{-,0}} \les \sigma $ that $\abs{A_2(\sigma,\eta)}  \les 1 $, the only difference being an extra power of $ u $ in the integral. Similarly, $ \abs{A_3(\sigma,\eta)}\les \sigma  $ due to the previously derived bound $ \abs{\dot{a}(\sigma,\eta)}\les \sigma  $.\par
We have shown then that $ \partial_{\sigma}\{a(\sigma,\eta)\}  $ satisfies a fixed point equation of the form
\begin{align*}
	 \partial_{\sigma}\{a_-(\sigma,\eta)\}  =O(1)+ \sigma^2\int\limits_0^\rho  u [K_1(u) - I_1^{-1} (\rho)K_1(\rho) I_1(u)]   \tilde{V}_-( \sigma u ) I _1(u) \partial_\sigma\{ a_- (\sigma, \sigma u )\} \:du  
\end{align*}
and since we have already established, via \eqref{ltone}, that the last term is bounded in terms of $ \sigma\sup_{\eta\in [0,\delta)} \abs{\partial_\sigma\{ a_- (\sigma, \sigma u )\}} $, we may iterate this equation to find that, for $ \sigma $ sufficiently small, $ \partial_{\sigma} \{a(\sigma,\eta)\} $ is uniformly bounded independent of $ \sigma $ on the domain under consideration.\par
For the mixed $ \sigma $ and $ \eta $ derivative, we first compute that
\begin{align}\label{sigmadot}
	&\partial_{\sigma}\{\dot{a}(\sigma,\eta)\} =\sigma^{-1}\dot{a}(\sigma,\eta)-\sigma^{-1}\eta\ddot{a}_{-}(\sigma,\eta)\\
	&-\sigma \partial_{\rho} \{I_1^{-1}(\rho)K_1(\rho)\}\int\limits_{0}^{\rho} u^2I_1^2(u)\tilde{V}'(\sigma u)(1+a_-(\sigma,\sigma u))\: du \nn \\
	&-\sigma \partial_{\rho} \{I_1^{-1}(\rho)K_1(\rho)\}\int\limits_{0}^{\rho} u^2I_1^2(u)\tilde{V}(\sigma u)\dot{a}_-(\sigma,\sigma u)\: du \nn \\
	&-\sigma \partial_{\rho} \{I_1^{-1}(\rho)K_1(\rho)\}\int\limits_{0}^{\rho} uI_1^2(u)\tilde{V}(\sigma u)\partial_\sigma \{a_-(\sigma,x)\}|_{x=\sigma u} \: du\,. \nn
\end{align}
The first two terms are bounded by $ \sigma^{-1} $ and the third and fourth are bounded by a constant times
\begin{align*}
	\sigma\rho^{-3} \Span{\rho}^{3} e^{-2\rho}\int_{0}^{\rho} u^4\Span{u}^{-3} e^{2u} \: du \les \sigma\Span{\rho}^{2} \les \sigma^{-1}\,.
\end{align*}
Similarly, the last term is $ O(\sigma^{-1}) $ from the fact that $ \partial_{\sigma}\{a(\sigma,\eta)\}  $ is bounded.\par
The second $\sigma $-derivative is now estimated by differentiating each $A_i$ for $i=1,2,3,4$. It is easy to see from the bounds we have already developed that
\begin{align*}
	\abs{\partial_{\sigma}\{A_1(\sigma,\eta)\}} \les \sigma^{-2}\,,
\end{align*}
the dominant term being $ \sigma^{-1}\eta \partial_{\sigma}\{\dot{a}(\sigma,\eta) \}  $. Furthermore,
\begin{align*}
	\abs{\partial_{\sigma} \{A_2(\sigma,\eta)\}} &\les \sigma^{-1} \abs{A_2(\sigma,\eta)}+\Big|\eta\partial_{\sigma}\{I_1^{-1}(\rho)K_1(\rho)\} \int\limits_{0}^{\rho} u^2I_1^2(u)\: du \Big|\\
 	&+\sigma^2 \Big|\int_0^\rho  u^3 [-K_1(u)I_1(u) + I_1^{-1} (\rho)K_1(\rho) I_1^2(u)]   \:  du \Big|
\end{align*}
which is in total bounded in terms of  $ \sigma^{-1} $. By the same token,
\begin{align*}
	\abs{\partial_{\sigma}\{A_3(\sigma,\eta)\}} &\les \sigma^{-1}|A_3(\sigma,\eta)|+\sigma\Big|\eta\partial_{\sigma}\{I_1^{-1}(\rho)K_1(\rho)\} \int\limits_{0}^{\rho} u^2I_1^2(u)\: du \Big|\\
 	&+\sigma \Big|\int\limits_0^\rho  u^2\la u \ra [K_1(u)I_1(u) - I_1^{-1} (\rho)K_1(\rho) I_1^2(u)]   \:  du\Big|\,,
\end{align*}
where we have used that $ \abs{\partial_{\sigma}\{\dot{a}(\sigma,\eta)\}  }\les \sigma^{-1}$. As before, we may argue that all three terms are bounded by a constant times $ \sigma^{-1} $. One can also show that
\begin{align*}
	A_4(\sigma,\eta)=O(\sigma^{-1})+ \sigma^2\int\limits_0^\rho  u [K_1(u) - I_1^{-1} (\rho)K_1(\rho) I_1(u)]   \tilde{V}_-( \sigma u ) I _1(u) \partial_\sigma^2\{ a_- (\sigma, \sigma u )\} \:du\,,
\end{align*}
so that 
\begin{align*}
	\partial_\sigma^2\{ a_- (\sigma, \sigma u )\}=O(\sigma^{-2})+ \sigma^2\int\limits_0^\rho  u [K_1(u) - I_1^{-1} (\rho)K_1(\rho) I_1(u)]   \tilde{V}_-( \sigma u ) I _1(u) \partial_\sigma^2\{ a_- (\sigma, \sigma u )\} \:du\,.
\end{align*}
Since the last term is bounded in terms of $ \sigma\sup_{\eta\in [0,\delta)} \abs{\partial_\sigma\{ a_- (\sigma, \sigma u )\}} $, as before we can iterate for small enough $ \sigma $ to see that in fact $ \partial_{\sigma}^2\{a_-(\sigma,\eta)\} =O(\sigma^{-2})$, as claimed.\par

We next consider $\partial_\sigma^2\{\dot{a}_-(\sigma,\eta)\}$. Differentiating \eqref{sigmadot} with respect to $\sigma$ and incorporating various bounds we obtained, we find that 
$$
\partial^2_{\sigma}\{\dot{a}_-(\sigma,\eta)\}= - \sigma^{-1} \eta \partial_{\sigma}\{\ddot{a}_-(\sigma,\eta)\} + O (\sigma^{-2})
$$

To bound the first term we differential  \eqref{secondderer} in $\sigma$ and obtain
\begin{multline*}
	\partial_\sigma\{\ddot{a}_- (\sigma, \eta)\}=  - \sigma^{-1} \rho \partial_\rho\{ \ddot{a}_-(\sigma, \eta)\} - \partial_{\rho} \{I_1^{-1} (\rho)K_1(\rho)\}\rho I^2_1(\rho)   \tilde{V}_-( \eta ) \partial_\sigma\{ a_- (\sigma, \eta ))\} 
  \\ -\partial^2_{\rho} \{I_1^{-1} (\rho)K_1(\rho)\}  \int_0^{\rho} u I_1^2(u) \partial_\sigma \{ \tilde{V}_-( \sigma u ) ( 1+ a_- (\sigma, \sigma u )) \} \:du\,.
\end{multline*}
By \eqref{derIK}, we have $|\partial_\rho\{ \ddot{a}_-(\sigma, \eta)\}| \les \rho^{-1} \la \rho \ra $.  Furthermore, as $\partial_{\sigma} \{ a_-(\sigma,\eta)\} | \les 1 $ we have 
$$
|\partial_\sigma\{\ddot{a}_- (\sigma, \eta)\}| \les \rho^{-1} \la \rho \ra + \rho^{-4}\la \rho \ra^{4} e^{-2\rho} \int_0^{\rho} u^4 \la u \ra^{-3} e^{2u} du \les \rho^{-1} \la \rho \ra + \la \rho \ra 
$$
Therefore, $\partial^2_{\sigma}\{\dot{a}_-(\sigma,\eta)\}= O(\sigma^{-3}) $.

Finally, we match the solution $\phi_-(\sigma,\eta)$ to $e(\sigma,x/\sigma^2)$. For any fixed $\sigma$, as $x \to 0$ we have 
 \begin{align*}
 \Big( \frac{\eta(x)}{ \eta^{\prime}(x)} \Big)^{\f 12} I_1(\rho) = \frac{\sqrt{2} x} { \sigma} (1+ O( x/\sigma))
 \end{align*} 
 as a consequence of \cite[(10.30.1)]{NIST}. Comparing this expansion to \eqref{smap} we obtain \eqref{exp:bessel}.
\end{proof}

Before we focus on the other regimes of $x$, we give the following bounds which will be useful for the oscillatory estimates.  
\begin{lemma} \label{lem:bb}For any $\delta\in (0,1)$ there exists $c>0$ and $\varepsilon>0$ so that
\begin{align}\label{besselbound}
 |\partial^j_\sigma \{ e(\sigma, r)\}| \les e^{-\frac{\varepsilon}{\sigma}} r,\,\,\,\ j=0,1,2
\end{align}
uniformly in the support of $\chi_c(\sigma) \chi_{\delta}(\sigma^2 r)$. 
\end{lemma}
\begin{proof} We will use throughout the proof that $x=\sigma^2r\leq \delta<1$. This allows us to apply Proposition \eqref{prop:bessel} on the interval $[0,\eta(\delta)]$ to obtain the representation
\begin{align}
e(\sigma,r)= \sqrt{2} \sigma^{-\f 12} [ e^{\frac {\pi}{\sigma}} -1]^{-\f 12} \Big( \frac{\eta(\sigma^2r)}{ \eta^{\prime}(\sigma^2 r)} \Big)^{\f 12} I_1 \Big(\frac{\eta(\sigma^2 r)}{\sigma}\Big) (1+a_-(\sigma,\eta(\sigma^2 r))\,,
\end{align}

 By series expansion, 
\begin{align*}
\eta(\sigma^2 r) = 2 (\sigma^2 r) ^{\f12} -{\f 13} (\sigma^2 r)^{\f 32} + O_2( (\sigma^2 r)^{\f 52})    
\end{align*} 
which gives
\begin{align*}
\Big( \frac{\eta}{ \eta^{\prime}} \Big)^{\f 12}(\sigma^2r) = \sqrt{2} \sigma r^{\f 12} ( 1+ O_2(\sigma^2 r))\,.
\end{align*}
Moreover, using \eqref{IKbounds} we obtain
\begin{align*}
|I_1 \Big(\frac{\eta(\sigma^2 r)}{\sigma}\Big) | \les \min \Big\{ 1, \Big| \frac{\eta(\sigma^2 r)}{\sigma}\Big|  \Big\} e^{\frac{\eta(\sigma^2 r)}{\sigma}} \les e^{\frac{\pi}{2\sigma} - \frac{\varepsilon}{\sigma}} \min\{1,r^{\f12} \}
\end{align*}
for $\epsilon<\frac{\pi}{2}-\eta(\delta)$. Hence,   
\begin{align*}
|e(\sigma, \sigma r)|\les  \sigma^{-\f 12}[ e^{\frac {\pi}{\sigma}} -1]^{-\f 12}  \sigma r^{\f 12} e^{\frac{\pi}{2\sigma} - \frac{\varepsilon}{\sigma}} \min\{1,r^{\f12} \} \les e^{-\frac{\varepsilon}{\sigma}} r^{\f 12} \min\{1,r^{\f12} \} \les e^{-\frac{\varepsilon}{\sigma}} r.
\end{align*}
This establishes \eqref{besselbound} for $j=0$. 

We observe that the estimates for $j=1,2$ follow in a straightforward manner from the previously employed estimates. In particular, by \eqref{IKbounds} and the bounds on $a_-(\sigma,\eta)$ from Proposition \eqref{prop:bessel}, one can see that the impact of differentiation in $\sigma$ on $e(\sigma,r)$ is either gaining powers of $\sigma r$ or losing powers of $\sigma$ i.e. one obtains $|\partial_j\{e(\sigma,r)\}|\les e^{-\frac{\varepsilon}{\sigma}}\sigma^{-n}(\sigma r)^m$ for some $n,m\geq0$. With this in mind, within the support of $\chi_c(\sigma) \chi_{\delta}(\sigma^2 r)$, one can utilize the following estimate ($n, m\geq 0  $)
 $$ e^{-\frac{\varepsilon}{\sigma}} \sigma^{-n} (\sigma r)^{m}\les e^{-\frac{\varepsilon}{\sigma}} \sigma^{-n-m+1} r \les e^{-\frac{\varepsilon}{\sigma}} r
 $$
 to conclude the proof.

\end{proof}

\subsection{Airy function approximation: $ x\sim 1 $}
Let $ Q(u)=u^{-1}-1 $ and define the \emph{Liouville-Green transform}
\begin{align*}
\zeta(x)=\sign(x-1)\left|\frac{3}{2}\int\limits_{1}^{x} \sqrt{\abs{Q(u)} } \: du \right|^{\frac{2}{3}}	\,.
\end{align*}
Its properties are summarized in the following lemmas:
\begin{lemma}\label{LGLem}
	The map $ x\mapsto \zeta(x) $ defines a smooth change of variables from $ (0,\infty)\rightarrow (-(\frac{3\pi}{4})^{\frac{2}{3}},\infty) $. Furthermore, $ \zeta $ has the explicit form given by
\begin{align}\label{zetaForm}
	\frac{2}{3}\zeta^{\frac{3}{2}}(x)&=\sqrt{x(x-1)} -\log(\sqrt{x} +\sqrt{x-1} ),\quad x\geq 1,\\
 	-\frac{2}{3}(-\zeta(x))^{\frac{3}{2}}&=\sqrt{x(1-x)} -\frac{1}{2}\arccos(2x-1),\quad x\leq 1.\nonumber
\end{align}

	The function $ q=-\frac{Q}{\zeta} $ is non-negative and satisfies $ \sqrt{q}=\frac{d\zeta}{dx}  $. Under the transformation $ w(\zeta)=q^{\frac{1}{4}}f $, the equation \eqref{xODE} becomes 
	\begin{align}\label{zetaODE}
		-\sigma^2\ddot{w}(\zeta)=(\zeta+\sigma^2V)w(\zeta)
	\end{align}
where
\begin{align*}
	V=-q^{-\frac{1}{4}}\frac{d^2q^{\frac{1}{4}}}{d\zeta^2}\,.
\end{align*}
Here and for the rest of the paper we use \ $\dot{} = \frac{d}{d\zeta}$.  
\end{lemma}
\begin{proof}
The smoothness of $ \zeta $ is clear away from $ x=1 $ and at this point it is a simple consequence of the fact that $ Q $ vanishes only to first order. Indeed, we may expand $ \sqrt{\abs{Q} }  $ into a series in powers of $ \sqrt{\abs{x-1} }  $ and integrate term by term to find that
\begin{align*}
	\int_{1}^{x} \sqrt{\abs{Q(u)} } \: du=\frac{2}{3}(x-1)^{\frac{3}{2}}(1+O(x-1))  
\end{align*}
from which the claim is immediate. We omit the proof of \eqref{zetaForm} and \eqref{zetaODE} as it can be verified by differentiation.

\end{proof}
For reference, we record all of the notation relevant to the Liouville-Green transform:
\begin{align*}
\begin{split}
	&Q(u)=u^{-1}-1,\,q=-\frac{Q}{\zeta}\\
	&\zeta(x)=\sign(x-1)\left|\frac{3}{2}\int\limits_{1}^{x} \sqrt{\abs{Q(u)} } \: du \right|^{\frac{2}{3}}\\
	&w=q^{\frac{1}{4}}f,\,V=-q^{-\frac{1}{4}}\frac{d^2q^{\frac{1}{4}}}{d\zeta^2}
	\end{split}\,.
\end{align*}
\begin{lemma} \label{smallV}Let $\zeta^* = - \big( 3\pi/4)^{\f 23 }$ and $ \zeta \in  (\zeta^*, 0]$ then we have  $ |\partial^j_{\zeta} V| \les 1 $ for $j=0,1,2\dots$. 
\end{lemma}
\begin{proof} We note that as $|x-1| <1$, one has $ \zeta \sim (x-1)$, and therefore, $q^{\f14} = \sum_{k=0}^{\infty} c_k \zeta^k$ for some $c_k \in \R$. This shows that 
$ |\partial^j_{\zeta} V| \les 1 $ in the range of $|\zeta|<1$.  
On the other hand, as $|x|<1$, one has $ (\zeta - \zeta^*)^{\f32} \sim x$, therefore, $q^{\f14} \sim (\zeta - \zeta^*)^{-\f38}$ and $V (\zeta) \in O_{\infty}(  (\zeta - \zeta^*)^{-2})$. This shows that $ |\partial^j_{\zeta} V| \les 1 $ as long as $|\zeta - \zeta^* |> \delta >0$. 
\end{proof}

We now construct a basis of solutions to \eqref{zetaODE} in terms of the Airy functions $ \Ai $ and $ \Bi $ whose properties may be found in \cite{Olver}.
\begin{prop}\label{smallexp} Let $ \delta>0 $. Then there exists $c>0$ such that for all $\sigma \in [0,c]$, a fundamental system of solutions to \eqref{zetaODE} in the range  $ \zeta^*+\delta < \zeta \leq 0$ is given by 
\begin{align}
	\begin{split}\label{airyBasis}
		\phi_1(\sigma,\zeta) = \Ai(\tau)( 1+{ \sigma}  a_1( \sigma,\zeta) ) \\
\phi_2(\sigma,\zeta) = \Bi(\tau)( 1+  {\sigma}a_2(\sigma,\zeta))
	\end{split}
\end{align}
where $\tau:= - \sigma^{-\frac{2}{3}} \zeta$ and $ a_1 $ and $ a_2 $ are smooth functions satisfying the bounds
\begin{align}\label{ajbound}
\begin{split}
& | a_j(\sigma,\zeta) | \les 1 ,\,\,\,| \dot{a}_j(\sigma,\zeta)| \les \sigma^{-\frac{1}{3}} ,\,\,\,\,   | \ddot{a}_j(\sigma,\zeta)| \les \sigma^{-\f 43},\,\,\\
& \abs{\partial_{\sigma}\{a_j(\sigma,\zeta)\}} \les \sigma^{-\frac{4}{3}},\,\,\,\abs{\partial_{\sigma}\{\dot{a}_j(\sigma,\zeta)\} }\les \sigma^{-\frac{7}{3}},\,\,\,\abs{\partial_{\sigma}^{2}\{a_j(\sigma,\zeta)\} }\les \sigma^{-\frac{10}{3}}, \\ 
&\abs{\partial^2_{\sigma}\{\dot{a}_j(\sigma,\zeta)\} }\les \sigma^{-\frac{13}{3}} \,.
\end{split}
\end{align}
for $j=1,2$ uniformly on $ [\zeta_*+\delta,0] $.
\end{prop}
\begin{rmk} The range of $ \zeta $ corresponds to $ x\in [\delta',1] $ for some $ \delta'>0 $ independent of $ \sigma $. The restriction is designed to avoid the singularity of $ V$ at $ \zeta=\zeta_* $, see the proof of Lemma~\ref{smallV}. Note also that this approximation is only possible because $ \tau>0 $ for $ \zeta<0 $ and thus the Airy functions do not have zeroes in this regime.
\end{rmk}
\begin{proof}
	Write $ \phi_{1,0}(\sigma,\zeta)=\Ai(\tau) $ and $ \phi_{2,0}(\sigma,\zeta)=\Bi(\tau) $. The variable $ \tau $ is chosen so that
\begin{align*}
	-\sigma^2\ddot{\phi}_{j,0}-\zeta \phi_{j,0}=0
\end{align*}
for each of $ j=1,2 $ where $ \dot{ }=\frac{\partial}{\partial\zeta} $. Therefore,
\begin{align*}
		-\sigma^2\ddot{\phi}_{j}-\zeta\phi_{j}=\left( \sigma^{3}\phi_{j,0}^2\dot{a}_j \right)^./\phi_{j,0}
\end{align*}
	and plugging the representations \eqref{airyBasis} into \eqref{zetaODE} yield the equation for $ a_j $
\begin{align}\label{ajEq}
		\left( \phi_{j,0}^2\dot{a}_j \right)^.=-\sigma^{-1}V\phi_{j,0}^2(1+\sigma a_j)\,.
\end{align}
	The solution to this equation for $ j=2 $ with $ a_2(\sigma,0)=0 $ and $ \dot{a}_2(\sigma,0)=0 $ is given by
\begin{align} \label{a2full}
	a_2(\sigma,\zeta) =  - \sigma^{ \frac{1}{3}}\int\limits_{0}^{ - \sigma^{-\frac{2}{3}} \zeta} \Bi^2(u) \Big[\int\limits_u^{ - \sigma^{-\frac{2}{3}} \zeta} \Bi^{-2} (v) \:dv \Big]  V(-\sigma^{\frac{2}{3}} u) ( 1+ \sigma a_2(\sigma,-\sigma^{\frac{2}{3}} u)) \:du \:. 
\end{align}

We now recall the following expansions of the Airy functions found in \cite{Olver}:
\begin{align} \label{asymptoticAB}
&\Bi(x) = \pi^{-\f1 2} x^{-\f 14} e^{{\f 23}x^{\f 32}} ( 1+ O (x^{-\f 32})) \,\, \text{as}\,\, x \to \infty  \\ 
&\Bi(x) \geq \Bi (0) > 0 , \,\,\,\text{ for $x\geq 0$} \nn \\ 
& \Ai(x) = \frac{1}{2}\pi^{-\f 1 2} x^{-\f 14} e^{-{\f 23}x^{\f 32}} ( 1+ O (x^{-\f 32})) \,\, \text{as}\,\,\,  x \to \infty \\ 
& \Ai(x) > 0   \,\,\,\text{ for $x\geq 0$}\,. \nn 
\end{align}
These asymptotics and the identity
\begin{align}\label{Bil}
   \int\limits_{x_0}^{x_1} \Bi^{-2}(y)\:dy=\pi^{-1}\left(\frac{\Ai}{\Bi}(x_0)-\frac{\Ai}{\Bi}(x_1)\right)\,,
\end{align}
which holds for $0\leq x_0<x_1$,
imply that for $x_0\geq 0$
\begin{align*}
    \left| \Bi^2(x_0)\int\limits_{x_0}^{x_1}\Bi^{-2}(y)\:dy\right|\les \Span{x_0}^{-\frac{1}{2}}\,,
\end{align*}
so that 
\begin{align}\label{Bil_2}
	\left|\int_{0}^{ - \sigma^{-\frac{2}{3}} \zeta} \Bi^2(u) \Big[\int\limits_u^{ - \sigma^{-\frac{2}{3}} \zeta} \Bi^{-2} (v) \:dv \Big]  f(u)  \:du  \right| \les \Span{\sigma^{-\frac{2}{3}}\zeta}^{\frac{1}{2}} \|f\|_{\infty}\,.
\end{align}
Moreover, 
\begin{align}\label{Bil2}
    \left|\Bi^{-2}(x_0)\int\limits_0^{x_0}\Bi^{2}(y)\:dy\right|\les \Span{x_0}^{-\frac{1}{2}}
\end{align}
which comes from inserting the above asymptotics into the integral and then computing for $x>1$
\begin{align*}
    x^{\frac{1}{2}}e^{-\frac{4}{3}x^\frac{3}{2}}\int\limits_1^{x}y^{-\frac{1}{2}}e^{\frac{4}{3}y^\frac{3}{2}}\:dy&\les x^{\frac{1}{2}}e^{-\frac{4}{3}x^\frac{3}{2}}\int\limits_1^{x^\frac{3}{2}}u^{-\frac{2}{3}}e^{\frac{4}{3}u}\:du= x^{\frac{1}{2}}e^{-\frac{4}{3}x^\frac{3}{2}}\left(\frac{3}{4}u^{-\frac{2}{3}}e^{\frac{4}{3}u}\bigg\vert_1^{x^\frac{3}{2}}+\frac{2}{3}\int\limits_1^{x^\frac{3}{2}}u^{-\frac{5}{3}}e^{\frac{4}{3}u}\:du\right)\\
    &\les x^{\frac{1}{2}}e^{-\frac{4}{3}x^\frac{3}{2}}\left(x^{-1}e^{\frac{4}{3}x^\frac{3}{2}}\right)=x^{-\frac{1}{2}}\,.
\end{align*}

To estimate \eqref{a2full}, we let 
\begin{align*}
	a_{2,0}(\zeta) := - \sigma^{ \f 13 }\int\limits_{0}^{ - \sigma^{-\frac{2}{3}} \zeta} \Bi^2(u) \Big[\int\limits_u^{ - \sigma^{-\frac{2}{3}} \zeta} \Bi^{-2} (v) \:dv \Big]  V(-\sigma^{\frac{2}{3}}u)  \:du
\end{align*}
be the leading term, where we have suppressed the $ \sigma $-dependence of $ a_2 $ for now. By Lemma~\ref{smallV} $ V(-\sigma^{\frac{2}{3}}u) $ is bounded on the domain of integration when $ \zeta\in [\zeta_*+\delta,0] $ so using \eqref{Bil_2} we have that
\begin{align} \label{a20est}
	|a_{2,0}(\zeta)|\les  \sigma^{\f 13} \Span{\sigma^{-\frac{2}{3}}\zeta}^{\frac{1}{2}}  \les 1\,.
\end{align}
 Now, a contraction argument shows that $|a_2(\zeta)| \les 1 $, as claimed.

We next consider the $\zeta$ derivative of $a_2(\zeta)$. One has that
\begin{align}\label{adot}
\dot a_2(\zeta) = \sigma^{-\frac{1}{3}} \Bi^{-2}(-\sigma^{-\frac{2}{3}} \zeta)\int\limits_0^{-\sigma^{-\frac{2}{3}} \zeta} \Bi^2(u) V(-\sigma^{\frac{2}{3}} u) ( 1+ \sigma a_2 (-\sigma^{\frac{2}{3}} u)) \:du
\end{align}
so that, since $ V $ is bounded and we have shown that $ a_2 $ itself is bounded, we see that by \eqref{Bil2}
\begin{align*} 
	\abs{\dot{a}_2(\zeta)} \les \sigma^{-\frac{1}{3}}\Span{\sigma^{-\frac{2}{3}}\zeta}^{-\frac{1}{2}}\,,
\end{align*}
 which is less than $ \sigma^{-\frac{1}{3}} $, as claimed. 
\par
For the second $ \zeta $-derivative, we use \eqref{ajEq} to write
\begin{align*}
	\ddot{a}_2(\zeta)=-\sigma^{-1}V(\zeta) (1+\sigma a_2(\zeta))-2\sigma^{-1}[\dot{\phi}_{2,0}\phi_{2,0}^{-1}](\zeta)\dot{a}_2(\zeta)\,.
\end{align*}
The first term is clearly bounded by $ \sigma^{-1} $ while the second is bounded in terms of
\begin{align}\label{BprimeB}
	\sigma^{-1}\abs{\Bi'(-\sigma^{-\frac{2}{3}}\zeta)\Bi^{-1}(-\sigma^{-\frac{2}{3}}\zeta)} \les \sigma^{-\frac{4}{3}}.
\end{align}
By \cite[(9.7.8)]{NIST}, $ \abs{\Bi'(-\sigma^{-\frac{2}{3}}\zeta)}\les \Span{\sigma^{-\frac{2}{3}}\zeta}^{\frac{1}{4}} e^{\frac{2}{3}\sigma^{-1}\abs{\zeta} ^{\frac{3}{2}}}$, which shows that $ \abs{\ddot{a}_2(\zeta)} \les \sigma^{-\frac{4}{3}} $.\par
Now, we consider the $ \sigma $-derivatives of $ a_2 $. Let $F(\sigma,u) = V(-\sigma^{\f23} u) ( 1+\sigma a_2(\sigma,-\sigma^{\f23}u))$, then by \eqref{a2full} we can compute
\begin{align*}
	\partial_{\sigma} \{a_2(\sigma,\zeta)\}  & =-\frac{1}{3 \sigma } [ a_2(\sigma,\zeta)-2\zeta \dot{a_2}(\sigma,\zeta) ] \\
	&-\sigma^{\f13} \int\limits_{0}^{-\sigma^{-\frac{2}{3}}\zeta} \Bi^2(u)\Big[\int\limits_{u}^{-\sigma^{-\frac{2}{3}}\zeta} \Bi^{-2}(v)\: dv\Big] \partial_{\sigma}\{F(u,\sigma)\}\: d\sigma 
 =:B_1(\sigma,\zeta)+B_2(\sigma,\zeta)\,.
\end{align*}
Using the bounds previously established for $a_2(\sigma,\zeta)$, we can deduce $|B_1(\sigma,\zeta)|\les \sigma^{-\frac{4}{3}} $. Moreover, we may write 
 \begin{align}
 \partial_{\sigma}\{F(u,\sigma)\}= O(\sigma^{-\f13} \la u \ra )+ \sigma V(-\sigma^{\f32}u)\partial_\sigma[a_2(\sigma,v)]_{v=-\sigma^{\f32}u}. 
 \end{align}
 Therefore, by \eqref{Bil_2} we have
\begin{align*}
B_2(\sigma,\zeta)= \sigma^{\f13} \la \sigma^{-\f23} \zeta \ra^{\f32}+ -\sigma^{\f43} \int\limits_{0}^{-\sigma^{-\frac{2}{3}}\zeta} \Bi^2(u)\Big[\int\limits_{u}^{-\sigma^{-\frac{2}{3}}\zeta} \Bi^{-2}(v)\: dv\Big] V(-\sigma^{\f32}u)\partial_\sigma[a_2(\sigma,v)]_{v=-\sigma^{\f32}u }\: du \,.
\end{align*}
Letting 
$$
T(a) := \sigma^{\f13} \int\limits_{0}^{-\sigma^{-\frac{2}{3}}\zeta} \Bi^2(u)\Big[\int\limits_{u}^{-\sigma^{-\frac{2}{3}}\zeta} \Bi^{-2}(v)\: dv\Big] V(-\sigma^{\f32}u)a(\sigma,u) \: du\,,
$$
we obtain 
\begin{align}
\partial_{\sigma} \{a_2(\sigma,\zeta)\} = \sigma^{-1} + \sigma T\Big( \partial_\sigma[a_2(\sigma,v)]_{v=-\sigma^{\f32}u }\Big)\,.
\end{align}
Now, by a contraction argument we obtain that $ \abs{\partial_{\sigma}\{a_2(\sigma,\zeta)\} }\les \sigma^{-\frac{4}{3}} $.\par
Proceeding onward, we differentiate \eqref{adot} with respect to $ \sigma $ to find that
\begin{align} \label{dersigmazeta}
	\partial_{\sigma}\{\dot{a}_2(\sigma,\zeta) \}&=-\frac{1}{3}\sigma^{-1}\dot{a}_2(\sigma,\zeta)+\frac{2}{3}\sigma^{-\frac{5}{3}}\zeta\Bi'(-\sigma^{-\frac{2}{3}}\zeta)\Bi^{-1}(-\sigma^{-\frac{2}{3}}\zeta)\dot{a}_2(\sigma,\zeta)\\
	&+\frac{2}{3}\sigma^{-2}\zeta V_+(\zeta)(1+\sigma a_2(\sigma,\zeta)) \nn  \\
	&+\sigma^{-\frac{1}{3}} \Bi^{-2}(-\sigma^{-\frac{2}{3}} \zeta)\int\limits_0^{-\sigma^{-\frac{2}{3}} \zeta} \Bi^2(u) \partial_{\sigma} \{F(u,\sigma)\} \:du\,. \nn
\end{align}
It is easy to see using previously derived inequalities that each term is at least $ O(\sigma^{-2}) $ except for the second term. For that, we write
\begin{align*}
	\abs{\sigma^{-\frac{5}{3}}\zeta\Bi'(-\sigma^{-\frac{2}{3}}\zeta)\Bi^{-1}(-\sigma^{-\frac{2}{3}}\zeta)\dot{a}_2(\sigma,\zeta)} \les \sigma^{-2}\zeta \Span{\sigma^{-\frac{2}{3}}\zeta}^{\frac{1}{2}} \abs{\dot{a}_2(\sigma,\zeta)} 
\end{align*}
and recall that $ \abs{\dot{a}_2(\sigma,\zeta) }\les \sigma^{-\frac{1}{3}}\Span{\sigma^{-\frac{2}{3}}\zeta}^{-\frac{1}{2}}  $. It follows then that $ \abs{\partial_{\sigma}[\dot{a}_2(\sigma,\zeta)\} }\les \sigma^{-2} $.\par
For the second $ \sigma $-derivative of $ a_2 $, we differentiate each of $B_j(\sigma,\zeta)$ $j=1,2$ separately. First, it is easy to see that
\begin{align*}
	\abs{\partial_{\sigma} \{B_1(\sigma,\zeta)\}  } \les \sigma^{-1}\abs{B_1(\sigma,\zeta)}+ \sigma^{-1}(\abs{ \partial_\sigma\{a_2(\sigma,\zeta)\}} + \sigma^{-1} \abs{\partial_{\sigma}\{\dot{a_2}(\sigma,\zeta)\}}
\end{align*}
which is in total $ O(\sigma^{-3}) $, the dominant term being the last term. Next, differentiating $B_2(\sigma,\zeta)$ we have
\begin{align}\label{B2der}
	\partial_{\sigma} \{B_2(\sigma,\zeta)\}  &= \frac{1}{3\sigma}B_2(\sigma,\zeta) +\frac{2}{3} \sigma^{-\frac{4}{3}} \zeta \Bi^{-2}(-\sigma^{-\frac{2}{3}} \zeta)\int\limits_0^{-\sigma^{-\frac{2}{3}} \zeta} \Bi^2(u) \partial_\sigma\{F(u,\sigma)\}\:du\\+&\sigma^{\f13}\int_{0}^{\tau}\Bi^2(u)\left[\int_{u}^{\tau}\Bi^{-2}(v) \: dv\right] \partial_\sigma^2\{F(u,\sigma)\}  \: du \: \nn.
\end{align}
Similar calculations to those employed in the estimation of $B_2(\sigma,\zeta)$ itself demonstrate that the initial two terms on the right-hand side of the equation in \eqref{B2der} are $ O(\sigma^{-\frac{7}{3}}) $. To estimate the last term, we compute
\begin{align*}
   \partial_\sigma^2\{F(u,\sigma)\} = O( \sigma^{-1} \la u \ra^{2} + \sigma^{2}) + \sigma V(-\sigma^{\f23}u)\partial_\sigma^2[a_2(\sigma,v)]_{v=-\sigma^{\f23}u} \:.
\end{align*}
Hence, using \eqref{Bil2}, we deduce the following expression:
\begin{align*}
   \partial_{\sigma} \{B_2(\sigma,\zeta)\} = O(\sigma^{-\f73}) + \sigma T\Big( \partial_\sigma^2[a_2(\sigma,v]_{v=-\sigma^{\f23}u} \Big).
\end{align*}
This, in turn, leads to:
\begin{align*}
    \partial_\sigma^2[a_2(\sigma,\zeta)]= O(\sigma^{-3})+ \sigma T\Big( \partial_\sigma^2[a_2(\sigma,v]_{v=-\sigma^{\f23}u} \Big)\,,
\end{align*}
which by means of a contraction argument, yields $|\partial_\sigma^2[a_2(\sigma,\zeta)]|\les \sigma^{-3}$.

We next focus on $\partial^2_{\sigma}\{\dot{a}_2(\sigma,\zeta)\}$. Differentiating \eqref{dersigmazeta}, and using the previously obtained bounds, we have 
\begin{align*}
	\partial^2_{\sigma}\{\dot{a}_2(\sigma,\zeta) \}&= \frac{2}{3}\sigma^{-\frac{5}{3}}\zeta \partial_\sigma\{\Bi'(-\sigma^{-\frac{2}{3}}\zeta)\Bi^{-1}(-\sigma^{-\frac{2}{3}}\zeta)\}\dot{a}_2(\sigma,\zeta)\\
	&+\sigma^{-\frac{1}{3}} \partial_\sigma\{ \Bi^{-2}(-\sigma^{-\frac{2}{3}} \zeta)\}\int\limits_0^{-\sigma^{-\frac{2}{3}} \zeta} \Bi^2(u) \partial_{\sigma} \{F(u,\sigma)\} \:du\,+ O(\sigma^{-3})
\end{align*}

By \eqref{BprimeB} and the last term in \eqref{dersigmazeta}, the second term is bounded by $ \sigma^{-\frac{12}{13}}$. Moreover, one has that 
\begin{align*}
|\partial_\sigma\{\Bi'(-\sigma^{-\frac{2}{3}}\zeta)\Bi^{-1}(-\sigma^{-\frac{2}{3}}\zeta)\}| \les \sigma^{-\f73}  \,.
\end{align*}
Therefore, we have 
$\partial^2_{\sigma}\{\dot{a}_2(\sigma,\zeta) \}=O(\sigma^{-\frac{13}{3}})$.

Having proven all of the stated bounds on $ a_2 $, we now turn to $ a_1 $. We make the reduction ansatz $\phi_1(\zeta) =g(\zeta) \phi_2(\zeta)$ and find that $ g $ solves $ (\phi_2^2 \dot{g}\dot{)}=0 $. To simplify the analysis that follows, we extend the functions $ \phi_1 $ and $ \phi_2 $, defined at the moment on $ [\zeta_*+\delta,0] $, to the interval $ (-\infty,0] $ in such a way that the previous bounds still hold. We then choose the solution $ g $ of the form 
\begin{align*}
 g(\zeta) &=  \pi\int\limits_{\tau}^{ \infty} \Bi^{-2}(u) ( 1+ a_2(-\sigma^{\f 23} u))^{-2}  \:du \,,
\end{align*}
 which yields 
 \begin{align} \label{phi2}
	 \phi_1(\zeta ) =  \pi \Bi(\tau)  ( 1+\sigma a_2(\zeta))\int\limits_{\tau}^{\infty}  \Bi^{-2}(u) ( 1+  \sigma a_2(-\sigma^{\f 23}u))^{-2}  \:du \,.
 \end{align}
 Here, we suppressed the $\sigma$ dependence of $a_2$. We now write $ (1+\sigma a_2)^{-2}=1+\sigma \tilde{a}_2 $ and for $ \sigma $ sufficiently small, $ \tilde{a}_2 $ satisfies all of the same bounds as $ a_2 $ because $ \abs{a_2} \les 1 $. Recalling the identities  
 \begin{align} \label{ftc}
	\frac{d}{dx}\Big\{\frac{\Ai}{\Bi}(x)\Big\} =-\pi^{-1}\Bi^{-2}(x),\,\,\frac{d}{dx}\Big\{\frac{\Bi}{\Ai}(x)\Big\} =\pi^{-1}\Ai^{-2}(x)
 \end{align}
and the fact that $ \Ai(u) $ and $ \Bi(u) $ are strictly positive for $ u\geq 0 $,  we integrate by parts to see that
 \begin{multline}
	 \pi\int\limits_{\tau}^{ \infty}  \Bi^{-2}(u) ( 1+\sigma \tilde{a}_2( -\sigma^{\f23} u) ) \:du=\\
	 \left[\frac{\Ai}{\Bi}(u)(1+\sigma\tilde{a}_2(-\sigma^{\frac{2}{3}}u))\right]\bigg\vert_{\infty}^{\tau}
- \sigma^{\frac{5}{3}}\int\limits_{\tau}^{\infty} \frac{\Ai}{ \Bi}(u)   \dot{\tilde{a}}_2( -\sigma^{\f23} u) \:du\,.
 \end{multline}
Therefore, 
\begin{align*}
	\phi_1(\zeta)=\Ai(\tau){(1+\sigma a_2(\zeta))}\Big[ ( 1+ \sigma \tilde{a}_2(\zeta))	- \sigma^{\f 53 } \frac{\Bi}{ \Ai }(\tau)  \int\limits_{\tau}^{ \infty}\frac{\Ai}{ \Bi}(u)   \dot{\tilde{a}}_2( -\sigma^{\f23} u) \:du \Big]\,.
\end{align*}
From this, we infer that 
\begin{align*}
	a_1(\zeta)&=a_2(\zeta)+(1+\sigma a_2(\zeta))\Big[ \tilde{a}_2(\zeta)+  \sigma^{\f 23 } \frac{\Bi}{ \Ai }(\tau)  \int\limits_{-\sigma^{-\f23} \zeta}^{\infty}\frac{\Ai}{ \Bi}(u)   \dot{\tilde{a}}_2( -\sigma^{\f23} u) \:du \Big]\\
	&:=a_2(\zeta)+(1+\sigma a_2(\zeta))\left[\tilde{a}_2(\zeta)+\tilde{a}_1(\zeta)\right]\,,
\end{align*}
so it suffices to control  
\begin{align*}
	\tilde{a}_1(\zeta) =\sigma^{\frac{2}{3}}\frac{\Bi}{\Ai}(\tau)\int_{\tau}^{ \infty}\frac{\Ai}{ \Bi}(u)   \dot{\tilde{a}}_2( -\sigma^{\f23} u) \:du \,.
\end{align*}
\par 
To begin with, we use that $ \abs{\dot{\tilde{a}}_2(\zeta) }\les \sigma^{-\frac{1}{3}} $ to write
\begin{align*}
	\abs{\tilde{a}_1(\zeta)} &\les \sigma^{\frac{1}{3}}e^{\frac{4}{3}\Span{\sigma^{-\frac{2}{3}}\zeta}^{\frac{3}{2}}}\int_{-\sigma^{-\frac{2}{3}}\zeta}^{\infty} e^{-\frac{4}{3}\Span{u}^{\frac{3}{2}} }\: du \\ 
 &\leq\sigma^{\frac{1}{3}}e^{\frac{4}{3}\Span{\sigma^{-\frac{2}{3}}\zeta}^{\frac{3}{2}} }\Span{\sigma^{-\frac{2}{3}}\zeta}^{-\frac{1}{2}}\int\limits_{-\sigma^{-\frac{2}{3}}\zeta}^{\infty} e^{-\frac{4}{3}\Span{u}^{\frac{3}{2}}}\Span{u}^{\frac{1}{2}} \: du \\
		&\les \sigma^{\frac{1}{3}}\Span{\sigma^{-\frac{2}{3}}\zeta}^{-\frac{1}{2}}  \,.
\end{align*}

We compute further that
\begin{align}\label{BZetaDer}
\dot{\tilde{a}}_1(\zeta) =-\pi^{-1}\sigma^{-\frac{2}{3}}(\Ai\Bi)^{-1}(\tau)\tilde{a}_1(\zeta)+\dot{\tilde{a}}_2(\zeta)\,.
\end{align}
We have already bounded the second term by $ \sigma^{-\frac{1}{3}} $, while the first term obeys this bound as well since $ \abs{(\Ai\Bi)^{-1}(x)} \les \Span{x}^{\frac{1}{2}}  $. Proceeding onward,
\begin{align*}
\partial_{\zeta}^2\{\tilde{a}_1(\zeta)\} =\pi^{-1}\sigma^{-\frac{4}{3}}\frac{d}{du}\left[(\Ai\Bi)^{-1}(u)\right]\vert_{u=\tau}\tilde{a}_1(\zeta)-\pi^{-1}\sigma^{-\frac{2}{3}}(\Ai\Bi)^{-1}(\tau) \dot{\tilde{a}}_1 (\zeta)+\ddot{\tilde{a}}_2(\zeta)\,.
\end{align*}
From the fact that $ \abs{\frac{d}{du}\{(\Ai\Bi)^{-1}(u)\} } \les \Span{u}^{-\frac{1}{2}}  $, it is easily checked as before that each term is bounded by at worst $ \sigma^{-\frac{4}{3}} $, thus establishing all of the desired bounds on the $ \zeta $-derivatives of $ \tilde{a}_1 $.\par
For the $ \sigma $-derivatives, it is convenient to first rewrite
\begin{align*}
	\tilde{a}_1(\sigma,\zeta)=\frac{\Bi}{\Ai}(\tau)\int_{-\zeta}^{\infty} \frac{\Ai}{\Bi}(\sigma^{-\frac{2}{3}}v)\dot{\tilde{a}}_2(\sigma,-v)\: dv \,,
\end{align*}
so that
\begin{align*}
	\partial_{\sigma} \{\tilde{a}_1(\sigma,\zeta)\} &=\frac{2}{3\pi}\sigma^{-\frac{5}{3}}\zeta(\Ai\Bi)^{-1}(\tau)\tilde{a}_1(\sigma,\zeta)\\
 &+\frac{2}{3\pi}\sigma^{-\frac{5}{3}}\frac{\Bi}{\Ai}(\tau)\int_{-\zeta}^{\infty} \Bi^{-2}(\sigma^{-\frac{2}{3}}v)v\dot{\tilde{a}}_2(\sigma,-v)\: dv\\
	&+\frac{\Bi}{\Ai}(\tau)\int_{-\zeta}^{\infty} \frac{\Ai}{\Bi}(\sigma^{-\frac{2}{3}}v)\partial_{\sigma} \{\dot{\tilde{a}}_2\}(\sigma,-v) \: dv\\ 
 &=: D_1(\sigma,\zeta)+ D_2(\sigma,\zeta)+ D_3(\sigma,\zeta)\,.
\end{align*}
Arguing as before, it is easy to see that $ \abs{D_1(\sigma,\zeta)}\les \sigma^{-\frac{4}{3}} $ whereas 
\begin{align}\label{D2est}
	\abs{D_2(\sigma,\zeta)} \les \sigma^{-\frac{2}{3}}\frac{\Bi}{\Ai}(\tau)\int\limits_{\tau}^{\infty}\Bi^{-2}(v)v \: dv =\sigma^{-\frac{2}{3}}\frac{\Bi}{\Ai}(\tau)\left[\pi^{-1}\frac{\Ai}{\Bi}(\tau)\tau-\int\limits_\tau^\infty \frac{\Ai}{\Bi}(u)\,du\right]\les \sigma^{-\frac{4}{3}}\,,
\end{align}
where the second term in brackets may be treated as the original estimate of $ \tilde{a}_1 $. By using that $ \abs{\partial_{\sigma} \{\dot{\tilde{a}}_2(\sigma,\zeta)\} }\les \sigma^{-2}  $, we may similarly argue that $\abs{D_3(\sigma,\zeta)} \les \sigma^{-\frac{4}{3}} $. Thus, we conclude that $ \abs{\partial_{\sigma} \{\tilde{a}_1(\sigma,\zeta)\} } \les \sigma^{-\frac{4}{3}} $.\par
For the mixed derivative $ \partial_{\sigma} \{\dot{\tilde{a}}_1(\sigma,\zeta)\}  $, we differentiate \eqref{BZetaDer} to find that
\begin{align}\label{dersigmazetatilde}
\partial_{\sigma} \{\dot{\tilde{a}}_1(\sigma,\zeta)\} & =\frac{2}{3\pi}\sigma^{-\frac{5}{3}}(\Ai\Bi)^{-1}(\tau)\tilde{a}_1(\sigma,\zeta) +\frac{2}{3\pi}\sigma^{-\frac{7}{3}}\zeta\frac{d}{du}\{(\Ai\Bi)^{-1}(u)\}\vert_{u=\tau}\tilde{a}_1(\sigma,\zeta)  \\
&-\pi^{-1}\sigma^{-\frac{2}{3}}(\Ai\Bi)^{-1}(\tau)\partial_{\sigma} \{\tilde{a}_1(\sigma,\zeta)\} 
	+\partial_{\sigma} \{\dot{\tilde{a}}_2(\sigma,\zeta)\} \nn \,,
\end{align}
and it is merely a matter of collecting previously derived bounds to deduce that each term is bounded by $ \sigma^{-2} $, except the third term which is bounded by  $\sigma^{-7/3}$.

Finally, to bound the second $ \sigma $-derivative of $ \tilde{a}_2 $, we comment on the derivatives of each $ D_i(\sigma,\zeta)$, $i=1,2,3$. The expression $ \partial_{\sigma}\{D_1(\sigma,\zeta)\}  $ is essentially the same as $ \partial_{\sigma} \{\dot{\tilde{a}}_2(\sigma,\zeta)\}  $ with the loss of an additional $ \sigma $ power, so it is bounded in terms of $ \sigma^{-3} $. For $ D_2(\sigma,\zeta) $, one differentiates to find that
\begin{align*}
	\abs{\partial_{\sigma} \{D_2(\sigma,\zeta) \} } \les \sigma^{-1}\abs{D_2(\sigma,\zeta) } +\sigma^{-\frac{5}{3}}(\Ai\Bi)^{-1}(\tau)\abs{D_2(\sigma,\zeta) } +\\
	\sigma^{-\frac{10}{3}}\frac{\Bi}{\Ai}(\tau)\int_{-\zeta}^{\infty} \frac{d}{du}\{\Bi^{-2}(u)\}\vert_{\sigma^{-\frac{2}{3}}v} v^2\: dv+\sigma^{-4}\frac{\Bi}{\Ai}(\tau)\int_{-\zeta}^{\infty} \Bi^{-2}(\sigma^{-\frac{2}{3}}v)v \: dv \,,
\end{align*}
where in the second integral we have used the bound $ \abs{\dot{\tilde{a}}_2(\sigma,\zeta)}\les \sigma^{-2} $.
Collecting bounds easily shows that the first term is bounded by $ \sigma^{-\frac{7}{3}} $ and the second by $ \sigma^{-\frac{10}{3}} $. The first integral is bounded in terms of
\begin{align*}
	\sigma^{-\frac{4}{3}}\frac{\Bi}{\Ai}(\tau)\int\limits_{\tau}^{\infty} \frac{d}{du}\{\Bi^{-2}(u)\}u^2 \: du =\sigma^{-\frac{4}{3}}\frac{\Bi}{\Ai}(\tau)\left[\Bi^{-2}(\tau)\tau^2-2\int\limits_\tau^\infty\Bi^{-2}(u)u\:du\right]\les \sigma^{-\frac{8}{3}}
	\end{align*}
where we may bound the second term as in \eqref{D2est}. Similarly, the second integral is bounded by $ \sigma^{-10/3} $.\par
We now compute
\begin{align*}
	\partial_{\sigma}\{D_3(\sigma,\zeta)\} =\frac{2}{3\pi}\sigma^{-\frac{5}{3}}\zeta(\Ai\Bi)^{-1}(\tau)D_3(\sigma,\zeta)+\frac{2}{3\pi}\sigma^{-\frac{5}{3}}\frac{\Bi}{\Ai}(\tau)\int_{-\zeta}^{\infty} \Bi^{-2}(\sigma^{-\frac{2}{3}}v)v \partial_{\sigma} \{\dot{\tilde{a}}_2(\sigma,-v\} \: dv \\
	+\frac{\Bi}{\Ai}(\tau)\int_{-\zeta}^{\infty} \frac{\Ai}{\Bi}(\sigma^{-\frac{2}{3}}v)\partial_{\sigma}^2\{\dot{\tilde{a}}_2(\sigma,-v)\} \: dv\,.
\end{align*}
The first two terms are easily bounded by $ \sigma^{-\frac{10}{3}} $. Using that $ \partial_{\sigma}^2\{\dot{\tilde{a}}_2 \}(\sigma,-v)=\frac{\partial}{\partial u}\{\partial_{\sigma}^2\{\tilde{a}_2(\sigma,u)\} \}\vert_{u=-v}   $, we integrate by parts in the last term to rewrite it as
\begin{align*}
	\frac{\Bi}{\Ai}(\tau)\left(\frac{\Ai}{\Bi}(\tau)\partial_{\sigma}^2\{\tilde{a}_2(\sigma,\zeta)\}-\pi^{-1}\sigma^{-\frac{2}{3}}\int_{-\zeta}^{\infty} \Bi^{-2}(\sigma^{-\frac{2}{3}}v)\partial_{\sigma}^2\{\tilde{a}_2(\sigma,-v)\} \: dv \right)\,,
\end{align*}
which is controlled by $ \partial_{\sigma}^2\{\tilde{a}_2(\sigma,\zeta)\}= O(\sigma^{-3}) $. It follows then that $ \partial_{\sigma}^2\{\tilde{a}_1\} =O(\sigma^{-\frac{10}{3}})  $. 

We finish the proof of this lemma by computing $\partial_\sigma^2\{\dot{\tilde{a}}_1(\sigma,\zeta)\} $. By \eqref{dersigmazetatilde} and collecting previous bounds we obtain 
\begin{align*}
\partial^2_{\sigma} \{\dot{\tilde{a}}_1(\sigma,\zeta)\} = +\frac{2}{3\pi}\sigma^{-4}\zeta^2  \frac{d^2}{du^2}\{(\Ai\Bi)^{-1}(u)\}\vert_{u=\tau}\tilde{a}_1(\sigma,\zeta) 
	+O(\sigma^{-\frac{13}{3}}) \nn \,,
\end{align*}
 Using $| \frac{d^2}{du^2}\{(\Ai\Bi)^{-1}(u)\}| \les \la u \ra^{-\f32}$
we see that the first term is also bounded by  $\sigma^{-\frac{13}{3}}$. 
This was the last bound we needed to demonstrate, so the proof of the lemma is complete.
\end{proof}

\begin{cor}\label{AiryCor}
	Let $ \sigma $ be sufficiently small. Then for $ x\in[\frac{1}{2},1+\delta] $,
\begin{align}\label{esmall}
	e(\sigma,r(x))&=A(\sigma) (q(x))^{-\frac{1}{4}}\phi_1(\sigma,x)+B(\sigma)(q(x))^{-\frac{1}{4}}\phi_2(\sigma,x)\,,\\
	A(\sigma)&= \sigma^{-\frac{1}{6}}(1+e_1(\sigma)),\,\,B(\sigma)= \sigma^{\frac{5}{6}}e^{-\sigma^{-1}(\frac{\pi}{2}-1)} (1+e_2(\sigma) ) \nn \,.
\end{align}
where 	
\begin{align}\label{lbound}
|e_j(\sigma)| \les \sigma,\,\, |e_j^{\prime  }(\sigma)| \les \sigma^{-\frac{1}{3}}, \,\, |e_j^{\prime \prime }(\sigma)| \les \sigma^{-\frac{7}{3}},\, \text{for}\,\ j=1,2.
\end{align}
\end{cor}
\begin{proof}
We match  the Bessel function approximation of $ e(\sigma,r(x)) $ to the basis $ \{q^{-\frac{1}{4}}\phi_1(\tau),q^{-\frac{1}{4}}\phi_2(\tau)\} $ at $ x=\frac{1}{2} $. We have that
\begin{align}\label{BesselHalf}
	e(\sigma,r(1/2))= \sqrt{2} \sigma^{-\f 12} [ e^{\frac {\pi}{\sigma}} -1]^{-\f 12}\sqrt{\eta_*}I_1(\sigma^{-1} \eta_*)(1+a_-(\sigma,\eta_*))\,,
\end{align}
where $ \eta_*=\eta(\frac{1}{2})= \frac{1}{2}- \frac{\pi}{4} $ and we have used that $ \eta'(\frac{1}{2})=1 $ because $ \eta'=\sqrt{Q}  $. Using \cite[(10.40.1)]{NIST}, we have from \eqref{BesselHalf} that
\begin{align*}
	e(\sigma,r(\frac{1}{2}))=
 \frac{1}{\sqrt{\pi}}e^{\frac{\eta_*}{\sigma}}[e^{\frac{\pi}{\sigma}}-1]^{-\frac{1}{2}}(1+l_1(\sigma))
=\frac{e^{\sigma^{-1}(\eta_*-\frac{\pi}{2})}}{\sqrt{\pi}}(1+l_2(\sigma))
\end{align*}
with $l_i=O_1(\sigma)$ and $|l_i^{\prime \prime}(\sigma)|\les \sigma^{-2}$ for $i=1,2$. Moreover, 
\begin{align*}
\partial_x[e(\sigma,r(x))]_{x=\frac{1}{2}} &=\sqrt{2} [ e^{\frac {\pi}{\sigma}} -1]^{-\f 12}\sqrt{\eta_*}\sigma^{-\frac{3}{2}}I_1'(\sigma^{-1}\eta_*)(1+l_3(\sigma)) \\
 &=\frac{1}{\sqrt{\pi}\sigma}e^{\frac{\eta_*}{\sigma}}[e^{\frac{\pi}{\sigma}}-1]^{-\frac{1}{2}}(1+l_4(\sigma))
=\frac{e^{\sigma^{-1}(\eta_*-\frac{\pi}{2})}}{\sqrt{\pi}\sigma}(1+l_5(\sigma))\,.
\end{align*}
for $l_j(\sigma)$, $j=3,4,5$ with the same bound as $l_1(\sigma)$.  We also have that
\begin{align*}
	\phi_1(\sigma,\zeta(\frac{1}{2}))=q(\frac{1}{2})^{-\frac{1}{4}}\Ai(\sigma^{-\frac{2}{3}}\zeta_*)(1+\sigma a_1(\sigma,\zeta_*))\\
	\phi_2(\sigma,\zeta(\frac{1}{2}))=q(\frac{1}{2})^{-\frac{1}{4}}\Bi(\sigma^{-\frac{2}{3}}\zeta_*)(1+\sigma a_2(\sigma,\zeta_*))\,,
\end{align*}
where $ \zeta_*=-\zeta(\frac{1}{2}) $. Because $ \zeta_*>0 $, we use the asymptotics \cite[(9.7.5) and (9.7.6)]{NIST} to see that
\begin{align*}
	&\phi_1(\sigma,\zeta(x))=\sigma^{\frac{1}{6}}\frac{e^{-\frac{2\zeta_*^{\frac{3}{2}}}{3\sigma}}}{2\sqrt{\pi} }(1+l_6(\sigma))\\
	&\phi_2(\sigma,\zeta(x))=\sigma^{\frac{1}{6}}\frac{e^{\frac{2\zeta_*^{\frac{3}{2}}}{3\sigma}}}{\sqrt{\pi} }(1+l_7(\sigma)) \,
\end{align*}
where $l_6,l_7$ hold the bounds in \eqref{lbound}. Moreover,  \cite[(9.7.6),(9.7.8)]{NIST} give 
\begin{align*}
	&\partial_x[\phi_1(\sigma,\zeta(x))]_{x=\frac{1}{2}}=-q^{-\frac{1}{4}}(\frac{1}{2})\zeta'(\frac{1}{2})\sigma^{-\frac{2}{3}}\Ai'(\sigma^{-\frac{2}{3}}\zeta_*)(1+l_8(\sigma))=\sigma^{-\frac{5}{6}}\frac{e^{-\frac{2\zeta_*^{\frac{2}{3}}}{3\sigma}}}{2\sqrt{\pi} }(1+l_9(\sigma))\\
	&\partial_x[\phi_2(\sigma,\zeta(x))]_{x=\frac{1}{2}}=-q^{-\frac{1}{4}}(\frac{1}{2})\zeta'(\frac{1}{2})\sigma^{-\frac{2}{3}}\Bi'(\sigma^{-\frac{2}{3}}\zeta_*)(1+l_{10}(\sigma)) =-\sigma^{-\frac{5}{6}}\frac{e^{\frac{2\zeta_*^{\frac{2}{3}}}{3\sigma}}}{\sqrt{\pi} }(1+l_{11}(\sigma))
\end{align*}
employing that $ q(\frac{1}{2})=\zeta_*^{-1}=[\zeta'(\frac{1}{2})]^2 $. Similarly, $l_9,l_{11}$ hold the bounds in \eqref{lbound} It follows that for some $l_{12}$ holding the bounds in \eqref{lbound}
\begin{align*}
	W\left[ \phi_1(\sigma,\zeta(\cdot)),\phi_2(\sigma,\zeta(\cdot)) \right]=-\frac{\sigma^{-\frac{2}{3}}}{\pi}(1+l_{12}(\sigma))\,
\end{align*}
where $ W $ is the Wronskian evaluated at $ x=\frac{1}{2} $, and also
\begin{align*}
	&W[e(\sigma,r(\cdot)),\phi_1(\sigma,\zeta(\cdot))]
 = \sigma^{\frac{1}{6}}e^{-\sigma^{-1}(\frac{\pi}{2}-1)} (1+l_{13}(\sigma)) \\
	&W[e(\sigma,r(\cdot)),\phi_2(\sigma,\zeta(\cdot))]= 
 -\frac{\sigma^{-\frac{5}{6}}}{\sqrt{\pi} }(1+l_{14}(\sigma))\,,
\end{align*}
where the last equality on each line follows from the facts that $ \eta_*-\frac{2}{3}\zeta_*^{\frac{3}{2}}=1 $ and $ \eta_*+\frac{2}{3}\zeta_*^{\frac{3}{2}}=\frac{\pi}{2} $. Here, $l_{13}$ and $l_{14}$ hold the bounds in \eqref{lbound}.   Therefore, 
\begin{align*}
	&A(\sigma)=\frac{W\left[e(\sigma,r(\cdot)),\phi_2(\sigma,\zeta(\cdot))\right]}{W\left[ \phi_1(\sigma,\zeta(\cdot)),\phi_2(\sigma,\zeta(\cdot)) \right]}= \sigma^{-\frac{1}{6}}(1+e_1(\sigma))\,,\\
	&B(\sigma)=-\frac{W\left[e(\sigma,r(\cdot)),\phi_1(\sigma,\zeta(\cdot))\right]}{W\left[ \phi_1(\sigma,\zeta(\cdot)),\phi_2(\sigma,\zeta(\cdot)) \right]}= \sigma^{\frac{5}{6}}e^{-\sigma^{-1}(\frac{\pi}{2}-1)} \left( 1+ e_2(\sigma) \right)\,.
\end{align*}
\end{proof}\subsection{Oscillatory Airy approximation: $ x\gg 1$ }
\begin{pr}\label{Vbound} 
When  $\zeta\geq 0 $, the potential $ V $ satisfies the bounds
\begin{align*}
	|\partial^j_\zeta V(\zeta)|\les \la \zeta \ra^{-2-j}\,\,\ j=0,1,2. 
\end{align*}
\end{pr}
\begin{proof}
Since by the above $ \zeta' $ is smooth and non-vanishing, the identity $ q=(\zeta')^2 $ shows that near $ 0 $, $ V=-q^{-\frac{1}{4}}\frac{d^2q^{\frac{1}{4}}}{d\zeta^2} $ is bounded, so we need only show that $ V $ has the claimed behavior as $ \zeta\rightarrow\infty $. With $ { }^.=\frac{d}{d\zeta} $, one computes first that
\begin{align}\label{Vfromq}
	q^{-\frac{1}{4}}\frac{d^2q^{\frac{1}{4}}}{d\zeta^2}=-\frac{3}{16}q^{-2} \dot{q}^{2}+\frac{1}{4}q^{-1}\ddot{q}\,,
\end{align}
so we need only find asymptotics for $ q $ in terms of $ \zeta $. Recalling definition of $ \zeta $ as a function of $ x $, we have that for $ x\geq 1 $ 
\begin{align*}
	\frac{3}{2}\zeta^{\frac{3}{2}}(x)=\int_{1}^{x} \sqrt{1-u^{-1}} \: du=\int_{1}^{x} 1+O(u^{-1})\: du=x+c+O(x^{-1}) \,.
\end{align*}
The chain rule applied to the above equality then shows that $ \zeta'(x)=O(\zeta^{-\frac{1}{2}}) $, where every derivative of $ O(\zeta^{-\frac{1}{2}}) $ loses a power of $ \zeta $, that is, it exhibits symbol behavior. It follows then that $ q=(\zeta')^2=O(\zeta^{-1}) $ for $ x $ large, from which the bound on $ V $ follows from \eqref{Vfromq} and one may obtain the bounds on $ V' $ and $ V^{''} $  by differentiating this equality.

\end{proof}
\begin{pr}\label{oscBasisPr}
	For any $ \sigma>0 $ sufficiently small, the following holds: in the range $ \zeta\geq 0 $, a basis of solutions to \eqref{zetaODE} is given by
	\begin{align}
		\begin{split}\label{oscBasis}
		\psi_+=(\Ai(\tau)+i\Bi(\tau))[1+\sigma b_+(\sigma,\zeta)]\\
		\psi_-=(\Ai(\tau)-i\Bi(\tau))[1+\sigma b_-(\sigma,\zeta)]\,,
	\end{split}
	\end{align}
with $ \tau=-\sigma^{-\frac{2}{3}}\zeta $,
and $ b_\pm  $ smooth functions satisfying the bounds
\begin{align}\label{bpmbound}
\begin{split}
	&\abs{b_{\pm}(\sigma,\zeta)} \les \Span{\zeta}^{-\frac{3}{2}},\,\, \abs{\dot{b}_{\pm}(\sigma,\zeta)} \les \sigma^{-\frac{1}{3}}\Span{\sigma^{-\frac{2}{3}}\zeta}^{-\frac{1}{2}}\Span{\zeta}^{-2},\,\, \abs{\ddot{b}_{\pm}(\sigma,\zeta)} \les \sigma^{-1}\Span{\zeta}^{-2}\\
	&\partial_{\sigma}[b_{\pm}(\sigma,\zeta)]\les \sigma^{-1}\Span{\zeta}^{-\frac{3}{2}},\,\,\partial_{\sigma}^2[b_{\pm}(\sigma,\zeta)]\les \sigma^{-3}, \,\,\, \partial_{\sigma}[\dot{b}_{\pm}(\sigma,\zeta)]\les \sigma^{-2}\Span{\zeta}^{-1},\\   & \partial^2_{\sigma}[\dot{b}_{\pm}(\sigma,\zeta)]\les \sigma^{-3} \,.
 \end{split}
\end{align}
\begin{rmk}
	This proposition is similar to \cite[Proposition 9]{Wmain}. Indeed, the bounds on $ b_{\pm} $ and $ \dot{b}_{\pm} $ are produced via the same proof, as the only inputs are the asymptotics of $ V $, which in this regime of $ \zeta $ are the same as in that paper. However, for our purposes we also require an additional derivative in $ \zeta $ and derivatives in the semi-classical parameter $ \sigma $ ($ \hbar $ in \cite{Wmain}).  Note that in both settings, a representation of the form \eqref{oscBasis} is only possible because $ \Ai $ and $ \Bi $ have no common zeroes, and therefore \eqref{oscBasis} does not fix the zeroes of any solution.
\end{rmk}
\end{pr}
\begin{proof}
	Let $ \psi_{\pm,0}(\zeta,\sigma)=\Ai(\tau)\pm i\Bi(\tau) $. Similar to \eqref{ajEq}, we obtain the equations
	\begin{align*}
		\left( \psi_{\pm,0}^2\dot{b}_\pm \right)^.=-\sigma^{-1}V\psi_{\pm,0}^2(1+\sigma b_\pm)
	\end{align*}
	whose solutions with $ b_\pm(\infty)=0 $ and $ \dot{b}_\pm(\infty)=0 $ are given by
	\begin{align}\label{voltEq}
	b_{\pm}(\zeta)=-\sigma^{-1}\int_{\zeta}^{\infty}\int_{\zeta}^{u} \psi_{\pm,0}^{-2}(v)\: dv\:\psi_{\pm,0}^2(u)V(u)(1+\sigma b_\pm(u))  \: du\,,
	\end{align}
	where for now we have suppressed the dependence on $ \sigma $ in the integrand. From \cite{Olver}, we have the asymptotic expansion
	\begin{align}\label{oscAiry}
	\Ai(-z)\pm i\Bi(-z)=\frac{1}{\sqrt{\pi} z^{\frac{1}{4}}}e^{\mp i(\frac{2}{3}z^{\frac{3}{2}}-\frac{\pi}{4})}(1+O(z^{-\frac{3}{2}}))\,,
	\end{align}
where the $ O(z^{-\frac{3}{2}}) $ term may be differentiated as a symbol. Thus, for $ 0<x_0<x_1 $ 
\begin{multline}\label{psiinverse}
	\int_{x_0}^{x_1} (\Ai(-z)\pm\Bi(-z))^{-2}\: dz \\ =\int_{x_0}^{x_1} z^{\frac{1}{2}}e^{\mp i\frac{4}{3}z^{\frac{3}{2}}}a(z)\: dz =\mp\frac{1}{2i}e^{\mp i \frac{4}{3}z^{\frac{3}{2}}}a(z)\bigg\vert_{x_0}^{x_1}\pm \frac{1}{2i}\int_{x_0}^{x_1} e^{\mp i \frac{2}{3}z^{\frac{3}{2}}}a'(z)\: dz
\end{multline}
for $ a(z)=1+O(z^{-\frac{3}{2}}) $. This shows that the above integral is $ O(1) $ for all $ x_0 $ and $ x_1 $. The main term with respect to $ \sigma $ in \eqref{voltEq}
\begin{align*}
	b_{\pm,0}(\zeta)=-\sigma^{-1}\int_{\zeta}^{\infty}\int_{\zeta}^{u} \psi_{\pm,0}^{-2}(v)\: dv\:\psi_{\pm,0}^2(u)V(u)\:du
\end{align*}
satisfies the bound
\begin{align*}
	\abs{b_{\pm,0}(\zeta)} \les \sigma^{\frac{1}{3}}\int_{\sigma^{-\frac{2}{3}}\zeta}^{\infty} \Span{u}^{-\frac{1}{2}}\abs{V(-\sigma^{\frac{2}{3}}u)} \: du \,,
\end{align*}
where we have changed variables and used that by \eqref{oscAiry},
\begin{align*}
	\abs{(\Ai(-z)+i\Bi(-z))^{2}} \les \Span{z}^{-\frac{1}{2}} \,.
\end{align*}
By Proposition~\ref{Vbound}, we see then that
\begin{align*}
	\abs{b_{\pm,0}(\zeta)} \les \sigma^{\frac{1}{3}}\int_{\sigma^{-\frac{2}{3}}\zeta}^{\infty} \Span{u}^{-\frac{1}{2}}\Span{\sigma^{\frac{2}{3}}u}^{-2}\: du \les \Span{\zeta}^{-\frac{3}{2}}\,,
\end{align*}
where the last inequality comes from bounding the integrand by $ \sigma^{-\frac{4}{3}}u^{-\frac{5}{2}} $ when $ \sigma^{-\frac{2}{3}}\zeta $ is large.  We can now extend this bound to $b_{\pm}$ by a contraction argument, as is explained in Proposition~\ref{prop:bessel}, by considering the linear operator $$ Ta= -\sigma^{-1}\int_{\zeta}^{\infty}\int_{\zeta}^{u} \psi_{\pm,0}^{-2}(v)\: dv\:\psi_{\pm,0}^2(u)V(u)a(u)  \: du$$ as a map on the weighted space $ \Span{\zeta}^{-\frac{3}{2}} L^\infty_\zeta $.

For the $ \zeta $-derivative bounds, we first write
\begin{align}\label{bdotInt}
	\dot{b}_\pm(\zeta)=\sigma^{-1}\psi_{\pm,0}^{-2}(\zeta)\int_{\zeta}^{\infty} \psi_{\pm,0}^2(u)V(u)(1+\sigma b_\pm(u))\: du 
\end{align}
and use \eqref{oscAiry} to see that $\psi_{\pm,0}^2(\zeta)=e^{i\frac{4}{3\sigma}\zeta^{\frac{3}{2}}}\omega(\sigma^{-\frac{2}{3}}\zeta) $ for some $ \omega(u) $ with $ \abs{\omega(u)} \les \Span{u}^{-\frac{1}{2}}$ and $ \abs{\omega'(u)} \les \Span{u}^{-\frac{3}{2}}$.  When $ \zeta >\sigma^{-\frac{2}{3}}$, we may exploit the oscillatory phase by integrating by parts in the above integral via 
\begin{align*}
	\psi_{\pm,0}^2=(2i\zeta^{\frac{1}{2}}/\sigma)^{-1}\omega(\sigma^{-\frac{2}{3}}\zeta)\frac{d}{d\zeta}[e^{i\frac{4}{3\sigma}\zeta^{\frac{3}{2}}}]
\end{align*}
to find that
\begin{align}
	\dot{b}_\pm(\zeta)=\frac{1}{2i}\psi^{-2}_{\pm,0}(\zeta)[u^{-\frac{1}{2}}\omega(\sigma^{-\frac{2}{3}}u)e^{i\frac{4}{3\sigma}\zeta^{\frac{3}{2}}}V(u)(1+\sigma b_{\pm}(u))]\bigg\vert_{\zeta}^\infty\label{bdot_1}\\
	-\frac{1}{2i}\psi_{\pm,0}^{-2}(\zeta)\int_{\zeta}^{\infty}e^{i\frac{4}{3\sigma}\zeta^{\frac{3}{2}}}\frac{d}{du}[u^{-\frac{1}{2}}\omega(\sigma^{-\frac{2}{3}}u)V(u)(1+\sigma b_{\pm}(u))] \: du \label{bdot_2}\,.
\end{align}
Using the bounds on $\omega$ we see that the term on the right hand side of the equality in \eqref{bdot_1} is bounded by
\begin{align*}
	\zeta^{-\frac{1}{2}}\Span{\sigma^{-\frac{2}{3}}\zeta}^{-\frac{1}{2}}\Span{\zeta}^{-2}\les \sigma^{-\frac{1}{3}}\Span{\sigma^{-\frac{2}{3}}\zeta}^{-\frac{1}{2}}\Span{\zeta}^{-2}\,,
\end{align*}
with the last inequality coming from the assumption on $ \zeta $. By differentiating the product in the integrand of \eqref{bdot_2}, one checks via the bounds on $ \omega $ and Proposition~\ref{Vbound} that every term other than the term in which the derivative falls on $ b_\pm $  is $ O(\sigma^{-\frac{1}{3}}\Span{\sigma^{-\frac{2}{3}}\zeta}^{-\frac{1}{2}}\Span{\zeta}^{-2}) $ so that we may write
\begin{align*}
\dot{b}_\pm(\zeta)= O(\sigma^{-\frac{1}{3}}\Span{\sigma^{-\frac{2}{3}}\zeta}^{-\frac{1}{2}}\Span{\zeta}^{-2}) 	-\frac{\sigma}{2i}\psi_{\pm,0}^{-2}(\zeta)\int_{\zeta}^{\infty}e^{i\frac{4}{3\sigma}\zeta^{\frac{3}{2}}}O(u^{-\frac{1}{2}} \Span{\sigma^{-\frac{2}{3}}u}^{-\frac{1}{2}}\Span{u}^{-2})\dot{b}_\pm(u) \: du\,.
\end{align*}
By iterating this equality, we see that the second term is better than the first one, so  we see that $| \dot{b}_{\pm}(\zeta)| \les\sigma^{-\frac{1}{3}}\Span{\sigma^{-\frac{2}{3}}\zeta}^{-\frac{1}{2}}\Span{\zeta}^{-2}$ . As for the case when $\sigma^{-\frac{2}{3}}\zeta \leq 1$, we simply write \eqref{bdotInt} as
\begin{align*}
	\sigma^{-1}\psi_{\pm,0}^{-2}(\zeta)\int_{\zeta}^{\sigma^{\frac{2}{3}}} \psi_{\pm,0}^2(u)V(u)(1+\sigma b_\pm(u))\: du\\ 
	+\sigma^{-1}\psi_{\pm,0}^{-2}(\zeta)\int_{\sigma^{\frac{2}{3}}}^{\infty} \psi_{\pm,0}^2(u)V(u)(1+\sigma b_\pm(u))\: du \,,
\end{align*}
where the first term is clearly bounded by $ \sigma^{-\frac{1}{3}} $ and the second one can be estimated similar to \eqref{bdot_1}. Thus, in either case we see that
\begin{align*}
	\abs{\dot{b}_{\pm}(\zeta)} \les \sigma^{-\frac{1}{3}}\Span{\sigma^{-\frac{2}{3}}\zeta}^{-\frac{1}{2}}\Span{\zeta}^{-2}
\end{align*}
as claimed. 

For the second $ \zeta$-derivative, we differentiate \eqref{bdotInt} to find that
\begin{align*}
	\ddot{b}_{\pm}(\zeta)=-\sigma^{-1} V(\zeta)(1+\sigma b_{\pm}(\zeta))-2\dot{\psi}_{\pm,0}(\zeta)\psi_{\pm,0}^{-1}(\zeta)\dot{b}_{\pm}(\zeta)\,.
\end{align*}
Since $ \abs{\psi_{\pm,0}^{-1}(\zeta) } \les \Span{\sigma^{-\frac{2}{3}}\zeta}^{\frac{1}{4}}$ and by \eqref{oscAiry} $ \abs{\dot{\psi}_{\pm,0}(\zeta)} \les \sigma^{-\frac{2}{3}}\Span{\sigma^{-\frac{2}{3}}\zeta}^{\frac{1}{4}} $, the previously derived bounds on $ b_\pm$ and $\dot{b}_{\pm} $ show that
\begin{align*}
	\abs{\ddot{b}_{\pm}(\zeta)} \les \sigma^{-1}\Span{\zeta}^{-2}\,.
\end{align*}\par
We now demonstrate the bounds on the $ \sigma $-derivatives of $ b_{\pm} $. To begin with, we rewrite \eqref{voltEq} as
\begin{align*}
	b_\pm(\sigma,\zeta)&=-\sigma^{\frac{1}{3}}\int_{\sigma^{-\frac{2}{3}}\zeta}^{\infty}\int\limits_{\sigma^{-\frac{2}{3}}\zeta}^{u} (\Ai\pm i\Bi)^{-2}(-v)\: dv (\Ai\pm i\Bi)^2(-u)V(\sigma^{\frac{2}{3}}u)(1+\sigma b(\sigma,\sigma^{\frac{2}{3}}u))\: du \\
	&=-\sigma^{\frac{1}{3}}\int_{\sigma^{-\frac{2}{3}}\zeta}^{\infty}\int\limits_{\sigma^{-\frac{2}{3}}\zeta}^{u} (\Ai\pm i\Bi)^{-2}(-v)\: dv (\Ai\pm i\Bi)^2(-u)F(u,\sigma) \:du\,,
\end{align*}
and then differentiate with respect to $ \sigma $ to find that
\begin{align*}
	\partial_{\sigma}[b_{\pm}(\sigma,\zeta)]&=\frac{1}{3}\sigma^{-1}b_{\pm}(\sigma,\zeta)\\
	&-\frac{2}{3}\sigma^{-\frac{4}{3}}\psi_{\pm,0}^{-2}(\zeta)\zeta\int_{\sigma^{-\frac{2}{3}}\zeta}^{\infty} (\Ai\pm i\Bi)^2(-u)F(u,\sigma)\: du \\
	&-\sigma^{\frac{1}{3}}\int_{\sigma^{-\frac{2}{3}}\zeta}^{\infty}\int\limits_{\sigma^{-\frac{2}{3}}\zeta}^{u} (\Ai\pm i\Bi)^{-2}(-v)\: dv (\Ai\pm i\Bi)^2(-u)\partial_\sigma[F(u,\sigma)] \:du  \\ 
 &=: F_1^{\pm}(\sigma,\zeta)+ F_2^{\pm}(\sigma,\zeta)+F_3^{\pm}(\sigma,\zeta).
\end{align*}

Our previous bound on $ b$ shows that $ |F^{\pm}_1(\sigma,\zeta)|\les \sigma^{-1}\Span{\zeta}^{-\frac{3}{2}} $. Also, from \eqref{bdotInt}, we see that $F_2^{\pm}(\sigma,\zeta)= -\frac{2}{3}\sigma^{-1}\zeta\dot{b}_{\pm}(\sigma,\zeta) $ and is therefore bounded by $ \sigma^{-\frac{4}{3}}\Span{\sigma^{-\frac{2}{3}}\zeta}^{-\frac{1}{2}}\Span{\zeta}^{-2}\zeta$, which is again bounded by $ \sigma^{-1}\Span{\zeta}^{-\frac{3}{2}} $. We compute that
\begin{align}
\begin{split}
\partial_\sigma[F(u,\sigma)]&=
\frac{2}{3}\sigma^{-\frac{1}{3}}V'(\sigma^{\frac{2}{3}}u)u(1+\sigma b_{\pm}(\sigma,\sigma^{\frac{2}{3}}u))+V(\sigma^{\frac{2}{3}}u)b_{\pm}(\sigma,\sigma^{\frac{2}{3}}u)\\
				    &+\frac{2}{3}\sigma^{\frac{2}{3}}V(\sigma^{\frac{2}{3}}u)\dot{b}_{\pm}(\sigma,\sigma^{\frac{2}{3}}u)+\sigma V(\sigma^{\frac{2}{3}}u)\partial_{\sigma}[b_{\pm}(\sigma,v)]_{v=\sigma^{\f23}u}
\end{split}\label{sigmaF}\\
				    &=O\left( \sigma^{-\frac{1}{3}}\Span{u}\Span{\sigma^{\frac{2}{3}}u}^{-3} \right)+ O\left( \Span{\sigma^{\frac{2}{3}}u}^{-\frac{7}{2}} \right) + O\left( \sigma^{\frac{1}{3}}\Span{u}^{-\frac{1}{2}}\Span{\sigma^{\frac{2}{3}}u}^{-4} \right) \nonumber\\
				    &+\sigma V(\sigma^{\frac{2}{3}}u)\partial_{\sigma}[b_{\pm}(\sigma,v)]_{v=\sigma^{\f23}u} \nonumber\\
				 &=O\left( \sigma^{-\frac{1}{3}}\Span{u}\Span{\sigma^{\frac{2}{3}}u}^{-3} \right)+\sigma V(\sigma^{\frac{2}{3}}u)\partial_{\sigma}[b_{\pm}(\sigma,v)]_{v=\sigma^{\f23}u}\nonumber\,.
\end{align}
Inserting the last line into $F_3^{\pm}(\sigma,\zeta)$ shows that
\begin{align*}
	F_3^{\pm}(\sigma,\zeta)&=\int\limits_{\sigma^{-\frac{2}{3}}\zeta}^{\infty} O\left( \Span{u}^{\frac{1}{2}}\Span{\sigma^{\frac{2}{3}}u}^{-3} \right) \: du +T\Big(\partial_\sigma[b_{\pm}(\sigma,v)]_{v=\sigma^{\f23}u}\Big)\\
 &=O\left( \sigma^{-1}\Span{\zeta}^{-\frac{3}{2}} \right) +T\Big(\partial_\sigma[b_{\pm}(\sigma,v)]_{v=\sigma^{\f23}u}\Big)\,,
\end{align*}
with
\begin{align*}
	T(a)=-\sigma^{\frac{4}{3}}\int_{\sigma^{-\frac{2}{3}}\zeta}^{\infty}\int\limits_{\sigma^{-\frac{2}{3}}\zeta}^{u} (\Ai\pm i\Bi)^{-2}(-v)\: dv (\Ai\pm i\Bi)^2(-u) V(\sigma^{\frac{2}{3}}u)a(\sigma^{\frac{2}{3}}u) \:du\,.
\end{align*}

Collecting the above bounds, we arrive that
\begin{align}
	\partial_{\sigma}[b_{\pm}(\sigma,\zeta)]&=O\left( \sigma^{-1}\Span{\zeta}^{-\frac{3}{2}}\right)+T\Big(\partial_\sigma[b_{\pm}(\sigma,v)]_{v=\sigma^{\f23}u}\Big)\,.\label{sigmabFixedPt}
\end{align}
For small enough $ \sigma $, $ T $ is a contraction on the weighted space $ \Span{\zeta}^{\frac{3}{2}}L^\infty_\zeta $ because
\begin{align*}
	\abs{T(\Span{\zeta}^{-\frac{3}{2}})} \les \sigma^{\frac{4}{3}}\int_{\sigma^{-\frac{2}{3}}\zeta}^{\infty} \Span{u}^{-\frac{1}{2}}\Span{\sigma^{\frac{2}{3}}u}^{-\frac{7}{2}}\: du \les \sigma \Span{\zeta}^{-3}\,,
\end{align*}
so we conclude that the first term in \eqref{sigmabFixedPt} bounds the second, \emph{i.e.}
\begin{align*}
	\abs{\partial_{\sigma}[b_{\pm}(\sigma,\zeta)]} \les \sigma^{-1}\Span{\zeta}^{-\frac{3}{2}}\,.
\end{align*}\par
The second $ \sigma $-derivative will require an estimate on the mixed derivative $ \partial_\sigma[\dot{b}(\sigma,\zeta)] $. This easily follows from the bounds we have in hand by differentiating in $ \zeta $ each of $F_i^{\pm}(\sigma,\zeta)$ for $i=1,2,3$. Clearly $ \partial_{\zeta}[F_1^{\pm}(\sigma,\zeta)] $ contributes $ \sigma^{-\frac{4}{3}}\Span{\sigma^{-\frac{2}{3}}\zeta}^{-\frac{1}{2}}\Span{\zeta}^{-2} $, while it is easy to see that $  $ 
\begin{align*}
	\partial_{\zeta} [ F_2^{\pm}(\sigma,\zeta) ]&=-2\psi^{-1}_{\pm,0}(\zeta)\dot{\psi}_{\pm,0}(\zeta)F_2^{\pm}( \sigma,\zeta) )+\zeta^{-1}F_2^{\pm}(\sigma,\zeta)+\frac{2}{3}\sigma^{-2}\zeta V(\zeta)(1+\sigma b_{\pm}(\sigma,\zeta))\\
	&=O\left( \sigma^{-2} \Span{\zeta}^{-1}  \right) \,.
\end{align*}
Finally, the bound we have obtained on $ \partial_{\sigma}[b(\sigma,\zeta)] $ may be used in \eqref{sigmaF} to show that
\begin{align}\label{F3}
	\partial_{\zeta}[ F_3^{\pm}(\sigma,\zeta)]&=\sigma^{-\frac{1}{3}}\psi^{-2}_{\pm,0}(\zeta)\int\limits_{\sigma^{-\frac{2}{3}}\zeta}^{\infty} (\Ai \pm i \Bi)^2(-u)\partial_\sigma[F(u,\sigma)]\: du \\
					    &\les \sigma^{-\frac{2}{3}} \Span{\sigma^{-\frac{2}{3}}\zeta}^{\frac{1}{2}}\int\limits_{\sigma^{-\frac{2}{3}}\zeta}^{\infty} \Span{u}^{\frac{1}{2}}\Span{\sigma^{\frac{2}{3}}u}^{-3} \: du \les \sigma^{-2}\Span{\zeta}^{-1} 
\end{align}
for a total bound of 
\begin{align*}
	\partial_{\sigma}[\dot{b}(\sigma,\zeta)]\les \sigma^{-2}\Span{\zeta}^{-1} \,.
\end{align*}
\par
For the second $ \sigma $-derivative, we proceed similarly and differentiate each of $F_i^{\pm}(\sigma,\zeta)$. First,
\begin{align*}
	\partial_\sigma[F_1^{\pm}(\sigma,\zeta)]=-\sigma^{-1}(F_1^{\pm}(\sigma,\zeta))+\frac{1}{3}\sigma^{-1}\partial_\sigma[b_\pm(\sigma,\zeta)]\les \sigma^{-2}\Span{\zeta}^{-\frac{3}{2}}\,,
\end{align*}
whereas
\begin{align*}
	\partial_{\sigma}[F_2^{\pm}&(\sigma,\zeta)]=-\frac{4}{3}\sigma^{-1}F_2^{\pm}(\sigma,\zeta)-\frac{4}{3}\sigma^{-\frac{5}{3}}\zeta\psi_{\pm,0}^{-1}(\zeta)(\Ai\pm i \Bi)'(-\sigma^{-\frac{2}{3}}\zeta)F_2^{\pm}(\sigma,\zeta)\\
					   &-\frac{4}{9}\sigma^{-3}\zeta^2V(\zeta)(1+\sigma b_{\pm}(\sigma,\zeta))-\frac{2}{3}\sigma^{-\frac{4}{3}}\psi^{-2}_{\pm,0}(\zeta)\zeta\int\limits_{\sigma^{-\frac{2}{3}}\zeta}^{\infty} (\Ai\pm\Bi)^2(-u)\partial_\sigma\{F(u,\sigma)\}\: du \,.
\end{align*}
The first term is $ O\left( \sigma^{-2}\Span{\zeta}^{-\frac{3}{2}}  \right) $ and second term is $ O\left( \sigma^{-3} \right)  $ because we have already shown that $|F_2^{\pm}(\sigma,\zeta)|\les \sigma^{-\frac{4}{3}}\Span{\sigma^{-\frac{2}{3}}\zeta}^{-\frac{1}{2}}\Span{\zeta}^{-2}\zeta $. The third term is also $ O(\sigma^{-3}) $ and the fourth is bounded by $ \sigma^{-\frac{5}{3}}\psi_{\pm,0}^{-2}(\zeta)\zeta F_3^{\pm}(\sigma,\zeta) $ which is again $ O(\sigma^{-3}) $ by the previously obtained bound on $F_3^{\pm}(\sigma,\zeta)$.\par
We compute further that
\begin{align*}
	\partial_\sigma[F_3^{\pm}(\sigma,\zeta)]&=\frac{1}{3}\sigma^{-1}F_3^{\pm}(\sigma,\zeta)-\frac{2}{3}\sigma^{-\frac{4}{3}}\psi_{\pm,0}^{-2}(\zeta)\int\limits_{\sigma^{-\frac{2}{3}}\zeta}^{\infty} (\Ai\pm i\Bi)^2(-u)\partial_\sigma[F(u,\sigma)]\: du\\
	&-\sigma^{\frac{1}{3}}\int_{\sigma^{-\frac{2}{3}}\zeta}^{\infty}\int\limits_{\sigma^{-\frac{2}{3}}\zeta}^{u} (\Ai\pm i \Bi)^{-2}(-v)\: dv (\Ai\pm i\Bi)^2(-u)\partial_\sigma^2[F(u,\sigma)]\: du \,.
\end{align*}
Arguing similarly, one may easily see that the first two terms are $ O(\sigma^{-3}) $. For the last term, we need to estimate $ \partial_\sigma^2[F(u,\sigma)] $. Series of elementary operations show that 

\begin{align*}
	\partial_{\sigma}^2[F(u,\sigma)] =O\left( \sigma^{-\frac{4}{3}}\Span{\sigma^{\frac{2}{3}}u}^{-2} \Span{u}^{}  \right) +\sigma V(\sigma^{\frac{2}{3}}u)\partial_{\sigma}^2[b_{\pm}(\sigma,v)]_{v=\sigma^{\frac{2}{3}}u}
\end{align*}
and therefore
\begin{align*}
	\partial_{\sigma}[F_3^{\pm}(\sigma,\zeta)]&=O(\sigma^{-3})\\
	&-\sigma^{\frac{1}{3}}\int_{\sigma^{-\frac{2}{3}}\zeta}^{\infty}\int\limits_{\sigma^{-\frac{2}{3}}}^{u} (\Ai\pm i \Bi)^{-2}(-v)\: dv (\Ai\pm i\Bi)^2(-u)O\left( \sigma^{-\frac{4}{3}}\Span{\sigma^{\frac{2}{3}}}^{-2} \Span{u}^{} \right)  \: du\\
	&+T\Big(\partial_{\sigma}^2[b_{\pm}(\sigma,v)]_{v=\sigma^{\frac{2}{3}}u}\Big)\,.
\end{align*}
The middle term is easily seen to be $ O\left( \sigma^{-2}\Span{\zeta}^{-\frac{1}{2}}  \right)  $ so that all of our estimates on the $ \sigma $-derivatives of $F_i^{\pm}(\sigma,\zeta)$ show that $ \partial_{\sigma}^2[b_{\pm}(\sigma,\zeta)] $ satisfies a fixed point equation of the form
\begin{align*}
	\partial_{\sigma}^2[b_{\pm}(\sigma,\zeta)]=O(\sigma^{-3})+T\Big(\partial_{\sigma}^2[b_{\pm}(\sigma,v)]_{v=\sigma^{\frac{2}{3}}u}\Big)\,.
\end{align*}
To conclude, we need only show that $ T $ is a contraction on $ L^\infty $ for $ \sigma $ sufficiently small, but this follows from the computation
\begin{align*}
	\abs{T(1)} \les \sigma^{\frac{4}{3}}\int_{\sigma^{-\frac{2}{3}}}^{\infty}\Span{u}^{-\frac{1}{2}} \Span{\sigma^{\frac{2}{3}}u}^{-2}  \: du\les \sigma\Span{\zeta}^{-\frac{3}{2}}
\end{align*}
so that, as before, it follows that $ \partial_\sigma^2[b_{\pm}(\sigma,\zeta)]=O(\sigma^{-3}) $. 

We finally consider $\partial^2_{\sigma}[\dot{b}_{\pm} (\sigma,\zeta)]$. For this, we need to differentiate $\partial_\zeta[F^{\pm}_i(\sigma,\zeta)]$ in $\sigma$ for $i=1,2,3$. A simple computation shows that $|\partial^2_{\sigma \zeta}[F^{\pm}_1(\sigma,\zeta)]| \les \sigma^{-3} \la \zeta \ra^{-1}$. We further compute
\begin{multline*}
\partial_{\sigma\zeta}[F_2(\sigma,\zeta)]= \partial_{\sigma}\{F_2^{\pm}(\sigma,\zeta)\} [-2\psi^{-1}_{\pm,0}(\zeta)\dot{\psi}_{\pm,0}(\zeta)+\zeta^{-1}] \\ + \f43 \sigma^{-1} \zeta \big(\psi^{-1}_{\pm,0}(\zeta)\dot{\psi}_{\pm,0}(\zeta)\big)^{\cdot}  F_2^{\pm}(\sigma,\zeta) + O(\sigma^{-3} \la \zeta \ra^{-1})
\end{multline*}

Recall that $F_2^{\pm}(\sigma,\zeta)= -\frac{2}{3}\sigma^{-1}\zeta\dot{b}{\pm}(\sigma,\zeta)$, and consequently, $|F_2^{\pm}(\sigma,\zeta)| \les \sigma^{-1} \zeta^{\frac{1}{2}}\langle \zeta \rangle^{-2}$. Furthermore, we estimate $\partial_{\sigma} \{F_2^{\pm}(\sigma,\zeta)\} \les \zeta \sigma^{-2}\langle \zeta \rangle^{-\frac{3}{2}}$, leading to the first term being bounded by $\sigma^{-3}$. Additionally, we may estimate $|\big(\psi^{-1}_{\pm,0}(\zeta)\dot{\psi}_{\pm,0}(\zeta)\big)^{\cdot}|\les \sigma^{-\frac{4}{3}} \langle \sigma^{-\frac{2}{3}} \zeta \rangle^{-\frac{3}{4}}$. Therefore, $|\partial^2_{\sigma\zeta}[F^{\pm}_2(\sigma,\zeta)]| \les \sigma^{-3}$. Finally, by differentiating \eqref{F3} and utilizing the previously obtained bounds, we find $\partial_{\sigma \zeta}[F^{\pm}_3(\sigma,\zeta)]= O(\sigma^{-3} \langle \zeta \rangle^{-1})$. This concludes the statement.
\end{proof}
We would now like to use the oscillatory basis to provide an approximation of $ e(\sigma,r) $ that is suitable for use inside the oscillatory integral defining $ K_t $. This requires detailed bounds on the function
\begin{align*}
	\zeta_r(\sigma):=\sigma^{-1}\frac{2}{3}\zeta^{\frac{3}{2}}(\sigma^2r) \,.
\end{align*}

\begin{pr}\label{langerBounds}
Fix constants $ k>2 $, $ 0<c<1 $, and $ s>kc^{-2} $. Then for any $ r\geq s $ the function $ \zeta_r(\sigma)= \sigma^{-1}\frac{2}{3}\zeta^{\frac{3}{2}}(\sigma^2r) $  satisfies the inequalities
\begin{align}
	&\zeta_r'(\sigma)\sim r\label{langerBound_1}\\
	&\abs{\zeta_r''(\sigma)} \les \frac{r}{\sigma}\\
	&\zeta_r''(\sigma)<0\label{langerBound_3}
\end{align}
uniformly for $ \sigma $ in the region $ [ks^{-\frac{1}{2}},c] $. Furthermore, for if $ s<r $ then
\begin{align}
	& r-s \les \zeta_r'(\sigma) - \zeta_s'(\sigma)\les r\label{langerDifferenceBound_0}\\
	&\abs{\zeta_r''(\sigma)-\zeta_s''(\sigma)} \les \frac{r-s}{\sigma}\label{langerDifferenceBound_1}\\
	&\zeta_r''(\sigma)-\zeta_s''(\sigma)<0\label{langerDifferenceBound_2}\,.
\end{align}
Here, all derivatives are with respect to $ \sigma $ and $ \les $ and $ \sim $ indicate bounds with respect to constants that depends only on $ k $ and $ c $ (i.e. not on $ r $ and $ s $).
\end{pr}
\begin{proof}
Recall that for $ z\geq 1 $ 
\begin{align*}
	\frac{2}{3}\zeta^{\frac{3}{2}}(z)=\int_{1}^{z} \sqrt{1-u^{-1}} \: du 
\end{align*}
or, explicitly,
\begin{align}\label{zetapower}
	\frac{2}{3}\zeta^{\frac{3}{2}}(z)=\sqrt{z(z-1)} -\log(\sqrt{z} +\sqrt{z-1} )\,.
\end{align}
One computes that
\begin{align*}
	\zeta_r'(\sigma)&=2r\sqrt{1-(\sigma^{2}r)^{-1}} - \frac{2}{3}\zeta^{\frac{3}{2}}(\sigma^{2}r)\sigma^{-2}\\
	&=r\sqrt{1-(\sigma^2r)^{-1}} +\sigma^{-2}\log(\sigma r^{\frac{1}{2}}+\sqrt{\sigma^{2}r-1})\,.
\end{align*}
Since $ \sigma\gtrsim s^{-\frac{1}{2}} $ on the regime in question, the log term is positive so the above is clearly greater than $ r\sqrt{1-k^{-2}} $. For the upper bound, one writes
\begin{align*}
	\log(\sigma r^{\frac{1}{2}}+\sqrt{\sigma^2r-1} )\leq \log(2\sigma r^{\frac{1}{2}})
\end{align*}
and then checks that as a function of $ \sigma $, $ \sigma^{-2}\log(\sigma a) $ has a global maximum of $ \frac{a^2}{2e} $. \par
To bound the second derivative, we first calculate
\begin{align*}
	\zeta_r''(\sigma)&=2\sigma^{-3}(1-(\sigma^2r)^{-1})^{-\frac{1}{2}}-2r\sigma^{-1}\sqrt{1-(\sigma^2r)^{-1}} +2\sigma^{-3}\frac{2}{3}\zeta^{\frac{3}{2}}(\sigma^2r)\\
	&=2\sigma^{-3}(1-(\sigma^2r)^{-1})^{-\frac{1}{2}}-2\sigma^{-3}\log(\sigma r^{\frac{1}{2}}+\sqrt{\sigma^2r-1} )
\end{align*}
and then observe that since $ \sigma^{-2}\les r $, 
\begin{align*}
	\abs{\zeta_r''}(\sigma) \les \sigma^{-1}(r+\sigma^{-2}\log(2\sigma r^{\frac{1}{2}})) \les \frac{r}{\sigma}\,.
\end{align*}
The negativity of $ \zeta_r'' $ follows from the above expression and the fact that $(1-(\sigma^2r)^{-1})^{-\frac{1}{2}}$ is a decreasing function of $ \sigma^2r $ while $ \log(\sigma r^{\frac{1}{2}}+\sqrt{\sigma^2r-1})  $ is an increasing function of $ \sigma^2r $ and one may verify that their difference is negative at, say, $ \sigma^2r=2 $.\par
For the estimate on the difference $ \zeta_r'-\zeta_s' $, we observe first that $ \zeta_r(\sigma)$ is an increasing function of $r$ for any fixed $\sigma$ and $\frac{ d \zeta^{\prime}_r}{dr}$   is uniformly bounded below and above for all allowed $ \sigma $ and $ s $. Hence, the mean value theorem implies \eqref{langerDifferenceBound_0} and \eqref{langerDifferenceBound_1}. To see  \eqref{langerDifferenceBound_2}, notice that $\zeta_r''$ is a negative decreasing function of $r$ for any fixed $\sigma$. 
\end{proof}
\begin{cor}\label{oscCor}
Let $ \sigma<c $ for $c>0$ sufficiently small. Then for all $ x\geq 1 $
\begin{align}\label{elarge}
	e(\sigma,r(x))=c_-(\sigma)q^{-\frac{1}{4}}\psi_+(\tau(x))+c_+(\sigma)q^{-\frac{1}{4}}\psi_-(\tau(x))\\
	c_-(\sigma)= \sigma^{-\frac{1}{6}}(1+e_3(\sigma)),\,\,c_+(\sigma)=\sigma^{-\frac{1}{6}}(1+e_4(\sigma)) \nn\,,
\end{align}
where 
\begin{align}\label{ejbound}
|e_j(\sigma)| \les \sigma ,\,\, |e_j^{\prime }(\sigma)| \les \sigma^{-\frac{1}{3}}, \,\, |e_j^{\prime \prime }(\sigma)| \les \sigma^{-\frac{7}{3}},\, \text{for}\,\ j=3,4.
\end{align}
Furthermore, for some constant $ C>0 $ large enough, we may write
\begin{align*}
	e(\sigma,r)=e^{i\zeta_r(\sigma)}a_+(\sigma,r)+e^{-i\zeta_r(\sigma)}a_-(\sigma,r)\,,
\end{align*}
where the functions $ a_\pm $ are smooth and satisfy the bounds
\begin{align*}
	&\abs{a_\pm(\sigma,r)} \les 1\\
	&\abs{\partial_{\sigma}[a_\pm(\sigma,r)]} \les \sigma^{-1} \les s^{\frac{1}{2}}\\
	&\abs{\partial^2_{\sigma}[a_\pm(\sigma,r)]} \les \sigma^{-\frac{7}{3}}\,,
\end{align*} 
uniformly for all $ 0<s\leq r $ and $ \sigma\in[Cs^{-\frac{1}{2}},c]$.
\end{cor}
\begin{proof}
		First, we connect the expression derived for $ e(\sigma,r) $ in Corollary \ref{AiryCor} to the basis constructed in Proposition \ref{oscBasisPr} at $ \zeta=0 $ (that is, $ x=1 $). Based on Corollary \ref{AiryCor}, by regarding $ r $ as a function of $ \zeta $, we have 	\begin{align*}
		&e(\sigma,r)|_{\zeta=0}=2 \sigma^{-\frac{1}{6}}\Ai(0)(1+e_1(\sigma))\\
		&\partial_{\zeta}\left[e(\sigma,r)\right]_{\zeta=0}=-2 \sigma^{-\frac{5}{6}}\Ai'(0)(1+O_2(\sigma^{\f23}))
	\end{align*} 
	since we may absorb the $ \Bi $ term in  into the $ e_1(\sigma) $ because the coefficient $ B(\sigma)=O(\sigma^{\infty}) $  and $ q(0)=1 $. Now, all of the Wronskians are evaluated at $ \zeta=0 $. It follows that
	\begin{align*}
		W\left[e(\sigma,r(x(\cdot))),\psi_\pm(\sigma,\cdot)\right]=\mp 2i \sigma^{-\frac{5}{6}}W[\Ai,\Bi](1+O(\sigma))=\mp\frac{i}{2\pi }\sigma^{-\frac{5}{6}}(1+l_1(\sigma))
	\end{align*}
	for some $l_1$ obeying the bounds in \eqref{ejbound} from Corolloary~\ref{AiryCor}. Similarly because
	\begin{align*}
		W[\psi_+(\sigma,\cdot),\psi_-\sigma,\cdot)]&=-\sigma^{-\frac{2}{3}}W[\Ai+i\Bi,\Ai-i\Bi](1+l_2(\sigma))=2i\sigma^{-\frac{2}{3}}W[\Ai,\Bi](1+l_2(\sigma))\\
		&=\frac{2i}{\pi}\sigma^{-\frac{2}{3}}(1+l_2(\sigma))\,,
	\end{align*}
	for some $l_2$ obeying the bounds in \eqref{ejbound}. From this, we see that the $ c_{\pm}(\sigma) $ have the form displayed in \eqref{elarge} with $e_3(\sigma)$ and $e_4(\sigma)$ satisfying the bounds in \eqref{ejbound}.
 
 Now, using the Airy function asymptotics \eqref{oscAiry} inside \eqref{elarge}, we may write
\begin{align*}
	e(\sigma,r)&= \frac{1}{2}(-Q)^{-\frac{1}{4}}(x)e^{-i(\zeta_r(\sigma)-\frac{\pi}{4})}(1+O\left(  (\zeta_r(\sigma))^{-1}\right)(1+\sigma b_+(\sigma))(1+e_3(\sigma))\\
	&+(-Q)^{-\frac{1}{4}}(x)e^{i(\zeta_r(\sigma)-\frac{\pi}{4})}(1+O\left(  (\zeta_r(\sigma))^{-1}\right)(1+\sigma b_-(\sigma))(1+e_4(\sigma))\\
	&=:e^{-i\zeta_r(\sigma)}a_-(\sigma,r)+e^{i\zeta_r(\sigma)}a_+(\sigma,r)\,,
\end{align*}
where each derivative of the $O(\cdot)$ terms loses one power. Thus, we are only left to show the bounds on $ a_\pm $ under the assumption that $ r\geq s>C \sigma^{-2} $. The bound $ \abs{a_\pm} \les 1 $ is immediate from the fact that for $ \sigma^2r>C $, $ \zeta_r(\sigma) $ is bounded below. For the derivatives, by making liberal use of the fact that $ \sigma^2r>C $, one computes from $Q(u)=u^{-1}-1$ that
\begin{align*}
	&\abs{(-Q(\sigma^2r))^{-\frac{1}{4}}} \les 1\\
	&\abs{\partial_{\sigma} [(-Q(\sigma^2r))^{-\frac{1}{4}}]} =\frac{1}{2}(1-(\sigma^2r)^{-1})^{-\frac{5}{4}}\sigma^{-3}r\les \sigma^{-1}\les s^{\frac{1}{2}}\\
	&\abs{\partial_{\sigma}^2[(-Q(\sigma^2r))^{-\frac{1}{4}}]} =\frac{5}{4}(1-(\sigma^2r)^{-1})^{-\frac{9}{4}}\sigma^{-6}r^{-2}+\frac{3}{2}(1-(\sigma^2r)^{-1})^{-\frac{5}{4}}\sigma^{-4}r^{-1}\les \sigma^{-2}\les s\,.
\end{align*}
By expanding \eqref{zetapower} in Taylor series, we also have $\zeta_r(\sigma) \sim \sigma r$ when $ \sigma^2 r >2$. Therefore, by \eqref{langerBound_1} we obtain 
\begin{align*}
	\abs{\partial_\sigma[(\zeta_r(\sigma))^{-1}]}=\abs{\zeta_r'(\sigma)(\zeta_r(\sigma))^{-2}} \les 1.
\end{align*}
Similarly
\begin{align*}
	\abs{\partial_{\sigma}^2[(\zeta_r(\sigma))^{-1}]} \les \abs{\zeta_r''(\sigma)(\zeta_r(\sigma))^{-2}} +\abs{(\zeta_r'(\sigma))^{2}(\zeta_r(\sigma))^{-3}} \les \sigma^{-1} \les s^{\frac{1}{2}}
\end{align*} 
again from Proposition \ref{langerBounds}. Now, one uses Proposition \ref{oscBasisPr} and Proposition \ref{langerBounds} to see that
\begin{align*}
	&\abs{b_{\pm}(\sigma,\zeta(\sigma^2r))} \les 1\\
	&\abs{\partial_\sigma[b_\pm(\sigma,\zeta(\sigma^2r))]}\les \abs{\sigma r\zeta'(\sigma^2r)\dot{b}_{\pm}(\sigma,\zeta(\sigma^2r))} +\partial_\sigma[b_{\pm}(\sigma,\zeta)]\vert_{\zeta=\zeta(\sigma^2r)}\\
	&\les (\sigma r)(\sigma^2r)^{-\frac{1}{3}}(\zeta(\sigma^2r))^{-\frac{5}{2}}+\sigma^{-1}\les \sigma^{-1}\les s^{\frac{1}{2}}\,,
\end{align*}
where $  \dot{ } $ represents the derivative with respect to the second variable and thus
\begin{align*}
	&\abs{1+\sigma b_{\pm}(\sigma,\zeta(\sigma^2r))} \les 1\,,\\
	&\abs{\partial_\sigma[\sigma b_\pm(\sigma,\zeta(\sigma^2r))]}\les \abs{b_\pm(\sigma,\zeta(\sigma^2)}+ \abs{\sigma^2r\zeta'(\sigma^2r)\dot{b}_{\pm}(\sigma,\zeta(\sigma^2r))} +\sigma\partial_\sigma[b_{\pm}(\sigma,\zeta)]\vert_{\zeta=\zeta(\sigma^2r)}\\
	&\les (\zeta(\sigma^2r))^{-\frac{3}{2}}+(\sigma^2r)(\sigma^2r)^{-\frac{1}{3}}(\zeta(\sigma^2r))^{-\frac{5}{2}}+(\zeta(\sigma^2r))^{-\frac{3}{2}}\les 1\,.
\end{align*}
 As per usual, suppressing the variable $\sigma^2r$ in $\zeta$, we obtain 
\begin{align*}
	\abs{\partial^2_\sigma[\sigma b_{\pm}(\sigma,\zeta)]}\les \abs{\partial_\sigma[b_{\pm}(\sigma,\zeta(\sigma^2 r))]} +\sigma\abs{\partial^2_\sigma [b_{\pm}(\sigma,\zeta(\sigma^2 r))]} \,.
\end{align*}
The first term is less than $ \sigma^{-1} $ by our previous computation and for the second we write 
\begin{align*}
	\partial^2_\sigma[b_\pm(\sigma,\zeta(\sigma^2 r)))]=\partial_{\sigma}^2[b_\pm(\sigma,\zeta(x))]_{x=\sigma^2 r} +(2r\zeta'(\sigma^2r)+4(\sigma r)^2\zeta''(\sigma^2r))\dot{b}_{\pm} (\sigma,\zeta) \\+2\sigma r\zeta'(\sigma^2r)\partial_{\sigma}[\dot{b}_{\pm}(\sigma,\zeta(x))]_{x=\sigma^2r}  +(2\sigma r\zeta'(\sigma^2r))^2\ddot{b}_\pm (\sigma,\zeta)
\end{align*}
and using various bounds from Proposition \ref{oscBasisPr} shows that  $\abs{\partial^2_\sigma[b_{\pm}(\sigma,\zeta(x))]}_{x=\sigma^2r} \les \sigma^{-3}$. 

This gives a total bound of
\begin{align*}	
	\abs{\partial^2_\sigma[b_{\pm}(\sigma,\zeta(\sigma^2 r))]} \les \sigma^{-2}\les s\,.
\end{align*}
\par
Summarizing all of these derivative computations, we have shown that
\begin{align*}
	&\abs{(-Q(\sigma^2r))^{-\frac{1}{4}}} \les 1,\,\, \abs{\partial_{\sigma} [(-Q(\sigma^2r))^{-\frac{1}{4}}]} \les \sigma^{-1} \les s^{\frac{1}{2}},\,\, \abs{\partial_{\sigma}^2[(-Q(\sigma^2r))^{-\frac{1}{4}}]} \les s\\
	&O((\zeta_r(\sigma))^{-1})\les 1,\,\, \abs{\partial_\sigma[(\zeta_r(\sigma))^{-1}]}\les 1,\,\,\abs{\partial_{\sigma}^2[(\zeta_r(\sigma))^{-1}]} \les \sigma^{-1} \les  s^{\frac{1}{2}}\\
	&\abs{1+\sigma b_\pm(\sigma, \zeta(\sigma^2r))} \les 1,\,\,\abs{\frac{d}{d\sigma}[\sigma b_{\pm}(\sigma, \zeta(\sigma^2 r))]} \les 1,\,\,\abs{\frac{d^2}{d\sigma^2}[\sigma b_{\pm}(\sigma,\zeta(\sigma^2 r))]} \les s\,.
\end{align*}
Combining these computations with the estimates on the derivatives of $e_i(\sigma)$ shows that
\begin{align*}
	\abs{\partial_{\sigma} a_{\pm}(\sigma,r)} \les \sigma^{-1} \les s^{\frac{1}{2}},\abs{\partial_{\sigma}^2a_\pm(\sigma,r)} \les \sigma^{-\frac{7}{3}}\,,
\end{align*}
the dominant term occurring when two derivatives fall on the $e_i(\sigma)$.
\end{proof}
\begin{rmk} \label{rmk:forstirling} We remark that  our connection in Corollary~\ref{oscCor} is consistent with the asymptotic behavior of $M_{\frac{i}{2 \sigma}, \f 12} ( 2 i \sigma r )$. Using \cite[(13.14.32) and (13.14.21)]{NIST} one can calculate  the following asymptotic behavior as $ \sigma r \to \infty$: 
\begin{align} \label{Mexp}
 M_{\frac{i}{2 \sigma}, \f 12} ( 2 i \sigma r )  \sim \frac{i e^{\frac{\pi}{4 \sigma}} }{ | \Gamma(1+\frac{i}{2 \sigma})|} \sin \Big( \sigma r  - \frac{1}{2 \sigma} \log (2 \sigma r)  + \theta(\sigma) \Big) \,,
\end{align}
where $\theta(\sigma):= \arg( \Gamma(1+i/(2\sigma))$. Therefore, by \eqref{erep} we have $e(\sigma, r  )\sim \sin \big( \sigma r   - \frac{1}{2 \sigma} \log (2 \sigma r  )  + \theta(\sigma)\big)$ as $ \sigma r \to \infty$. 
Here, we used the fact that   $ |\Gamma(i s)| = \sqrt{ \frac{\pi}{s \sinh(\pi s)}} $, see  (5.4.3) in \cite{NIST}.

On the other hand, by Stirling's formula \cite[(5.11.1)]{NIST}, we have that as $ \sigma \to 0 $, 
\begin{align}
\theta(\sigma) =  \frac{- \log (2 \sigma)}{2 \sigma}- \frac{1}{2 \sigma} + \frac{\pi}{4} - \frac{\sigma}{3} + O_2(\sigma^3) \nn 
\end{align}
and therefore, 
\begin{align*}
e^{\pm i( \sigma r - \frac{1}{2\sigma} \log(2 \sigma r) +\theta(\sigma))} = e^{\pm i( \sigma r - \frac{1}{2 \sigma} - \frac{ \log (4 \sigma^2 r )}{ 2 \sigma}+ \frac{\pi}{4})} ( 1+ O(\sigma)). 
\end{align*}


Furthermore, by \eqref{oscAiry}, one can see
\begin{align}
& c_{\mp}(\sigma) q^{-\frac{1}{4}}\psi_{\pm}(\tau(\sigma^2 r))=  \pi^{-\f 12} e^{\mp i( \sigma r - \frac{1}{2 \sigma} - \frac{ \log (4 \sigma^2 r )}{ 2 \sigma}- \frac{\pi}{4})} (1+ O( (\sigma r) ^{-1} \nn ))
\end{align} 
as $(\sigma r) \to \infty$.

\end{rmk}

We finish this section with the following lemma and its corollary. 
\begin{lemma}\label{lem:ab} Let $c$ be sufficiently small such that for all $\sigma <c$ \eqref{esmall} and \eqref{elarge} hold. Moreover let $ \f 3 4 \leq n <1$, $ m < \infty$ and define $ \chi_n^m := \tilde{\chi}_n \chi_m$. Then one has 
\begin{align}
|\partial_\sigma^j\{e(\sigma,  r )\}| \les   \sigma^{2-2j} r \,\,\ j=0,1,2\nn\,.
\end{align} 
in the support of $\chi_c(\sigma)\chi_n^m (\sigma^2 r)$. 
\end{lemma}
\begin{proof} Recall by \eqref{esmall}, and \eqref{elarge}, we have 
\begin{multline}\label{finalexp}
e(\sigma, r) = c_A(\sigma) \frac{\Ai(-\sigma^{-\f23} \zeta(\sigma^2 r)}{q^{\f14} (\sigma^2 r)} (1 + \sigma a(\sigma, \zeta(\sigma^2 r)) \\ + c
_B(\sigma) \frac{\Bi(-\sigma^{-\f23} \zeta(\sigma^2 r)}{q^{\f14} (\sigma^2 r)} (1 + \sigma b(\sigma, \zeta(\sigma^2 r))\,,
\end{multline}
where $a$ and $b$ obey the bounds in \eqref{ajbound} for $\sigma^2 r \leq 1$, and the bounds in \eqref{bpmbound} for $\sigma^2 r \geq 1$. Moreover, $C_A(\sigma) = \sigma^{-\f16}(1+e_3(\sigma))$; and $C_B(\sigma) =\sigma^{-\f16}(1+e_4(\sigma))$ for $\sigma^2 r>1$, $C_B=  \sigma^{\f56} e^{-\sigma^{-1} (\frac{\pi}{2} -1)}(1+e_2(\sigma))$ for $ \sigma^2 r \leq 1$. 

We show the statement first for the leading terms in \eqref{finalexp}. We note that by the definition of $\chi_c$, given $ n \geq \f 34 $, the cut-off $\tilde{\chi}_n \chi_m(x) $ is supported for $x \geq \f12$. Therefore, we consider \eqref{finalexp} for $ \zeta(n) \geq \zeta \geq \zeta(\frac{1}{2})>0$.

Let \ $\dot{}$ \ refer the $\zeta$-derivatives as per usual. By expansions \eqref{asymptoticAB},  \eqref{oscAiry} and the fact that $ -{\f 32} ( \zeta(\f12))^{\f 23} =\frac{\pi}{4} - {\f 12}$, we obtain the following inequalities 
\begin{align} 
&|\chi_c(\sigma) \chi_n^m(x) \Bi(-\sigma^{-\f 23} \zeta(x))| \les e^{\frac{1}{\sigma}( \frac{\pi}{4} - {\f 12})}  \chi(\sigma^2 r \leq 1) + \sigma^{\f16} \chi(\sigma^2 r \geq 1)\,,\label{Biders} \\ 
&|\chi_c(\sigma)\chi_n^m(x) \dot{ \Bi}(-\sigma^{-\f 23} \zeta(x))| \les e^{\frac{1}{\sigma}( \frac{\pi}{4} - {\f 12})} \chi(\sigma^2 r \leq 1) + \sigma^{-\f5 6} \chi(\sigma^2 r \geq 1)\,,\nn\\ 
&|\chi_c(\sigma)\chi_n^m(x) \ddot{ \Bi}(-\sigma^{-\f 23} \zeta(x))| \les e^{\frac{1}{\sigma}( \frac{\pi}{4} - {\f 12})}\chi(\sigma^2 r \leq 1) + \sigma^{-\f {11}{6}}\chi(\sigma^2 r \geq 1)\,. \nn  
\end{align}
 Similarly,
\begin{align}
&|\chi_c(\sigma)\chi_n^m(x)\Ai(-\sigma^{-\f 23} \zeta(x))| \les  \sigma^{\f16} \label{Aiders}\,, \\
&|\chi_c(\sigma)\chi_n^m(x)\dot{ \Ai}(-\sigma^{-\f 23} \zeta(x))| \les  \sigma^{-\f5 6}\,, \nn \\ 
&|  \chi_c(\sigma)\chi_n^m(x)\ddot{ \Ai}(-\sigma^{-\f 23} \zeta(x))| \les  \sigma^{-\f {11}{6}} \nn . 
\end{align}
Furthermore, by the definition of $q$, we have for $x\geq {\f12}$ 
\begin{align}
&\partial_x^j\{\zeta(x)\} = \la x \ra^{2/3-j} ,\,\,\  \partial_x^j\{ q^{-\f14}(x)\} = \la x \ra ^{1/6-j} \,.\label{qders}
\end{align}

Consequently, by utilizing \eqref{Biders}, \eqref{qders}, and incorporating $c_B(\sigma)$, we obtain 
\begin{align}\label{variable}
\Big| \chi_c(\sigma) \chi_n^m (\sigma^2 r) c_B(\sigma) \frac{\Bi(-\sigma^{-\f23} \zeta(\sigma^2 r))}{ q^{\f14}(\sigma^2 r)}\Big| \les  \chi_n^m (\sigma^2 r) [e^{ - \frac{1}{\sigma}( \frac{\pi}{4}-{\f 12})} \la \sigma^2 r \ra^{\f 16} +  \la \sigma^2 r \ra^{\f 16}]\les\chi_n^m\sigma^{2} r\,.
\end{align}
 In the last inequality we used the fact that $ n \leq \sigma^2 r $.
 
 To shorten the notation, for the rest of the proof we suppress the variable $\sigma^2 r$ in $\zeta$. We continue with estimating the $\sigma$ derivative of the leading term. We compute
\begin{multline}\label{middleBider}
| \chi_c \chi_n^m \partial_{\sigma} \{ c_B(\sigma) q^{-\f14}(\zeta) \Bi (-\sigma^{-\f23}\zeta)\} | \les \chi_c \chi_n^m  [ c_B^{\prime}(\sigma)(q(\sigma^2 r))^{-\f14} \Bi(-\sigma^{-\f23}\zeta) \\ + c_B(\sigma)\frac{\partial (q(\sigma^2 r))^{-\f14}}{\partial \sigma } \Bi (-\sigma^{-\f23}\zeta)+ c_B(\sigma)(q(\sigma^2 r))^{-\f14} \dot{ \Bi}(-\sigma^{-\f23}\zeta) \frac{d \zeta}{d \sigma} ] \,.
 \end{multline}
Using \eqref{Biders},\eqref{qders} and  $ \big| \frac{d \zeta(\sigma^2 r)}{d \sigma}\big| \les \sigma r$, we estimate $|\eqref{middleBider}| \les  r$.

Similarly, using  \eqref{Biders},\eqref{qders} and  $ \big| \frac{d \zeta(\sigma^2 r)}{d \sigma}\big| \les \sigma r$, one can compute
 \begin{align}\label{middleBisecder}
| \chi_c \chi_n^m  \partial^2_{\sigma} \{ c_B(\sigma) q^{-\f14}(\zeta) \Bi(-\sigma^{\f23}\zeta)\} | \les  \sigma^{-2} r.
\end{align}

As anticipated from the computation in \eqref{middleBider}, the restricting term $\sigma^{-2} r$ in \eqref{middleBisecder} arises  when both of the two derivatives fall on $\Bi(-\sigma^{-\f23}\zeta)$.  This situation leads us to the bound $ \sigma^{-2} (\sigma r)^2 $. However, due to the constraint $r \leq m \sigma^{-2}$, we can simplify this estimate to $ \sigma^{-2}r$.
With a similar observation and using \eqref{Aiders},  \eqref{qders} we  have
\begin{align}\label{middleAi}
 | \chi_c  \chi_n^m  \partial^j_{\sigma} \{ c_A(\sigma) q^{-\f14}(\zeta) \Ai(-\sigma^{\frac{3}{2}} \zeta)\} | \les \sigma^{2-2j} r\,,\quad j=0,1,2\,.
\end{align}
Hence, we obtain the bounds for the leading terms. 

We then proceed to estimate the error terms. We start with $ \sigma^2 r \leq 1$. For that we use Proposition~\ref{smallexp}. We conduct the essential computations only for $a(\sigma,\zeta)$ in \eqref{finalexp}, subsequent computations for $b(\sigma,\zeta)$ follows a similar procedure. Note that $|a_1|\les 1$, and 

\begin{align}\label{firstderer}
|\chi_c \chi_n^m \partial_{\sigma}\{a(\sigma, \zeta(\sigma^2r)\}|\les \chi_c \chi_n^m[ | \partial_\sigma\{a(\sigma, \zeta(x)\}|_{x=\sigma^2 r}| + |\dot a(\sigma, \zeta) \frac{d \zeta}{d \sigma}|] \les \sigma^{-\f43}\,.
\end{align}
Furthermore, 
\begin{align}\label{secondderer}
| \chi_c \chi_n^m\partial^2_{\sigma}\{a_1(\sigma, \zeta(\sigma^2r)\}|\les \chi_c \chi_n^m[ | \partial^2_\sigma\{a_1(\sigma, \zeta(x))\}|_{x=\sigma^2 r}| + |\partial_\sigma\{\dot a_1(\sigma, \zeta)\}\frac{d \zeta}{d \sigma}| \\ +| \ddot a_1(\sigma,\zeta) \big(\frac{d \zeta}{d \sigma}\big)^2| + |\dot a_1(\sigma, \zeta) \frac{d^2 \zeta}{d \sigma^2} |] \les \chi_n^m \sigma^{-\f 43} r  \les \sigma^{-\f {10}{3}} \:.\nn 
\end{align}

We finally, comment on the error term when $ \sigma^2 r \geq 1$. For that we use Proposition~\ref{oscBasisPr} and estimate $\sigma b_{\pm}$. In fact, using $\sigma b_{\pm}$ in \eqref{firstderer} and \eqref{secondderer} instead of $a_1$, we see the same bounds hold for the sigma derivatives of $\sigma b_{\pm}$ in the support of $ \chi_c \chi_n^m$. Combining these bounds with \eqref{variable}, \eqref{middleBider}, \eqref{middleBisecder}, \eqref{middleAi} we obtain the statement. Here for the discontinuity at $ \sigma^2 r=1$, we note that originally we converge to $M_{\frac{i}{2\sigma},{\f12}}(2i \sigma r)$ and this function is analytic in the vicinity of turning point. 
\end{proof}
The following corollary is due to Lemma~\ref{lem:bb} and Lemma~\ref{lem:ab}. 
\begin{cor}\label{cor:low} Fix $ c >0$ sufficiently small and $ k < \infty$. Then for any $ \beta \geq 0 $ we have 
\begin{align}
|\chi_c(\sigma)\chi_k (\sigma^2 r)  \partial_\sigma^j\{e(\sigma,  r )\}| \les \chi_c(\sigma) r [ \sigma^{\beta}\chi_{\f12}(\sigma^2 r) + \sigma^{2-2j}  \chi^k_{\f12}(\sigma^2 r)],\,\,\, j=0,1,2
\end{align} 
\end{cor}

\section{Eigenfunction approximation: Large Energies} 
In this section, we will consider the energies $\sigma \geq c >0 $ for some fixed $ c < 1$.  Recall that we have 
\begin{align*}
e(\sigma,  r) &=  -i \sigma^{-\f 12} [ e^{\frac {\pi}{\sigma}} -1]^{-\f 12}    M_{\frac{i}{2 \sigma}, \f 12} ( 2 i \sigma r)\,.  
\end{align*}
Using \cite[(13.14.2)]{NIST}, we may rewrite this as
\begin{align}\label{erep}
-i \sigma^{-\f 12} [ e^{\frac {\pi}{\sigma}} -1]^{-\f 12} e^{-i \sigma r} ( 2 i \sigma r) M\Big( 1 - i/(2 \sigma), 2, 2 i \sigma r \Big)\,,
\end{align}
where $M(a,b,z)$ is the Kummer function of the first kind. In this section we prove Proposition~\ref{prop:lexp}:
\begin{prop} \label{prop:lexp} Fix $k>0$. Then for $\sigma <c$ and $\sigma r \leq k$, we have that
\begin{align}\label{lsest}
     e(\sigma, r)= O_2(\sigma r).
\end{align}
Furthermore, for $\sigma r>k$, we have that
\begin{align}\label{llest}
  e(\sigma, r) =&  -i (\pi)^{-\f12} [ e^{ i(\sigma r - \frac{ \log(2\sigma r)}{2\sigma}+ \theta(\sigma))}  +e^{- i(\sigma r - \frac{ \log(2\sigma r)}{2\sigma}+\theta(\sigma))}]  + e^{i\sigma r}b_+(\sigma,r) + e^{-i\sigma r}b_-(\sigma,r)
\,
\end{align} 
where $\theta(\sigma)= \arg(\Gamma(1+i/(2 \sigma)))$,  and
\begin{align}
b_{\pm}(\sigma,r)=O_2(\sigma^0)\,.
\end{align}
\end{prop} 
Before we start the proof, we state a couple of expressions for $M(a,b,z)$. We will use \eqref{sumexp} to prove \eqref{lsest} and \eqref{intrep} to prove the \eqref{llest}. One has by \cite[(13.2.2)]{NIST}
\begin{align} \label{sumexp}
M(a,b,z) = 1+ \sum_{s=1}^{\infty} \frac{a(a+1)...(a+(s-1))}{ b(b+1)...(b+(s-1)) s!} z^s 
\end{align}
for any nonpositive integer $b$ and by \cite[(13.4.1)]{NIST}
\begin{align} \label{intrep}
M(a, b, z)= \frac{1}{ \Gamma(a) \Gamma(b-a)}\int\limits_0^{1} e^{z s } s^{a-1} (1-s)^{b-a-1}\: ds 
\end{align}
 for $\Re(b) > \Re(a) >0$. In the integral, the powers correspond to the principal branch of the logarithm.
 
 We start with the following lemma which analyzes the integral in \eqref{intrep} for $a= 1-\frac{i}{2\sigma}$, $b=2$ and $z= 2i\sigma r$. 
\begin{lemma} \label{lem:Mintest} For $\sigma>c$ and $\sigma r>k$, we have the expansion:
\begin{multline} \label{eq:Mintest}
\int\limits_0^{1} e^{2 i \sigma r s} s^{- \frac{i}{ 2\sigma}} (1-s)^{\frac{i}{2 \sigma}} ds \\ 
= (2i\sigma r)^{-1}\Bigg[  e^{2i \sigma r} (2i \sigma r )^{-\frac{i}{2\sigma}} \Gamma\Big(1+ \frac{i}{ 2\sigma}\Big) +(-2 i \sigma r )^{\frac{i}{2\sigma}}  \Gamma\Big(1- \frac{i}{ 2\sigma}\Big)  \Bigg]  \\  +  \Big[ \tilde{b}_+(\sigma,r) + e^{2 i \sigma r} \tilde{b}_-(\sigma,r) \Big] \,,
\end{multline}
 where $  \tilde{b}_{\pm}(\sigma,r)  = O_2(\sigma^{-1}r^{-1})$. 
\end{lemma}

\begin{proof}

Using the contour below, we obtain \eqref{contour}: 

\begin{minipage}{0.3\textwidth}
\begin{tikzpicture}
\draw[black,thick,->](0,0)-- (0,6.5) node[anchor=west]{$y$};
\draw[black,thick,->](0,0)-- (3.5,0) node[anchor=west]{$x$};
\draw[black,ultra thick] (0.5,0) arc
    [
        start angle=0,
        end angle=90,
        x radius=0.5cm,
        y radius =0.5cm
    ] ;
    \draw[black,thick,->] (0.5,0) arc
    [
        start angle=0,
        end angle=45,
        x radius=0.5cm,
        y radius =0.5cm
    ] ;
    \draw[black,ultra thick]  (3,0.5) arc
    [
        start angle=90,
        end angle=180,
        x radius=0.5cm,
        y radius =0.5cm
    ] ;
     \draw[black,thick,->]  (3,0.5) arc
    [
        start angle=90,
        end angle=145,
        x radius=0.5cm,
        y radius =0.5cm
    ] ;
 \draw[black, ultra thick] (0,0.5) -- (0,6);
  \draw[black,thick,->] (0,0.5) -- (0,3);
 \draw[black,thick,->] (3,6) -- (3,3);
 \draw[black, ultra thick] (3,0) -- (3,6);
  \draw[black,thick,->] (0,6) -- (1.5,6);
  \draw[black, ultra thick] (0,6) -- (3,6);
    \draw[dotted] (3,0)--(3,0.5);
  \node at (0.6,0.6) {$c_\varepsilon$};
    \node at (2.4,0.6) {$\tilde{c}_\varepsilon$};
     \node at (-.4,2) {$l_1$};
      \node at (3.4,2) {$l_2$};
       \node at (2,6.3) {$l_\mathcal{R}$};
       \node at (-.2,6) {$\mathcal{R}$};
         \node at (3,-0.2) {$1$};
         \node at(-.2,-.2){$0$};
         
\end{tikzpicture}
\end{minipage}
\begin{minipage}{0.7\textwidth}
\begin{multline}\label{contour}
\int_\varepsilon^{1-\varepsilon} e^{2 i \sigma r s} s^{- \frac{i}{ 2\sigma}} (1-s)^{\frac{i}{2 \sigma}} \:ds=\int\limits_{c_\varepsilon} e^{2 i \sigma r s} s^{- \frac{i}{ 2\sigma}} (1-s)^{\frac{i}{2 \sigma}} ds \\ +\int\limits_{l_1}e^{2 i \sigma r s} s^{- \frac{i}{ 2\sigma}} (1-s)^{\frac{i}{2 \sigma}} \:ds +\int\limits_{l_{\mathcal{R}}}e^{2 i \sigma r s} s^{- \frac{i}{ 2\sigma}} (1-s)^{\frac{i}{2 \sigma}} \:ds  \\
 +\int\limits_{l_2}e^{2 i \sigma r s} s^{- \frac{i}{ 2\sigma}} (1-s)^{\frac{i}{2 \sigma}} \: ds + \int\limits_{\tilde{c}_\varepsilon}e^{2 i \sigma r s}  s^{- \frac{i}{ 2\sigma}} (1-s)^{\frac{i}{2 \sigma}} \:ds\,.
 \end{multline}
\end{minipage}
Note that if $ \sigma \geq c >0$, then 
\begin{align*}
&\Big|\int\limits_{c_\varepsilon} e^{2 i \sigma r s} s^{- \frac{i}{ 2\sigma}} (1-s)^{\frac{i}{2 \sigma}} \:ds \Big|= \varepsilon  \Big|\int\limits_{0}^{\pi/2} e^{2i \sigma r \varepsilon e^{is}} ( \varepsilon e^{is})^{-\frac{i}{2\sigma}} (1- \varepsilon e^{is})^{\frac{i}{ 2\sigma}} e^{is}\:ds \Big| \leq 2 \varepsilon e^{\frac{\pi}{2 c}}\,, \\
&\Big|\int\limits_{\tilde{c}_\varepsilon} e^{2 i \sigma r s} s^{- \frac{i}{ 2\sigma}} (1-s)^{\frac{i}{2 \sigma}} \:ds \Big|= \varepsilon  \Big|\int\limits_{\pi/2}^{\pi} e^{2i \sigma r (1+\varepsilon e^{is})}  (1+ \varepsilon e^{is})^{-\frac{i}{ 2\sigma}} (- \varepsilon e^{is})^{\frac{i}{2\sigma}}e^{is}\:ds \Big| \leq 2\varepsilon e^{\frac{\pi}{2 c}}\,.
\end{align*}
Hence, the first and last term on the right side of the equality in  \eqref{contour} go to zero as $\varepsilon \to 0$. Moreover, as $ \sigma r \geq k$, we have that 
\begin{align*}
\Big|\int\limits_{l_{\mathcal{R}}}e^{2 i \sigma r s} s^{- \frac{i}{ 2\sigma}} (1-s)^{\frac{i}{2 \sigma}} \:ds \Big|= \Big|\int\limits_{0}^{1} e^{2i \sigma r (i \mathcal{R}+s)} (i\mathcal{R}+ s)^{-{\frac{i}{2 \sigma}}} (1-(i\mathcal{R}+s))^{\frac{i}{2 \sigma}} \: ds \Big| \leq e^{-2\mathcal{R}k} e^{\frac{\pi}{2 c}}\,,
\end{align*}
so that the third term on the right-hand side of also \eqref{contour} goes to zero as $\mathcal{R} \to \infty$. 

Therefore, letting $\varepsilon \to 0$ and $\mathcal{R} \to \infty$ we may write 
\begin{align} \label{a1a2sum}
\int\limits_0^{1} e^{2 i \sigma r s} s^{- \frac{i}{ 2\sigma}} (1-s)^{\frac{i}{2 \sigma}} ds=A_1+A_2\,,
\end{align} 
where 
\begin{align*}
&A_1(\sigma,r):=  \int\limits_0^{\infty} e^{2 i \sigma r (is)} (is)^{- \frac{i}{ 2\sigma}} (1-is)^{\frac{i}{2 \sigma}} i\:ds\\
&A_2(\sigma,r):=\int\limits_0^{\infty} e^{2 i \sigma r (1+i s)} (1+is)^{- \frac{i}{ 2\sigma}} (-is)^{\frac{i}{2 \sigma}} (-i)\:ds\,.
\end{align*}
We decompose
\begin{align}\label{eq:A1}
A_1(\sigma,r)&=\int\limits_0^{\infty} e^{-2 \sigma r s} (is)^{- \frac{i}{ 2\sigma}}  i\:ds +\int\limits_0^{\infty} e^{-2  \sigma rs } (is)^{- \frac{i}{ 2\sigma}} [(1-is)^{\frac{i}{2 \sigma}}-1] i\:ds  \\
&=(-2 i \sigma r )^{ -1 + \frac{i}{2\sigma}}  \Gamma(1-i/(2\sigma))+ \int\limits_0^{\infty}  e^{-2 \sigma r s}  (is)^{- \frac{i}{ 2\sigma}}[(1-is)^{\frac{i}{2 \sigma}}-1 ] i\:ds \,,\nn
\end{align}
and similarly 
\begin{align*}
A_2(\sigma, r) &= e^{2 i \sigma r}\int\limits_0^{\infty} e^{-2\sigma rs} (1+is)^{-\frac{i}{2\sigma}} (-is)^{\frac{i}{2\sigma}} (-i) \:ds \\ 
& =  e^{2 i \sigma r}\int\limits_0^{\infty} e^{-2\sigma rs}  (-is)^{\frac{i}{2\sigma}} (-i) ds  -e^{2 i \sigma r}\int\limits_0^{\infty} e^{-2\sigma rs} [(1+is)^{-\frac{i}{2\sigma}}-1] (-is)^{\frac{i}{2\sigma}} i\:ds \nn  \\
& = e^{2i \sigma r} (2 i\sigma r )^{-1-\frac{i}{2\sigma}} \Gamma(1+ i/(2\sigma)) - e^{2 i \sigma r}\int\limits_0^{\infty} e^{-2\sigma rs} [(1+is)^{-\frac{i}{2\sigma}}-1] (-is)^{\frac{i}{2\sigma}} i\:ds \nn \,.
\end{align*}
If we let
\begin{align*}
&\tilde{b}_+(\sigma,r):=  \int\limits_0^{\infty}  e^{-2 \sigma r s} (is)^{- \frac{i}{ 2\sigma}}[(1-is)^{\frac{i}{2 \sigma}}-1 ] i\:ds\,, \\
&\tilde{b}_-(\sigma,r):=- \int\limits_0^{\infty} e^{-2\sigma rs} [(1+is)^{-\frac{i}{2\sigma}}-1] (-is)^{\frac{i}{2\sigma}} i\:ds
\end{align*}
then comparing to \eqref{eq:Mintest}, we see that it is enough to demonstrate the claimed bounds on $\tilde{b}_{\pm}$. Below, we only show the bounds for $\tilde{b}_+$, as the bounds for $\tilde{b}_-$ follow similarly. 

 Observe that $|(is)^{- \frac{i}{ 2\sigma}}|$ and $|(1-is)^{\frac{i}{2 \sigma}}-1 |$ are each bounded in terms of $e^{\frac{\pi}{2\sigma}}$, which is $O(1)$ for $\sigma>c$. 
Therefore, we obtain 
\begin{align}\label{b1est}
| \tilde{b}_+(\sigma,r)| \les \int\limits_0^\infty e^{-2r\sigma s}\:ds=(2\sigma r)^{-1}\,.
\end{align}
We next estimate $\partial_{\sigma}b_+(\sigma,r)$. We compute
\begin{align*}
    \partial_{\sigma}[(is)^{-\frac{i}{2\sigma}}]=\frac{i}{2\sigma^2}(is)^{-\frac{i}{2\sigma}}\log(is),\,\,\partial_{\sigma}[(1-is)^{\frac{i}{2\sigma}}]=-\frac{i}{2\sigma^2}(1-is)^{\frac{i}{2\sigma}}\log(1-is),
\end{align*}
to obtain the bounds
\begin{align*}
    |\partial_{\sigma}[(is)^{-\frac{i}{2\sigma}}]|\les \sigma^{-2}\Span{\log{s}},\,\,|\partial_{\sigma}[(1-is)^{\frac{i}{2\sigma}}]|\les \sigma^{-2}\min\{s,\abs{\log(s)}\}\,.
\end{align*}
Using these bounds and the fact that
\begin{align*}
    | (1-is)^{\frac{i}{2 \sigma}} - 1|\les s\chi(s<1)+\chi(s\geq 1),
\end{align*}
allows us to write, for some $0<\alpha<1$,
\begin{align*}
    |\partial_\sigma[\tilde{b}_+(\sigma,r)]|&\les \max\{r,\sigma^{-2}\}\int\limits_0^\infty e^{-2\sigma rs}s\:ds+\sigma^{-2}\left(\int\limits_{0}^{1}e^{-2\sigma rs}s^{1-\alpha}\:ds+\int\limits_{1}^\infty e^{-2\sigma rs}s^{\alpha}\:ds\right)\\
    &\les \max\{r,\sigma^{-2}\}(\sigma r)^{-2}+\sigma^{-2}(\sigma r)^{\alpha-2}\int_0^\infty e^{-2s}s^{1-\alpha}\:ds+\sigma^{-2}(\sigma r)^{-\alpha-1}\int\limits_0^\infty e^{-2s}s^\alpha\:ds\,,
\end{align*}
which is $O(\sigma^{-2}r)$, as desired. Note that we are justified in differentiating under the integral defining $\tilde{b}_+$ because the resulting integral is absolutely convergent uniformly for $\sigma$ in some compact set away from $0$.

To complete the proof, we estimate the second derivative of $b_+(\sigma,r)$. First, it is easy to see that
\begin{align*}
    |\partial_\sigma^2[(is)^{-\frac{i}{2\sigma}}]|\les \sigma^{-4}\Span{\log(s)}^2+\sigma^{-3}\Span{\log(s)},\,\,|\partial_\sigma^2[(1-is)^{\frac{i}{2\sigma}}]|\les \sigma^{-3}s\,.
\end{align*}
Based on this and the previously derived bounds, one has that the six terms arising from taking the second $\sigma$-derivative of the integrand may be estimated as follows:
\begin{align*}
    &|\partial^2_\sigma[e^{-2\sigma rs}](is)^{-\frac{i}{2\sigma}}((1-is)^{\frac{i}{2\sigma}}-1)|\les r^2e^{-2\sigma rs}s^2\\
    &|e^{-2\sigma rs}\partial^2_\sigma[(is)^{-\frac{i}{2\sigma}}]((1-is)^{\frac{i}{2\sigma}}-1)|\les \sigma^{-3} e^{-2\sigma rs}(s^{1-\alpha}\chi(s<1)+s^\alpha\chi(s\geq 1))\\
    &|e^{-2\sigma rs}(is)^{-\frac{i}{2\sigma}}\partial_\sigma^2[((1-is)^{\frac{i}{2\sigma}}-1)]|\les \sigma^{-3} e^{-2\sigma rs}s\\
    &|\partial_\sigma[e^{-2\sigma rs}]\partial_\sigma[(is)^{-\frac{i}{2\sigma}}]((1-is)^{\frac{i}{2\sigma}}-1)|\les \sigma^{-2}re^{-2\sigma rs}(s^{2-\alpha}\chi(s<1)+s^{1+\alpha}\chi(s\geq 1))\\
    &|\partial_\sigma[e^{-2\sigma rs}](is)^{-\frac{i}{2\sigma}}\partial_\sigma[((1-is)^{\frac{i}{2\sigma}}-1)]|\les \sigma^{-2}re^{-2\sigma rs}s^2\\
    &|e^{-2\sigma rs}\partial_\sigma[(is)^{-\frac{i}{2\sigma}}]\partial_\sigma[((1-is)^{\frac{i}{2\sigma}}-1)]|\les \sigma^{-4}e^{-2\sigma rs}(s^{1-\alpha}\chi(s<1)+s^{\alpha}\chi(s\geq 1))\,.
\end{align*}
Every term may be integrated as before for a bound of $\sigma^{-3}r^{-1}$, so the proof is complete.
\end{proof}
\begin{proof}[Proof of Proposition~\ref{prop:lexp}]
We first prove \eqref{llest}. By \eqref{intrep} we have 
\begin{align*} 
 M(1-i/(2 \sigma), 2, 2 i \sigma r)= \frac{1}{ \Gamma(1+ i/(2 \sigma))\Gamma(1-i/(2\sigma))}\int_0^{1} e^{2 i \sigma r s} s^{- \frac{i}{ 2\sigma}} (1-s)^{\frac{i}{2 \sigma}} ds\,,
\end{align*}
which, on the support of $\tilde{\chi}_c(\sigma)\tilde{\chi}_{k}(\sigma r)$, when combined with \eqref{eq:Mintest} yields
\begin{align}\label{intrepsigma}
\begin{split}
    e(\sigma,r)=-i\sigma^{-\frac{1}{2}}[e^\frac{\pi}{\sigma}-1]^{-\frac{1}{2}}\left[e^{i\sigma r}\frac{(2i\sigma r)^{-\frac{i}{2\sigma}}}{\Gamma(1-i/(2\sigma))}+e^{-i\sigma r}\frac{(-2i\sigma r)^{\frac{i}{2\sigma}}}{\Gamma(1+i/(2\sigma))}\right]\\
    -i\sigma^{-\frac{1}{2}}[e^\frac{\pi}{\sigma}-1]^{-\frac{1}{2}}|\Gamma(1+i/2\sigma)|^{-2}(2i\sigma r)\left[e^{-i\sigma r}\tilde{b}_+(\sigma,r)+e^{i\sigma r}\tilde{b}_+(\sigma,r)\right],
    \end{split}
\end{align}
where we have used that $\overline{\Gamma(1+\frac{i}{2\sigma})}=\Gamma(1-\frac{i}{2\sigma})$.
Moreover, 
\begin{align*}
\Gamma(1\pm i/(2\sigma)) = |\Gamma(1+ i/(2\sigma))| e^{\pm i \theta(\sigma)}
\end{align*}
where $\theta(\sigma) = \arg (\Gamma(1+ i/(2\sigma))$ and $$(\mp 2i\sigma r)^{\pm \frac{i}{2\sigma} }= e^{\frac{\pi}{4\sigma}} e^{\pm i \frac{\log(2\sigma r)}{2\sigma}}.$$
One may also compute using \cite[(5.4.3)]{NIST} that
\begin{align*}
| \Gamma(1+ \frac{i}{2 \sigma})|^{-1} = (\pi)^{-\f 12}\sigma^{\f12}  [ e^{\frac{\pi}{\sigma}}-1]^{\f 12}e^{-\frac{\pi}{4\sigma}}.
\end{align*}
Therefore, inserting these equalities into \eqref{intrepsigma}, we have that
\begin{align*}
e(\sigma,  r) &=  -i(\pi)^{-\f12} [ e^{i(\sigma r - \frac{ \log(2\sigma r)}{2 \sigma} + \theta(\sigma))} + e^{-i(\sigma r - \frac{ \log(2\sigma r)}{2\sigma} + \theta(\sigma))} ] \\ 
&  +\frac{2}{\pi} r \sigma^{\frac{3}{2}} [1-e^{-\frac{\pi}{\sigma}}]^\frac{1}{2}\cdot [ e^{-i\sigma r} \tilde{b}_+(\sigma,r) + e^{i\sigma r} \tilde{b}_-(\sigma,r)] \nn\,.
\end{align*}
Letting $b_\pm (\sigma,r):=  \frac{2}{\pi} r\sigma^{\frac{3}{2}}  [1- e^{-\frac{\pi}{\sigma}}]^{\f 12}  \tilde{b}_{\mp}(\sigma,r)$, we see that $b_\pm (\sigma,r)$ satisfy the required bounds, in particular because $1-e^{-\pi/\sigma} = O_2(\sigma^{-1})$ for large $\sigma $.\par
To establish \eqref{lsest}, we write
\begin{align*}
	M(1-\frac{i}{2\sigma},2,2i\sigma r)=1+\sum_{s=1}^\infty C_s(\sigma)(2i\sigma r)^s,\,\,C_s(\sigma):=\frac{(1-\frac{i}{2\sigma})(2-\frac{i}{2\sigma})\cdots(s-\frac{i}{2\sigma})}{(s+1)!s!}.
\end{align*}
Using that
\begin{align*}
    \partial_\sigma[C_s(\sigma)]=\sigma^{-1}\sum_{j=1}^s\frac{(1-\frac{i}{2\cdot 1\sigma})(1-\frac{i}{2\cdot 2\sigma})\cdots (\frac{i}{2\cdot j\sigma})\cdots(1-\frac{i}{2\cdot s\sigma})}{(s+1)!},
\end{align*}
one easily checks that for $ \sigma>c $,  
\begin{align*}
	|\partial^j_{\sigma} [ C_s(\sigma)]|\leq \sigma^{-j}\frac{C^{s-j}}{(s-j)!},\,\,s\geq j,\,\,j=0,1,2,
\end{align*}
for some constant $ C>0 $. Clearly then for $ \sigma r<k $, we have that $ |M(1-\frac{i}{2\sigma},2,2i\sigma r)|\les 1 $ and the sum converges absolutely and uniformly for $ \sigma $ in any compact subset away from $ 0 $. Therefore, we can differentiate to find that
\begin{align*}
	|\partial_{\sigma}[M(1-\frac{i}{2\sigma},2,2i\sigma r)]|=|\sum_{s=1}^{\infty} C_s(\sigma)\frac{(2i\sigma r)^s}{\sigma}+C_s'(\sigma)(2i\sigma r)^{s}|\les \sigma^{-1}
\end{align*}
and similarly
\begin{align*}
	|\partial^{2}_{\sigma}[M(1-\frac{i}{2\sigma},2,2i\sigma r)]|\les \sigma^{-2}.
\end{align*}
Moreover, using that $ r\leq \sigma^{-1} $ in this regime, we have that $ e^{-i\sigma r }(2i\sigma r)=O_2(\sigma r) $. Finally, the term $ \sigma^{-\frac{1}{2}}[e^{\frac{\pi}{\sigma}}-1]^{-\frac{1}{2}} $ may be expanded as $(\pi)^{-\f12}+O_\infty(\sigma^{-1}) $ for large $\sigma$. Combining these computations yields \eqref{lsest}.

\end{proof}

 
 \section{Proof of Theorem~\ref{thm.main}}\label{proof}
 In this section, we estimate the kernel given by:
  $$ K_t(r,s) = \frac{1 }{2rs}\int\limits_0^{\infty} e^{it q^2 \sigma^2} e(q\sigma, r) e(q\sigma, s) \: d\sigma\,,$$ 
as  $ \sup_{r,s}  | K_t(r,s)| \les t^{-\f32}$ for $r,s \geq0$ and $t \geq 1$. Recall that this bound is sufficient to establish the validity of Theorem~\ref{thm.main} as one has
\begin{align*}
\| e^{itH_{0,q}} g \|_{L^{\infty}(\R^3)}= \Big\|\int\limits_0^{\infty} K_t(r,s) s^2 g(s) ds \Big\|_{L_r^{\infty}(\R^3)} \les \sup_{r,s}  | K_t(r,s)| \| s^2 g(s)\|_{ L^{1}_{\R^+}}. 
\end{align*}
  We normalize $q=1$ and chose $c <1$ sufficiently small in order to write
\begin{align*}
K_t(r,s)&= \frac{1}{rs}\int\limits_0^{\infty} e^{it \sigma^2} \chi_c(\sigma) e(\sigma,  r) e(\sigma,  s) \: d\sigma + \frac{1}{rs}\int\limits_0^{\infty} e^{it \sigma^2} \tilde{\chi}_c(\sigma) e(\sigma,  r) e(\sigma,  s) \: d\sigma \\
&=: K_t^l(r,s) +  K_t^h(r,s) \nn\,.
\end{align*}
We treat each of these integrals separately in the following subsections.
 \subsection{ Estimation of $K_t^l(r,s)$}
 In this subsection, we will prove that 
 \begin{prop}\label{prop:low} For any $c>0$ sufficiently small, we have that
\begin{align*}
 \sup_{r,s} |  K_t^l(r,s) | \les t^{-\f 32}\,.
\end{align*}
  \end{prop}
 
We prove Proposition~\ref{prop:low} with a series of lemmas. Fix some constant $ k \geq 4$ and write
  \begin{align} \label{maintermlow}
 K_t^l(r,s)= & \frac{1}{rs}\int_0^{\infty} e^{it \sigma^2} \chi_c(\sigma)\chi_{k}(\sigma^2 r)  \chi_{k}(\sigma^2 s) e(\sigma,r) e(\sigma,s) \: d\sigma  \\
 &+ \frac{1}{rs}\int_0^{\infty} e^{it \sigma^2} \chi_c(\sigma)  [\chi_{k}(\sigma^2 r) \tilde{\chi}_k(\sigma^2 s) + \tilde{\chi}_{k}(\sigma^2 r) \chi_k(\sigma^2 s)  ] e(\sigma,r) e(\sigma,s)\: d\sigma \nn \\
  &+ \frac{1}{rs}\int_0^{\infty} e^{it \sigma^2}  \chi_c(\sigma) \tilde{\chi}_k(\sigma^2 s) \tilde{\chi}_k(\sigma^2 r)e(\sigma,r) e(\sigma,s) \nn \\
  =&  K_1(r,s; t) + K_2(r,s; t) +K_3(r,s; t)  \nn \,.
 \end{align}
 By symmetry, we may always assume that $ r\geq s $. Furthermore, observe that the support of $ \chi_c(\sigma)\tilde{\chi}_k(\sigma^2r) $ is empty unless $ r\geq \frac{k}{c^2}>1 $ so we are free to assume that $ r\geq s >1 $ when considering  $ K_3 $.\par
 \begin{lemma}We have that 
$
| K_1(r,s;t)| \les t^{-\f 32}.
$
 \end{lemma}
 \begin{proof} Let $a(\sigma;s,r)= (rs)^{-1}\chi_c(\sigma)\chi_k(\sigma^2 r)  \chi_k(\sigma^2 s) e(\sigma,r) e(\sigma,s)$. With $^\prime$ denoting the derivative respect to $\sigma$, it is easy to see that, as a function of $\sigma$, $\chi_k(\sigma^2 s) = O_{\infty}(\sigma^{0})$ from the computation $  \chi^{\prime}_k(\sigma^2 r)  = \chi^{\prime}(\sigma^2 r) (2 \sigma r)$ and the fact that $\sigma^2 r \sim1$ on the support of $\chi^{\prime}(\sigma^2 r)$. Therefore, we may use the bounds from Corollary \ref{cor:low} to see that
\begin{align} \label{smsmain}
|\partial_\sigma^j\{ a(\sigma;r,s) \}| \les  \sigma^{4-2j}, \,\,\ j=0,1,2 \, .
\end{align}

Integrating by parts via the identity $ e^{it\sigma^2}=(2it\sigma)^{-1}\frac{d}{d\sigma}[e^{it\sigma^2}] $ and suppressing the variables $r$ and $s$, we obtain
\begin{align*}
K_1(r,s;t) = \frac{1}{2it}\int\limits_0^{\infty} e^{it \sigma^2} \Big(\frac{a(\sigma)}{\sigma}\Big)^{\prime} d \sigma =\frac{1}{2it}\int\limits_{\sigma < t^{-\f12}} e^{it \sigma^2} \Big(\frac{a(\sigma)}{\sigma}\Big)^{\prime} d \sigma +\frac{1}{2it}\int\limits_{\sigma \geq t^{-\f12}} e^{it \sigma^2} \Big(\frac{a(\sigma)}{\sigma}\Big)^{\prime} d \sigma \,.
\end{align*}
By  \eqref{smsmain}, we have $\Big|\Big(\frac{a(\sigma,r,s)}{\sigma}\Big)^{\prime}\Big| \les \sigma $, therefore, the first term is bounded by $t^{-\f 32}$. We apply another integration by parts to the second term to bound it by 
\begin{align} \label{nobound}
 t^{-2}\int\limits_{\sigma \geq t^{-\f12}} \Big| \frac{1}{\sigma} \Big(\frac{a(\sigma)}{\sigma}\Big)^{\prime \prime}\Big| + \Big| \frac{1}{\sigma^2}\Big(\frac{a(\sigma)}{\sigma}\Big)^{\prime} \Big| \: d\sigma \les t^{-2}\int\limits_{\sigma \geq t^{-\f12}} \sigma^{-2} \: d\sigma \les t^{-\f32} \,.
\end{align} 
Here, we omit the boundary term arising from the integration by parts since it is bounded by the integral in \eqref{nobound}.  This finishes the proof.

 \end{proof}

To estimate the other terms in \eqref{maintermlow}, we prove the following oscillatory integral estimate.

\begin{lemma} \label{lem:osc2} Suppose that for all $ r>\frac{k}{c^2} $, $ \omega_r(\sigma):[0,\infty)\rightarrow\bbR $ is a $ C^2 $ phase function and $\delta_r:\bbR\rightarrow \bbR$ is a weight function satisfying
\begin{enumerate}
	\item $ 0<\delta_r \les \omega_r'(\sigma)\les r $
	\item $ \omega_r''(\sigma)<0 $ and $ \abs{\omega_r''(\sigma)} \les \frac{\delta_r}{\sigma}$ 
 \end{enumerate}
and $ a_r(\sigma):[0,\infty)\rightarrow \bbC $ is a $ C^2 $ amplitude function satisfying 
\begin{enumerate} 
\item $\abs{a_r(\sigma)} \les \frac{\sigma^2}{r}$
\item $\abs{a_r'(\sigma)} \les \sigma^2$
\item $\int_{0}^{\infty} \sigma^{-1}(\abs{a_r''(\sigma)}+r\abs{a_r'(\sigma)})\chi(\sigma)  \: d\sigma \les 1$
\end{enumerate}
uniformly for $ \sigma\in [k^{\frac{1}{2}} r^{-\f12},c] $ and $ r>0 $, with $ \chi(\sigma) := \chi_c(\sigma) \tilde{\chi}_k(\sigma^2 r)$. Then we have that
\begin{align*}
 I^{\pm}(r,s;t) := \int\limits_0^{\infty} e^{i (t \sigma^2 \pm \omega_r(\sigma))} \chi(\sigma) a_r(\sigma) d\:\sigma = O( t^{-\f32})\,,
   \end{align*}
   with an implicit constant that does not depend on $ a_r $ or $ \omega_r $.

\end{lemma}
\begin{rmk}
	In the application of this lemma, the phase and amplitude may depend additionally on the parameter $ s $. The last sentence of the statement indicates that as long as the bounds on $ \omega_r $ and $ a_r $ hold uniformly in $ s $, the conclusion will also hold uniformly in $ s $. 
\end{rmk}
\begin{proof}
	As before, we first integrate by parts via $ e^{it\sigma^2}=(2it\sigma)^{-1}\frac{d}{d\sigma}[e^{it\sigma^2}] $ to find that
\begin{align*}
	I^{\pm}(r,s;t) &=\frac{1}{2it}\int\limits_{0}^{\infty} e^{i(t\sigma^2\pm\omega_r(\sigma))} [ b_{\pm} (\sigma;r) + \tilde{b} (\sigma;r) ]\chi(\sigma) \: d\sigma \\
		 &=I^{\pm}_1+I ^{\pm}_2
\end{align*}
for
\begin{align*}
b_{\pm} (\sigma;r) = \pm i\sigma^{-1} \omega_r^{\prime}(\sigma)   a_r(\sigma) \,\,\,\ \tilde{b} (\sigma;r) = \sigma^{-1} [\chi(\sigma) a_r(\sigma)]^{\prime}-\sigma^{-2} a_r(\sigma)]\,.
\end{align*}
We first estimate $ I_2^\pm $, splitting the integral as
\begin{align*}
	I^{\pm}_2 =  \frac{1}{2it} \int\limits_{0}^{t^{-\frac{1}{2}}} e^{i(t\sigma^2\pm \omega_r(\sigma)} \tilde{b}(\sigma;r)  \:d\sigma +  \frac{1}{2it}\int\limits_{t^{-\frac{1}{2}}}^\infty e^{i(t\sigma^2 \pm \omega_r(\sigma))} \tilde{b}(\sigma;r) \:d\sigma \,,
\end{align*}
and observe that the assumptions on $ a_r $ imply that
\begin{align*}
& |\tilde{b} (\sigma; r) |\les \sigma^{-1}  \abs{a^{\prime}_r(\sigma)} +\sigma^{-2} \abs{a_r(\sigma)}   \les  \sigma,  \\ 
&| \tilde{b}^{\prime}(\sigma; r)| \les \sigma^{-1} |a_r^{\prime \prime}(\sigma)| + \sigma^{-2}  |a_r^{\prime}(\sigma)| + \sigma^{-3} |a_r(\sigma)| \les \sigma^{-1} [ |a_r^{\prime \prime}(\sigma)|+ r |a_r^{\prime}(\sigma)|] \,,
\end{align*}
where for the final term in the second line we have used that $ \sigma^{-3}\abs{a_r(\sigma)} \les 1/(\sigma r)\les 1  $. Therefore, the first integral is bounded by $ t^{-\frac{3}{2}} $ and for the second we apply another integration by parts (ignoring the easily estimated boundary term) to bound it by 
 \begin{multline*}
 t^{-2}\int\limits_{t^{-\f12}}^{\infty}  \Big( | [\sigma^{-1}\tilde{b} (\sigma;r)]^{\prime} |+ |\sigma^{-1}\tilde{b}(\sigma;r) \omega^{\prime}_r(\sigma)| \Big)\chi(\sigma) \:d\sigma\\  \les t^{-\frac{3}{2}}\int\limits_{0}^\infty \sigma^{-1} (\abs{a_r''(\sigma)}+ r\abs{a_r'(\sigma)})\chi(\sigma)  \: d\sigma \les t^{-3/2}\,.
\end{multline*}

\par
We now turn our attention to $ I^{\pm}_1 $. Here, we must treat the $ \omega_r(\sigma) $ term as part of the phase so we write 
\begin{align*}
	I^{\pm}_1=(2it)^{-1}\int_{0}^{\infty} e^{it\Phi^{\pm}_{r,t}(\sigma)} b_{\pm}(\sigma;r) \chi(\sigma)\: d\sigma \,\,\ \text{with} \,\,\,\ 
	\Phi_{r,t}(\sigma):=\sigma^{2}\pm t^{-1} \omega_r(\sigma)\,.
\end{align*}
We have  $(\Phi^{\pm}_{r,t})'(\sigma)=2\sigma\pm t^{-1}\omega^{\prime}_r(\sigma)$, and $(\Phi^{\pm}_{r,t})''(\sigma)=2\pm t^{-1}\omega_r^{\prime \prime}(\sigma)$. As $\omega^{\prime}_r>0$ and $\omega_r^{\prime \prime}<0$, only $ \Phi_{r,t}^{-}$ has a stationary point and it is automatically non-degenerate.  \par
Since the phase in $I_1^+$ is non-stationary, the integral is easily estimated. As before, the integrand is $ O(\sigma) $ so we may split the domain of integration at $ t^{-\frac{1}{2}} $ and integrate by parts to find that
\begin{align}\label{ibp2}
	I_1^{+}\les t^{-\frac{3}{2}} +t^{-2}\int\limits_{t^{-\frac{1}{2}}}^{\infty}\left| \frac{ b^{\prime}_+(\sigma;r) }{(\Phi^{+}_{r,t}(\sigma))'} + b_+(\sigma;r)  \frac{d}{d\sigma}[((\Phi^{+}_{r,t})^{\prime})^{-1}(\sigma)] \right|\chi(\sigma) \: d\sigma \,.
\end{align}
Now by applying various properties of $ a_r $ and $ \omega_r $, we observe that
\begin{align}
& |b_{\pm}(\sigma;r) |\les \sigma^{-1} r a_r(\sigma) \les \sigma, \label{bbounds} \\ 
& | b_{\pm}^{\prime}(\sigma;r)| \les  \sigma^{-2} |[\omega_r^{\prime} a_r](\sigma) |+ \sigma^{-1} |[\omega_r^{\prime \prime }a_r](\sigma)|  + \sigma^{-1} |[\omega_r^{\prime} a_r^{\prime}](\sigma)| \les 1+ \sigma^{-1} r |a_r^{\prime}(\sigma)| \nn
\end{align}
so that, as $|(\Phi^{+}_{r,t})^{\prime}(\sigma)|^{-1} \les \sigma^{-1}$, we have 
\begin{align*}
	\left| \frac{ b^{\prime}_+(\sigma;r) }{(\Phi^{+}_{r,t}(\sigma))'}\right| \les \sigma^{-1}  + \sigma^{-2} r \abs{a_r^{\prime}(\sigma)} \,,
\end{align*}
which makes an admissible contribution. Furthermore, observe that 
\begin{align*}
 \frac{|(\Phi_{r,t}^{+})^{\prime \prime}(\sigma)|}{  |(\Phi_{r,t}^{+})^{\prime}(\sigma)|^2} \les \frac{2+\delta_r/(t\sigma)}{(2\sigma +\delta_r/t)^2}
\end{align*} 
so that
\begin{align*}
    \left|b_+(\sigma;r)\frac{d}{d\sigma}[(\Phi^+_{r,t})^{-1}](\sigma)\right|\les (2\sigma+\delta_r/t)^{-1}< \sigma^{-1}.
\end{align*}
Integrating now shows that $|I_1^+| \les t^{-\frac{3}{2}} $. \par
We now treat $I_1^{-}$. Since the stationary point is not explicitly calculable, some care is required. Because of the lower bound on $ \omega_r' $, we may find $ C $ depending only on $ c $ so that if $ t<\delta_r C $ then $ \Phi_{r,t}'<-1 $ uniformly on $ \supp \chi $. Using this, we break into cases depending on the value of $ t $:\\
 \underline{Case 1: $ t < \delta_r C $ }\\
 Due to the lower bound on $ \abs{\Phi_{r,t}'}  $, the phase is non-stationary and therefore the integral may be estimated similarly to $I_1^+$.\\
\underline{Case 2: $ t\geq \delta_r C $  }\\
In this regime, the phase may become stationary, however any stationary point will be non-degenerate because $( \Phi^{-}_{r,t})'\geq 2 $ on $ \text{supp} \chi $ uniformly in $ r $ by the properties of $\omega_r$.  Indeed,  because the second derivative is bounded below away from $ 0 $, we claim that we may always find some $ \sigma_*\in \text{supp} \chi $ so that $ |(\Phi^-_{r,t})'(\sigma)| \geq 2 |\sigma-\sigma_*| $ on $ \text{supp} \chi $. If $(\Phi^-_{r,t})' $ vanishes at some $ \sigma_* $ then this is immediate from the mean value theorem. Otherwise, we know that $ \Phi^-_{r,t} $ is increasing on $[a,b]=\text{supp}\chi$ so we must have that either $ (\Phi^-_{r,t})'(a)>0 $ or $ (\Phi^-_{r,t})'(b)<0 $ if $ \Phi^-_{r,t} $ does not vanish. In either case, the claim is easily seen to hold with $ \sigma_*=a $ or $ b $, respectively.\par
Splitting $ I_1^- $ as above, we write 
\begin{align*}
	I^-_1&=\frac{1}{2it}\int\limits_{\abs{\sigma-\sigma^*} <t^{-\frac{1}{2}}} e^{it\Phi^-_{r,t}(\sigma)} b_-(\sigma;r)\chi(\sigma) \: d\sigma +\frac{1}{2it}\int\limits_{\abs{\sigma-\sigma_*} >t^{-\frac{1}{2}}}^{\infty} e^{it\Phi^-_{r,t}(\sigma)} b_-(\sigma;r)\chi(\sigma)\: d\sigma\,.
\end{align*}
As before, the integrand of the first integral is bounded, so that by integrating by parts in the second, we see that we need only estimate 
\begin{align*}
 t^{-2}\int\limits_{\abs{\sigma-\sigma_*} >t^{-\frac{1}{2}}}^{\infty }  \Bigg( \frac{ b_-^{\prime}(\sigma;r) }{(\Phi^-_{r,t})^{\prime}(\sigma)} - \frac{b_-(\sigma;r)(\Phi^-_{r,t})''(\sigma)}{ ((\Phi^-_{r,t})'(\sigma))^2 }  +\frac{b_-(\sigma;r)}{(\Phi^-_{r,t})'(\sigma)} \chi'(\sigma)\Bigg) \,d \sigma .
\end{align*}
The term with $ \chi' $ is easily seen to be admissible, whereas the rest of the integral can be bounded as
\begin{align*}
	& t^{-2}\int\limits_{\abs{\sigma-\sigma_*}>t^{-\frac{1}{2}} }^{} \left(\frac{|b_-^{\prime}(\sigma;r) |}{\abs{\sigma-\sigma_*} }+ \frac{|b_- (\sigma;r)| \sup\abs{(\Phi^{-}_{r,t})''(\sigma)} }{\abs{\sigma-\sigma_*}^2}\right)\chi(\sigma)\: d\sigma \\
	&\les t^{-\frac{3}{2}}\int\limits_0^{\infty}    \abs{b_-'(\sigma;r) } \chi(\sigma) \:d\sigma  + t^{-2}\int_{\abs{\sigma-\sigma_*}>t^{-\frac{1}{2}} }^{} \frac{ |b_-(\sigma;r)| \sup| \Phi_{r,t}''(\sigma)| }{ |\sigma - \sigma_*|^2} d \sigma\,.
 \end{align*}
The bounds in $\eqref{bbounds}$ show that the first integral is bounded by $ t^{-\frac{3}{2}} $ whereas for the second we use that
\begin{align*}
	\abs{b_-(\sigma;r)\Phi_{r,t}''(\sigma)} \les t^{-1}\sigma\abs{\omega_r''(\sigma)} \les \frac{\delta_r}{t}\les 1
\end{align*}
to conclude. This finishes the proof.
\end{proof}

 We are now ready to show that
 \begin{lemma} We have that \label{lem:sml}$
 | K_2(r,s;t)| \les  t^{-\f32}
$. 
 \end{lemma}
 \begin{proof}
	 Since $ \tilde{\chi}_k(\sigma^2r)\chi_k(\sigma^2s) $ only has non-empty support if $ r\gtrsim s $, it suffices to consider 
	\begin{align*}
		\int_{0}^{\infty} e^{it\sigma^2}\frac{e(\sigma,r)}{r}\frac{e(\sigma,s)\chi_k(\sigma^2s)}{s}\chi(\sigma)\: d\sigma 
	\end{align*}
	where $\chi(\sigma)=\chi_c(\sigma)\tilde{\chi}_k(\sigma^2r) $  as in Lemma \ref{lem:osc2}. By Corollary~\ref{oscCor}, on the support of $ \chi $, $ e(\sigma,r)= e^{i\zeta_r(\sigma)}a_{+}(\sigma,r)+e^{- i\zeta_r(\sigma)}a_{-}(\sigma,r) $.  Therefore, we need to estimate the integrals
 \begin{align} \label{smalllarge}
	 I_{\pm}:=\int_{0}^{\infty} e^{i(t\sigma^2\pm\zeta_r)}\frac{a_{\pm}(\sigma,r)}{r}\frac{b_s(\sigma)}{s}\chi(\sigma)\: d\sigma\,,
 \end{align} 
 where $ b_s(\sigma):=\chi_{k}(\sigma^2s)e(\sigma,s) $.\par
 We verify the conditions of Lemma~\ref{lem:osc2}. Proposition \ref{langerBounds} shows that $ \zeta_r $ satisfies the hypotheses of the lemma so we need only check that $ a_{r,s}(\sigma)=\frac{a_{-}(\sigma)}{r}\frac{b_s(\sigma)}{s} $ satisfies the hypotheses as well, uniformly in $ s $. From Corollary \ref{oscCor}, we obtain
 \begin{align*}
 	\abs{a_-(\sigma,r)} \les 1,\,\,\abs{a_-'(\sigma,r)} \les \sigma^{-1},\,\,\abs{a_-''(\sigma,r)} \les \sigma^{-\frac{7}{3}}
 \end{align*}
 while from Corollary \ref{cor:low}, it follows that
 \begin{align*}
	 \abs{b_s(\sigma)} \les s\sigma^2,\,\,\abs{b_s'(\sigma)} \les s[\sigma^2\chi_{\frac{1}{2}}(\sigma^2s)+\chi_{\frac{1}{2},k}(\sigma^2s)],\,\,\abs{b_s''(\sigma)} \les s[\sigma^2\chi_{\frac{1}{2}}(\sigma^2s)+\sigma^{-2}\chi_{\frac{1}{2},k}(\sigma^2s)],\,\,
 \end{align*}
 where $ \chi_{\frac{1}{2},k}=(\tilde{\chi}_{\frac{1}{2}}\cdot\chi_{k})(\sigma^2s) $. Clearly $ \abs{a_{r,s}(\sigma)} \les \frac{\sigma^2}{r} $ and furthermore
\begin{align*}
	\abs{a_{r,s}'(\sigma)} \les r^{-1}[\sigma^2\chi_{\frac{1}{2}}(\sigma^2s)+\chi_{\frac{1}{2},k}(\sigma^2s)]+\frac{\sigma}{r}\les\sigma^2
\end{align*}
since $ r^{-1}\les \sigma^2 $ on the support of $ \chi $. Moreover, one can also check that
\begin{align*}
	\abs{a_{r,s}''(\sigma)} \les \sigma\chi_{\frac{1}{2}}(\sigma^2s)+\chi_{\frac{1}{2},k}(\sigma^2s)\,.
\end{align*}
It now follows from the computation
\begin{align*}
	\int_{0}^{\infty}\sigma^{-1}\chi_{\frac{1}{2},k}(\sigma^2s)\chi \: d\sigma\leq \log(2k)/2
\end{align*}
that $\int\limits_{0}^{\infty}\sigma^{-1}(\abs{a_{r,s}''(\sigma)} +r\abs{a_{r,s}'(\sigma)})  \: d\sigma\les 1  $, so we conclude from Lemma \ref{lem:osc2} that $ I_1=O(t^{-\frac{3}{2}}) $.
\end{proof}

We next prove that 
 \begin{lemma} $
 |K_3(r,s;t)| \les t^{-\f32}.  
$
 \end{lemma}
 \begin{proof}
  By Corollary~\ref{oscCor}, in this $ \sigma $ regime, $e(\sigma,r) e(\sigma,s)$ can be written as a sum of the terms 
  \begin{align}
   e^{ - i (\zeta_r(\sigma)\pm \zeta_s(\sigma))} a_{- }(\sigma,r) a_{\pm }(\sigma,s), \,\,\,\  e^{ i (\zeta_r(\sigma)\pm \zeta_s(\sigma))} a_{+}(\sigma,r) a_{\pm}(\sigma,s)  
  \end{align}
Therefore, it suffices to bound
 \begin{align} \label{largelarge}
 &\int\limits_0^{\infty} e^{it \sigma^2 \pm i (\zeta_r(\sigma)- \zeta_s(\sigma))}  \frac{a_{\pm }(\sigma,r)}{r} \frac{a_{\mp }(\sigma,s) \tilde{\chi}_k(\sigma^2 s)}{s}\chi(\sigma)\:d \sigma\,, \\
 &\int_0^{\infty} e^{it \sigma^2 \pm i (\zeta_r(\sigma)+ \zeta_s(\sigma))} \frac{a_{\pm}(\sigma,r)}{r} \frac{a_{\pm }(\sigma,s) \tilde{\chi}_k(\sigma^2 s)}{s}\chi(\sigma)\:d \sigma\,,
\end{align}    
where $ \chi(\sigma):=\chi_c(\sigma) \tilde{\chi}_k(\sigma^2 r) $. In order to apply Lemma \ref{lem:osc2}, we must verify that the phases $ \zeta_r+\zeta_s $ and $ \zeta_r-\zeta_s $ and the amplitude $a_{r,s}(\sigma)= \frac{a_+(\sigma,s)}{r}\frac{a_+(\sigma,r)\tilde{\chi}_k(\sigma^2s)}{s} $ satisfies the hypotheses of the lemma (the latter being sufficient because $ a_+ $ and $ a_- $ obey the same bounds). When $r=s$, the phase $\zeta_r-\zeta_s$ degenerates to $0$. Ignoring this easily treated case, the conditions on the phases are satisfied by Proposition \ref{langerBounds} so we consider the amplitude $ a_{r,s}(\sigma) $. From Corollary \ref{oscCor}, we see that
\begin{align*}
\abs{a_{r,s}(\sigma)} \les \frac{\tilde{\chi}_k(\sigma^2s)}{rs}\les \frac{\sigma^2}{r}
\end{align*}
and also that 
\begin{align*}
	\abs{a^{\prime}_{r,s}(\sigma)} \les \frac{\sigma^{-1}\tilde{\chi}_k(\sigma^2s)}{rs}\les \frac{1}{rs^{\frac{1}{2}}}\,,
\end{align*}
which is indeed less than $ \sigma^2 $ on the domain in question. Furthermore, we have that
\begin{align*}
	\abs{a''_{r,s}(\sigma)} \les \frac{\sigma^{-\frac{1}{3}}}{r}\,,
\end{align*}
from which one can easily check that $\int\limits_{0}^{\infty} \sigma^{-1}(\abs{a''_{r,s}(\sigma)} +r\abs{a'_{r,s}(\sigma)}) \chi(\sigma)\: d\sigma \les 1 $.\par
Applying Lemma \ref{lem:osc2} now completes the proof.
 \end{proof}

 \begin{proof}[Proof of Proposition~\ref{prop:low}]
 Combining the bounds for $K_1$, $K_2$ and $K_3$ we obtain the statement. 
 \end{proof}
 
\subsection{Estimate of $K_t^h(r,s)$}
We will prove the following Proposition. 

 \begin{prop} \label{prop:high}Let $c<1$ and  $ t \geq 1$.  Then we have 
$
 \sup_{r,s} |  K_t^h(r,s) | \les t^{-\f 32}.
$
  \end{prop}

  As in the previous subsection, we will prove Proposition~\ref{prop:high} with a series of lemmas. We let $k \geq 4 $  and write  
    \begin{align} \label{maintermhigh}
  K_t^h(r,s)= & \frac{1}{rs}\int\limits_0^{\infty} e^{it \sigma^2} \tilde{\chi}_c(\sigma)\chi_{k}(\sigma r)  \chi_{k}(\sigma s) e(\sigma,r) e(\sigma,s)\: d\sigma  \\
 &+ \frac{1}{rs}\int\limits_0^{\infty} e^{it \sigma^2} \tilde{\chi}_c(\sigma)  [\chi_{k}(\sigma r) \tilde{\chi}_k(\sigma s) + \tilde{\chi}_{k}(\sigma r) \chi_k(\sigma s)  ] e(\sigma,r) e(\sigma,s)\: d\sigma \nn \\
  &+ \frac{1}{rs}\int\limits_0^{\infty} e^{it \sigma^2}  \tilde{\chi}_c(\sigma) \tilde{\chi}_k(\sigma r) \tilde{\chi}_k(\sigma s)e(\sigma,r) e(\sigma,s) \: d\sigma\nn \\ 
  =:& \tilde{K}_1(r,s;t)+\tilde{K}_2(r,s;t)+ \tilde{K}_3(r,s;t). \nn 
 \end{align}
We start by estimating the first term, in which the integrand does not oscillate. 
\begin{lemma}\label{lem:tildeK1}
$
| \tilde{K}_1(r,s;t)| \les  t ^{-\f 32}.
$
 \end{lemma}
\begin{proof} Let $a(\sigma;r,s)= (rs)^{-1} \tilde{\chi}_c(\sigma) \chi_{k}(\sigma r)  \chi_{k}(\sigma s) e(\sigma,r) e(\sigma,s)$, be the integrand of $\tilde{K}_1(r,s;t)$ so that using \eqref{lsest} and the fact that $\chi_k(\sigma r)$ and $\chi_k(\sigma s)$ are $O_\infty(\sigma^0)$, we have that
\begin{align*}
a(\sigma;r,s)=O_2(\sigma^2)\,.
\end{align*}
Hence, by integrating by parts twice, we obtain
\begin{align} \label{ibptwice}
 \Big| \int\limits_0^{\infty} e^{i t \sigma^2} a(\sigma;r,s) d \sigma\Big| \les t^{-2}\int\limits_0^{\infty}   \Big| \Big(\frac{1}{\sigma} \Big( \frac{ a(\sigma;r,s)}{\sigma} \Big)^{\prime} \Big)^{\prime} \Big|\: d \sigma  
 \les t^{-2}\int\limits_{c}^{\infty} \sigma^{-2} \:d\sigma \les t^{-2}\,.
 \end{align}
 We obtain no boundary terms since the support of the integral is away from both zero and infinity. 
\end{proof}

We next prove the following oscillatory integral lemma which will be useful for estimating the remaining terms. It's worth noting that owing to the dependency on $r$, Lemma~\ref{lem:osc3} differs from the Van der Corput lemma.
\begin{lemma}\label{lem:osc3} Let $a_s(\sigma)=  \tilde{\chi}_c(\sigma) O_2( \sigma s )$, $\varphi_r(\sigma):\bbR\rightarrow\bbR$. Then for $t>1$, one has  
\begin{align*}
I(r,s;t)= \frac{1}{rs} \int\limits_0^{\infty} e^{it \sigma^2 \pm i \varphi_{r}(\sigma)} \tilde{\chi}_k(\sigma r) a_s(\sigma) \: d \sigma =O(t^{-\f32}),  
\end{align*}
provided that the following conditions hold for $\varphi_{r}(\sigma)$ within the support of the integral:
\begin{align}\label{condu}
0\leq \varphi^{\prime}_{r}(\sigma) \les r, \quad \varphi^{\prime \prime}_{r}(\sigma) \leq 0, \quad |\varphi^{\prime \prime}_{r}(\sigma)| \les \sigma^{-2} r\,.
\end{align}
Furthermore, if $\varphi_r(\sigma)$ depends upon an additional parameter and \eqref{condu} holds uniformly in this parameter then $|I(r,s;t)|\les t^{-\f23}$ uniformly as well. The same is true if $a_s(\sigma)=O_2(\sigma s)$ uniformly in some additional parameter.
\end{lemma}

\begin{proof}
As in the proof of Lemma~\ref{lem:osc2}, we start with an integration by parts and write 
\begin{align*}
I^{\pm}(r,s;t) &= \frac{1}{2 it rs}\int\limits_0^{\infty} e^{it \sigma^2 \pm i \varphi_{r}(\sigma)}[ b_1(\sigma;r,s) + b_2(\sigma;r,s)] \:d \sigma \\
&=: I^{\pm}_1 + I^{\pm}_2 \,,\nn 
\end{align*}
with
\begin{align*}
b_1(\sigma;r,s):= \sigma^{-1} ( \tilde{\chi}_k(\sigma r) a_s(\sigma))^{\prime} - \sigma^{-2}  \tilde{\chi}_k(\sigma r) a_s(\sigma), \,\,\,\,\,\ 
b_2(\sigma;r,s):= \pm i  \sigma^{-1} \tilde{\chi}_k(\sigma r) a_s(\sigma) \varphi_{r}^{\prime}(\sigma)\,.
\end{align*}
One may check that $|b_1(\sigma;r,s)|\les \sigma^{-1}s$ using that $\tilde{\chi}_k(\sigma r)=O_{\infty}
(\sigma^{0})$. Moreover,  $\tilde{\chi}_k(\sigma r)\les \sigma  r$, shows that $|b_1'(\sigma;r,s)|\les \sigma^{-1}rs$. Thus,
we can apply another integration by parts to $I^{\pm}_1$ to bound it by 
\begin{align*} 
\frac{1}{t^2 rs} \int\limits_c^{\infty} | (\sigma^{-1} b_1(\sigma;r,s))^{\prime}|+ \sigma^{-1}   | b_1(\sigma;r,s) \varphi_{r}^{\prime}(\sigma)| \: d \sigma \,,
\end{align*}
which is bounded by $t^{-2}$ since $ \varphi_{r}^{\prime}(\sigma) \les r$.

We next focus on $I^{\pm}_2$. We let $ \Phi_{\pm}(\sigma) = \sigma^2 \pm t^{-1} \varphi_{r} (\sigma) $, which implicitly depends on $r$ via $\varphi_r(\sigma)$.  Note that the conditions on $\varphi_{r}$ are arranged so that only $\Phi_{-}(\sigma)$ may have a critical point. However, this critical point can only be non-degenerate because of the condition $\varphi^{\prime \prime}_{r}(\sigma) \leq 0$. As in the proof of Lemma~\ref{lem:osc2}, because $\Phi_-''(\sigma)\geq 2$, we may find some $\sigma_*>0$ such that  $| \Phi_-^{\prime}(\sigma) | \gtrsim |\sigma - \sigma^*|$ on the support of $\tilde{\chi}_c(\sigma)$.

With the above in mind, we first focus on $I^{-}_2 $ and write:
\begin{align*}
I^{-}_2 = (2itrs)^{-1}\int\limits_{|\sigma - \sigma_*| \leq t^{-\f12}} e^{it \Phi_{-}(\sigma)} b_2(\sigma; r,s)\: d \sigma + (2itrs)^{-1}\int\limits_{|\sigma - \sigma_*| > t^{-\f12}} e^{it \Phi_{-}(\sigma)} b_2(\sigma; r,s) \:d \sigma. 
\end{align*}
Since $|b_2(\sigma;r,s)| \les  r s$, the first term in $I^{-}_2$ is bounded by $t^{-\f32}$. We apply another integration by parts to bound the second term in $I^{-}_2$ by 
\begin{align} \label{lastterm}
\frac{1}{t^2 rs} \int_{|\sigma - \sigma^*| > t^{-\f12}} \frac{|b_2^{\prime}(\sigma;r,s)|}{| \Phi_-^{\prime}(\sigma)|} + \frac{|b_2(\sigma;r,s)|| \Phi_-^{\prime \prime }(\sigma)|}{| \Phi_-^{\prime}(\sigma)|^2} \: d \sigma \,,
\end{align}
where we omit the boundary term since it will be no worse than \eqref{lastterm}. As for $b_1(\sigma;r,s)$, one may argue that $|b_2^{\prime}(\sigma;r,s)| \les \sigma^{-1} r s$. Therefore,
$$
\frac{|b_2^{\prime}(\sigma;r,s)|}{| \Phi^{\prime}(\sigma)|} \les \frac{rs} {\sigma |\sigma - \sigma_*|} \les \frac{r s}{\sigma^{2}} + \frac{rs}{|\sigma - \sigma_*|^2}\,,
$$ 
which yields the bound:
\begin{align}\label{finest}
\int_{|\sigma - \sigma^*| > t^{-\f12}} \frac{|b_2^{\prime}(\sigma;r,s)|}{| \Phi_-^{\prime}(\sigma)|} \:d \sigma \les \int_0^{\infty} \tilde{\chi}_c(\sigma) \frac{rs}{\sigma^{2}} d \sigma + \int_{|\sigma - \sigma^*| > t^{-\f12}}  \frac{r s}{|\sigma - \sigma_*|^2} d \sigma \les r s\: \la t \ra ^{\f12}. 
\end{align}

To estimate the second term in \eqref{lastterm}, we recall that $|b_2(\sigma;r,s)| \les s| \varphi^{\prime}_{r}(\sigma)| \les rs$. Furthermore, as  $ |\Phi_-^{\prime \prime }(\sigma)| \leq 2 + t^{-1} |\varphi_{r}^{\prime \prime} |$ and  $| \varphi_{r}^{\prime \prime}|  \les \sigma^{-2} r $, we have 
$$
\frac{|b_2(\sigma;r,s)|| \Phi_-^{\prime \prime }(\sigma)|}{| \Phi_-^{\prime}(\sigma)|^2} \les \frac{ rs}{| \Phi_-^{\prime}(\sigma)|^2} + \frac{| \varphi^{\prime}_{r}(\sigma)| r s}{ t \sigma^2| \Phi_-^{\prime}(\sigma)|^2} \les \frac{ r s}{|\sigma - \sigma_*|^2}  + \frac{ | \varphi^{\prime}_{r }(\sigma)| r s}{t \sigma^2 | \Phi_-^{\prime}(\sigma)|^2 }\,. 
$$ 

We will show that the expression $\frac{|\varphi_r'(\sigma)|}{t\sigma^2|\Phi'_-(\sigma)|^2}$ is bounded by $\sigma^{-2}+|\sigma-\sigma_*|^{-2}$. Similarly to \eqref{finest}, this will allow us to conclude that
$$
\int_{|\sigma - \sigma^*| > t^{-\f12}} \frac{|b_2(\sigma;r,s)|| \Phi_-^{\prime \prime }(\sigma)|}{| \Phi_-^{\prime}(\sigma)|^2} \: d \sigma  \les r s \: \la t \ra ^{\f12}\,, 
$$
and together with \eqref{finest} this gives $|I^{-}_2| \les  t^{-3/2}$.  To see this, for $r$ and $t$ fixed, let $z(\sigma)=(2t)^{-1}\varphi_r'(\sigma)$ so that $\Phi_-'(\sigma)=2(\sigma-z(\sigma))$. When $|z(\sigma)|\leq 2$, it is clear that
$$
\frac{|\varphi_r'(\sigma)|}{t\sigma^2|\Phi'_-(\sigma)|^2}\les |\sigma-\sigma_*|^{-2}\,,
$$
on the support of $\tilde{\chi}_c(\sigma)$. On the other hand, when $ z(\sigma)>2$, we first consider $|\sigma - z(\sigma)| \leq z(\sigma)/2$, in which case $\sigma \sim z(\sigma)$, so that we have  
$$ \frac{ | \varphi^{\prime}_{r}(\sigma)|}{t \sigma^2 | \Phi_-^{\prime}(\sigma)|^2 }  \les \frac{z(\sigma)}{ (z(\sigma))^2 |\sigma - \sigma_*|^2} \les |\sigma - \sigma_*|^{-2}\,. 
$$ 
Additionally, when  $|\sigma - z(\sigma)| \geq z(\sigma)/2$ we have
$$ \frac{ | \varphi^{\prime}_{r}(\sigma)|}{t \sigma^2 | \Phi_-^{\prime}(\sigma)|^2 }  \les \frac{z(\sigma)}{ \sigma^2   |\sigma-z(\sigma)|^2} \les \frac{z(\sigma)}{\sigma^2(z(\sigma))^2} \les\sigma^{-2}\,. 
$$
With this, we have completed the estimate of $I_2^-$.

Finally, we consider $I^{+}_2$ where the phase is non-stationary. We apply another integration by parts to this integral to see 
\begin{align*}
|I^+_2| \les \frac{1}{t^2 rs}  \int_c^{\infty} \frac{|b_2^{\prime}(\sigma;r,s)|}{| \Phi_+^{\prime}(\sigma)|} + \frac{|b_2(\sigma;r,s)|| \Phi_+^{\prime \prime }(\sigma)|}{| \Phi_+^{\prime}(\sigma)|^2}\: d \sigma .
\end{align*}
By using the estimates
$$
\frac{|b_2^{\prime}(\sigma;r,s)|}{| \Phi_+^{\prime}(\sigma)|} \les \frac{r s}{\sigma |\sigma + \sigma^*(\sigma) | }\les \sigma^{-2}rs\,,
$$
and 
$$
\frac{|b_2(\sigma;r,s)|| \Phi_+^{\prime \prime }(\sigma)|}{| \Phi_+^{\prime}(\sigma)|^2} \les \frac{s\varphi_r'(\sigma)(1+\sigma^{-2}t^{-1}r)}{(\sigma+z(\sigma))^2}\les \frac{r s}{\sigma^2}+ \frac{rsz(\sigma) }{\sigma^3 z(\sigma) } \les \sigma^{-2} r s\,,
$$
we obtain $|I^+_2| \les t^{-2}  $. 

\end{proof}

With this lemma in hand, we proceed by estimating the term $\tilde{K}_2(r,s;t)$. 
\begin{lemma} $|\tilde{K}_2(r,s;t)| \les t^{-\f32} $. 
\end{lemma}
\begin{proof} Taking into account the symmetry in $\tilde{K}_2(r,s;t)$ with respect to $r$ and $s$, it suffices to estimate the following expression:
$$(rs)^{-1}\int\limits_0^{\infty}e^{it\sigma^2} \tilde{\chi}_c(\sigma) \tilde{\chi}_k(\sigma r)  \chi_k(\sigma s) e(\sigma,r) e(\sigma,s) \: d\sigma. $$  Using Proposition~\ref{prop:lexp}, $ \tilde{\chi}_k(\sigma r) \chi_k(\sigma s) e(\sigma,r) e(\sigma,s)$ can be represented as a sum of the following terms for $\sigma>c$: 
\begin{align} \label{hhhl}
  e^{\pm i ( \sigma r - \frac{ \log(2 \sigma r)}{2\sigma})} e^{ \pm i \theta(\sigma)}  \tilde{\chi}_k(\sigma r)    \chi_k(\sigma s) O_2(\sigma s),\,\,\,\,\ e^{\pm i \sigma r} b_{\pm}(\sigma,r)\tilde{\chi}_k(\sigma r) \chi_k(\sigma s) O_2(\sigma s)   \,.
\end{align}
 We begin by focusing on the  first expression in \eqref{hhhl}. By writing $\varphi_r(\sigma) = \sigma r - \frac{ \log(2 \sigma r)}{2\sigma}$ and $a^{\pm}_s(\sigma)=\tilde{\chi}_c(\sigma) e^{ \pm i \theta(\sigma)} \chi_k(\sigma s) O_2(\sigma s)$, the contribution of this term to the integral is
\begin{align} \label{highsmlg}
(rs) ^{-1}\int\limits_0^{\infty} e^{it \sigma^2 \pm  \varphi_r(\sigma)}  \tilde{\chi}_k(\sigma r) a^{\pm}_s(\sigma) \: d \sigma. 
 \end{align}
 
 We verify that $\varphi_r(\sigma)$ and $a_s^\pm(\sigma)$ satisfy the conditions of Lemma~\ref{lem:osc3}.
 Note that $\varphi_r^{\prime}(\sigma) = r +\frac{[\log(2 \sigma r) - 1]}{2\sigma^2}$. Within the support of $\tilde{\chi}_k(\sigma r)$, it holds that $ \log(2 \sigma r) > 2$, leading to the inequality $ 0 \leq  \varphi_r^{\prime}(\sigma) \les r $. Additionally, we find that $ \varphi_r^{\prime \prime}(\sigma)  = -\sigma^{-3}[  \log(2 \sigma r) - 3/2]$. Therefore, in the support of $\tilde{\chi}_k(\sigma r)$, it's evident that $ \varphi_r^{\prime\prime} (\sigma) <0$ and $|\varphi_r^{\prime \prime}(\sigma)| \les \sigma^{-2} r$. Therefore, by  Lemma~\ref{lem:osc3}, we can conclude that $|\eqref{highsmlg}|\les  t^{-\frac{3}{2}}$ provided that $a^{\pm}_s(\sigma)=  \tilde{\chi}_c(\sigma) O_2(\sigma s )$.
 
 To validate this bound, let us analyze the derivatives of $\theta(\sigma)$. By (5.7.3) of \cite{NIST}, we have the following expansion, 
 \begin{align*}
 \pm \theta(\sigma)= - \text{Arctan}(\pm (2\sigma)^{-1}) \pm \frac{1-\gamma}{2\sigma} + \Im \Big( \sum_{k=1}^{\infty} \frac{(-1)^k ( \zeta(k) -1)}{ ( \pm 2 i \sigma)^k k} \Big)
 \end{align*}
where $\zeta(k)$ is the zeta function and the infinite sum is convergent for $|\sigma| \geq {\f14}$. Since $\theta(\sigma)$ is $C^2$ for $\sigma\in[c,\frac{1}{4}]$, we may conclude from this expansion that $ | \theta^{\prime}|  \les  \sigma ^{-2}$ and $ |\theta^{\prime \prime}| \les \sigma ^{-3}$. Consequently, by \eqref{llest}, we observe that $a^{\pm}_s(\sigma)= \tilde{\chi}_c(\sigma) O_2(\sigma s)$.

We next consider the last term in \eqref{hhhl}. That is, we need to bound 
\begin{align} \label{epsilon}
(rs)^{-1}\int\limits_0^{\infty} e^{i t \sigma^2 \pm i \sigma r}  \tilde{\chi}_k(\sigma r) \tilde{a}_s^{\pm}(\sigma) \: d\sigma\,,
\end{align}
for $\tilde{a}^{\pm}_s(\sigma)= \tilde{\chi}_c(\sigma) \chi_k(\sigma s)b_{\pm}(\sigma,r) O_2(\sigma s) $. Using \eqref{llest}, we obtain $\tilde{a}^{\pm}_s(\sigma)= \tilde{\chi}_c(\sigma)O_2(\sigma s)$. Consequently, applying Lemma~\ref{lem:osc2} with $\varphi_r(\sigma)= \sigma r$, we deduce that $|\eqref{epsilon}| \les t^{-\frac{3}{2}}$. This finishes the proof. 
\end{proof}
Finally, we prove: 
\begin{lemma} $|\tilde{K}_3(r,s;t)| \les t^{-\f32}  $. 
\end{lemma}
\begin{proof} We start by expanding the integrand in $\tilde{K}_3(r,s;t)$ via our approximations. By Proposition~~\ref{prop:lexp}, we may write $\tilde{\chi}_c(\sigma)\tilde{\chi}_k(\sigma r)e(\sigma,r) \tilde{\chi}_k(\sigma s)e(\sigma,s)$ as the sum of the following terms (up to irrelevant constants)
\begin{align}\label{hhhh}
&e^{ \pm i [ \sigma (r+s) - \frac{ \log(2 \sigma r)}{2\sigma}]} \big[e^{ \mp i \frac{ \log(2 \sigma s)}{2\sigma} } e^{ \pm 2 i \theta(\sigma)} + b_{\pm}(\sigma,s) e^{\pm i \theta(\sigma)}\big] \widetilde{\chi}_c(\sigma ) \widetilde{\chi}_k(\sigma r) \widetilde{\chi}_k(\sigma s), \\ 
&e^{ \pm i [ \sigma (r-s) - \frac{ \log(2 \sigma r)}{2\sigma}]}\big[e^{ \pm i \frac{ \log(2 \sigma s)}{2\sigma} }  + b_{\mp}(\sigma,s)e^{\pm i \theta(\sigma)}\big]\widetilde{\chi}_c(\sigma ) \widetilde{\chi}_k(\sigma r) \widetilde{\chi}_k(\sigma s) ,\nn \\
& e^{\pm i \sigma(r+s)} \big[e^{\mp i \frac{ \log(2 \sigma s)}{2\sigma}} e^{ \mp i  \theta(\sigma)}+ b_{\pm}(\sigma,s)\big]b_{\pm}(\sigma,r) \widetilde{\chi}_c(\sigma )\widetilde{\chi}_k(\sigma r) \widetilde{\chi}_k(\sigma s), \nn \\
& e^{\pm i \sigma(r-s)} \big[e^{\pm i \frac{ \log(2 \sigma s)}{2\sigma}} e^{ \mp i  \theta(\sigma)}+ b_{\mp}(\sigma,s)\big]b_{\pm}(\sigma,r) \widetilde{\chi}_c(\sigma )\widetilde{\chi}_k(\sigma r) \widetilde{\chi}_k(\sigma s) \, .\nn 
\end{align}

Each term's contribution to $\tilde{K}_3(r,s;t)$ will be determined by employing Lemma~\ref{lem:osc3}. As per usual, we assume without loss of generality that $r\geq s$.

To analyze the initial two lines within \eqref{hhhh}, we define $\varphi^{\pm}_{r, s}(\sigma) = \sigma (r \pm s) - \frac{ \log(2 \sigma r)}{2\sigma}$ and consider the following two integrals 
\begin{align}
 &(rs)^{-1}\int\limits_0^{\infty} e^{it\sigma^2 \pm i \varphi^+_{r,s}(\sigma)} \tilde{\chi}_c(\sigma)  \tilde{\chi}_k(\sigma r)a^{\pm}_{1,s}(\sigma) \: d \sigma \label{rs1}  \\  & (rs)^{-1}\int\limits_0^{\infty} e^{it\sigma^2 \pm i \varphi^-_{r,s}(\sigma)} \tilde{\chi}_c(\sigma)  \tilde{\chi}_k(\sigma r)a^{\pm}_{2,s}(\sigma) \: d \sigma\,, \label{rs2} 
  \end{align}
where 
\begin{align*}
&a^{\pm}_{1,s}(\sigma) =  \widetilde{\chi}_k(\sigma s) \big[e^{\mp i \frac{\log(2 \sigma s)}{2 \sigma}} e^{\pm 2 i \theta(\sigma)} + b_{\pm}(\sigma,s) e^{\pm i \theta(\sigma)}\big] \\ 
&a^{\pm}_{2,s}(\sigma) =  \widetilde{\chi}_k(\sigma s) \big[e^{\pm i \frac{\log(2 \sigma s)}{2 \sigma}} + b_{\mp}(\sigma,s) e^{\pm i \theta(\sigma)}\big].\, 
 \end{align*}

To confirm the validity of \eqref{condu} for $\varphi^{\pm}_{r,s}$, we start by computing
\begin{align*}
    (\varphi^{\pm}_{r,s})'(\sigma) = (r\pm s)+ (2\sigma)^{-2} [\log(2 \sigma r) -1]\,,
\end{align*}
which yields $0\leq (\varphi^{\pm}_{r, s})' (\sigma) \les  r$. Additionally, $ (\varphi^{\pm }_{r, s})''  (\sigma)= -\sigma^{-3}[\log(2\sigma r) - 3/2]<0$, and $ |(\varphi^{\pm }_{r, s})'' (\sigma)| \les \sigma^{-2}r$. Hence, we can apply Lemma~\ref{lem:osc3}, provided that $a^{\pm}_{1,s}(\sigma)=O_2(\sigma s)$ and $a^{\pm}_{2,s}(\sigma)=O_2(\sigma s)$. For this, we just compute that $\partial_\sigma^j\{ e^{\pm i\frac{ \log(2 \sigma s)}{2\sigma}}\} |=O_2(\sigma s)$ and recall that we have shown that $e^{i\theta(\sigma)}=O_2(\sigma r)$. Combining these estimates with \eqref{llest} yields  $a^{\pm}_{1,s}(\sigma)=O_2(\sigma s)$, $a^{\pm}_{2,s}(\sigma)=O_2(\sigma s)$ in the support of $\tilde{\chi}_c(\sigma)$ implying that $|\eqref{rs1}|, |\eqref{rs2}|\les t^{-\f32}$.

We next consider the last two lines within \eqref{hhhh}. Now we write, 
\begin{align}
 & (rs)^{-1}\int\limits_0^{\infty} e^{it\sigma^2 \pm i \sigma (r+s)} \tilde{\chi}_c(\sigma)  \label{r,s1}\tilde{\chi}_k(\sigma r)a^{\pm}_{1,s}(\sigma) \: d \sigma\,, \\ 
 & (rs)^{-1}\int\limits_0^{\infty} e^{it\sigma^2 \pm i \sigma (r-s)} \tilde{\chi}_c(\sigma)  \tilde{\chi}_k(\sigma r)a^{\pm}_{2,s}(\sigma) \: d \sigma \label{r,s2}\,,
  \end{align}
where 
\begin{align*}
&a^{\pm}_{1,s,r}(\sigma) =  \tilde{\chi}_k(\sigma r)\widetilde{\chi}_k(\sigma s) \big[e^{\mp i \frac{\log(2 \sigma s)}{2 \sigma}}  e^{\pm i \theta(\sigma)} + b_{\pm}(\sigma,s)\big]b_{\pm}(\sigma,r), \\ 
&a^{\pm}_{2,s,r}(\sigma) =  \tilde{\chi}_k(\sigma r)\widetilde{\chi}_k(\sigma s) \big[e^{\pm i \frac{\log(2 \sigma s)}{2 \sigma}}  e^{\mp i \theta(\sigma)}+ b_{\mp}(\sigma,s) \big] b_{\mp}(\sigma,r)
 \end{align*}
The fulfillment of requirement \eqref{condu} in Lemma~\ref{lem:osc3} by $\sigma(r\pm s)$ is evident. Additionally, $b_\pm(\sigma,r)$ are bounded independently of $r$ so employing a similar argument as used for $a^{\pm}_{1,s}$ and $a^{\pm}_{2,s}$ shows that $a^{\pm}_{1,s,r}(\sigma)=O_2(\sigma s)$ and $a^{\pm}_{2,s,r}(\sigma)=O_2(\sigma s)$. Thus, we may use Lemma~\ref{lem:osc3} to see that $|\eqref{r,s1}| , |\eqref{r,s2}|  \les t^{-\f32}$.

This establishes the statement. 
\end{proof}

 \begin{proof}[Proof of Proposition~\ref{prop:high}]
 Combining the bounds for $\tilde{K}_1$, $\tilde{K}_2$ and $\tilde{K}_3$, we obtain the statement. 
 \end{proof}

\appendix

\section{Kernel of the Coulomb evolution} 
The analysis above stems from the following explicit representation of the time evolution of $ H_{0,q} $ for $ q>0 $:
\begin{align} \label{evolrep}
[e^{itH_{0,q} } f](r) = -\frac{q}{2r}  \int\limits_0^{\infty}\int\limits_0^{\infty} e^{it \sigma^2} M_{\frac{iq}{2\sigma},\f12}(2 i\sigma  r)  M_{\frac{iq}{2\sigma},\f12}(2i \sigma  s) sf(s) \sigma^{-1} [e^{\frac{q \pi}{\sigma}} -1]^{-1} \:  d\sigma \:ds 
\end{align}
where $M_{\frac{iq}{2 \sigma}, \f 12} (\cdot)$ is the Whittaker M function (see \cite[Ch~13]{NIST}). Upon substituting $q\sigma$ for $\sigma$, the integral in \eqref{evolrep} transforms into:
\begin{align*}
\frac{q}{2r}\int\limits_0^{\infty}\int\limits_0^{\infty} e^{it q^2 \sigma^2} e(q \sigma,r) e(q \sigma,s) sf(s) \: d\sigma \: ds\,.
\end{align*}
Here, the function $e(q\sigma,r)$ is defined as:
\begin{align*}
e(q\sigma,r) :=-i\sigma^{-\frac{1}{2}} [e^{\frac{\pi}{\sigma}} - 1]^{-\frac{1}{2}} M_{\frac{i}{2 \sigma}, \frac{1}{2}} (2iq\sigma r).
\end{align*}

This representation is obtained by diagonalizing $ r H_{0,q} r^{-1} = -\frac{d^2}{dr^2} +\frac{q}{r}$ via the \emph{distorted Fourier transform}. The purpose of this section is to explain the proof of \eqref{evolrep}, namely
\begin{theorem}\label{dftThm}
	For all $ f\in rC^\infty_{0,\mathrm{rad}}(\bbR^3) $, the equality \eqref{evolrep}holds.  
\end{theorem}
\par
\subsection{Review of Weyl-Titchmarsh theory}
We begin by briefly recalling some of the basic spectral theory of half-line Schr\"odinger operators with a regular left endpoint. In particular, we summarize the construction of the distorted Fourier transform. This theory is well-known and more details may be found in: \cite{CW}, \cite[Ch.9]{CL}, \cite[Sect. XIII.5]{DS}, \cite[Ch. 2]{EK}, \cite{Ev}, \cite{GZ}, \cite[Ch. 10]{Hi}, \cite{Ko}, \cite{Le}, \cite[Ch. 2]{LS}, \cite[Ch. VI]{Na}, \cite[Ch. 6]{Pe}, \cite[Ch.X]{RS}, \cite[Chs. II, III]{Ti}, \cite[Sects. 7–10]{W}. While our potential is not regular at $ 0$ due to the $ \frac{1}{r} $ singularity, the theory of strongly singular potentials is developed in parallel to the treatment of the regular case given here. In particular, we do not use the Herglotz property of the $m$-function, which is standard for the classical theory but may not hold in the strongly singular setting.
\par

Consider the symmetric Schr\"{o}dinger operator 
\begin{align*}
	H=-\frac{d^2}{dx^2} + V(x),\,\,\,\,V=\overline{V}\in L^1_{\text{loc}}(\bbR_+)
\end{align*}
with domain $ \calD(H)=C^2_{0}(\bbR_+) $. We assume that $ V\in L^1(0,1)$  and that it is \emph{limit point} at $ \infty $, that is, for any $ z\in \bbC\setminus \bbR $, the space of solutions to $ Hf=zf $ that are $ L^2 $ at $ \infty $ is at most $ 1 $-dimensional. For instance, it is sufficient (but by no means necessary) to assume that $ V $ is bounded at $ \infty $. For $ \alpha \in [0,\pi)$, let $H_{\alpha}$ be the self-adjoint extension of $H$ with the domain 
\begin{align*}
\calD_{\alpha}:=\{ g \in H^2(\R_+) \mid  \sin(\alpha) g^{\prime}(0) + \cos(\alpha) g(0) =0 \}. 
\end{align*}

We first define $ \phi_\alpha(z,x)$ and $ \theta_\alpha(z,x)$ as the fundamental system of solutions to $H_{\alpha}f = -z^2f$, for $z \in \mathbb{C}$, that satisfy 
\begin{align} \label{phitheta}
\phi_{\alpha}(z,0)= - \theta_{\alpha}^{\prime}(z,0) = -\sin(\alpha),\,\,\,\ \phi^{\prime}_{\alpha}(z,0)=  \theta_{\alpha}(z,0)= \cos(\alpha), \,\,\,\ W(\phi(z, \cdot), \theta(z,\cdot)=1.
\end{align}

 Because $ V $ is $ L^1 $ near $ 0 $, the existence of $ \phi_\alpha$ and $\theta_\alpha $ is assured by Picard iteration as is their analyticity as functions of $ z $. Furthermore, they are real-valued for $ z^2\in\R $.\par
 We next define a \emph{Weyl solution} $ \psi_{\alpha}(z,\cdot)$  near infinity (or zero) to be a non-zero solution to $H_{\alpha}f = -z^2f$ that is $L^2$ near infinity (or zero) for  $\Im \{ z^2 \} \neq 0$. We note that, as long as $V$ is continuous and real valued in $(0,\infty)$, there exist at least one Weyl solution near infinity and at least one Weyl solution near zero, see Theorem~X.6 of \cite{RS}. Since $H$ is in the limit point case at infinity the Weyl solution near infinity is unique (up to scaling), whereas because $H$ is in the limit circle case at zero, all solutions are  Weyl solutions near zero.
Hence, we can uniquely characterize the Weyl solution at infinity as
\begin{align*}
	\psi_\alpha(x,z)=\theta_\alpha(x,z)+m(z)\phi_{\alpha}(x,z)
\end{align*}
where $m(z)= W(\theta_{\alpha}(z,\cdot), \psi_{\alpha}(z,\cdot))$, is the Weyl-m function, which is analytic for $ z^2 \in \mathbb{C} \backslash \bbR$ and in fact Herglotz. Note that this representation is possible because the $ \theta_\alpha $ coefficient of $ \psi_\alpha $ cannot vanish or $ \psi_\alpha $ would be an eigenfunction for non-real $ z^2 $. 

The significance of the Weyl solution is that it allows us write the resolvent kernel or Green's function via
\begin{align}\label{Res}
(H_{\alpha}+z^2)^{-1} f(x) =\int\limits_0^{\infty} [\phi_\alpha(x,z)\psi_\alpha(y,z)\chi_{[0<x<y]} + \phi_\alpha(y,z)\psi_\alpha(x,z)\chi_{[x>y>0]} ]f(y) \:dy. 
\end{align}

With these objects in hand, we are ready to define the distorted Fourier transform:
\begin{pr} 
For $ f\in C_0([0,\infty)) $ let 
\begin{align*}
	[U_{\alpha}f](\lambda)=\int_{0}^{\infty} f(x)\phi_\alpha(x,\lambda)\: dx \,.
\end{align*}
Then we have the following Plancherel theorem
\begin{align*}
	\|f\|_{L^2(\bbR_+)}=\|U_\alpha f\|_{L^2(\bbR,\rho)}
\end{align*}
and inversion formula
\begin{align*}
	f(x)=\lim_{b\to\infty}\int_{-b}^{b} [U_\alpha f](\lambda)\phi_\alpha(x,\lambda)\: \rho(d\lambda)\,.
\end{align*}
Furthermore, for any $F \in C(\R)$ and $ f \in C_0^{\infty}([0,\infty))$, we have
\begin{align}\label{diag}
 [F(H_{\alpha}) f](\cdot)=\int\limits_{\sigma(H_{\alpha})}\int_0^{\infty} F(\sigma^2) \phi_{\alpha}(\sigma,\cdot) \phi_{\alpha}(\sigma,x)  f(x) \:dx  \: \rho(d\sigma).  
\end{align}

\end{pr}
We refer to $ \phi_\alpha $ as the distorted Fourier basis and to $ \rho $ as the associated measure. Typically, one determines $\rho$ via the Herglotz representation of $m(z)$, but we instead insert \eqref{Res} into Stone's formula.  Recall that Stone's formula is given for $\lambda^2 \in\mathbb{R}$ and $ f \in C_0^{\infty}([0,\infty))$ as 
\begin{align*}
 \lim_{\eps\to0+} \frac{1}{\pi i}\int\limits_a^b \la [R(\lambda^2+i\eps)- R(\lambda^2-i\eps) ]f,f \ra \: \lambda \: d\lambda = \la[E(a,b)+\frac12(E(\{a\})+E(\{b\}))] f,f \ra\,,
\end{align*} 
 where $-\infty\leq a\leq b\leq \infty$, $E(\cdot)$ is the spectral resolution of $ H_{\alpha} $, $R(z):=(H_{\alpha}-z)^{-1}$ is the resolvent operator, and we adopt the convention $E(\{\pm\infty\})=0$. The result then follows from the fact that $ \rho(d\lambda) $ is recoverable via the weak-$*$ limit as $ \varepsilon\rightarrow 0  $ of 
 \begin{align*}
	 \frac{\lambda }{\pi i}[m(\lambda^2 +i\varepsilon)-m(\lambda^2-i\varepsilon)]\:d\lambda\,.
 \end{align*}

 \subsection{Proof of Theorem \ref{dftThm}}
This section closely follows \cite{GZ}, which in turn relies on the idea of \cite{HS} to determine the spectral measure via Stone's formula. Let $ \calL_q $ be the half-line Schr\"{o}dinger operator that is unitarily equivalent to $ H_{0,q} $, which we recall is the restriction of the Coulomb Hamiltonian $ H $ to the radial sector. In order to apply the above scheme to $ \calL_q $, we must first make sense of it as a self-adjoint operator. First, recall the following simple consequence of the Kato-Rellich theorem:
 \begin{lemma}[Theorem~X.15 in \cite{RS}]
	Suppose that $ V:\bbR^3\rightarrow \bbR $ is equal to $ V_1+V_2 $ where $ V_1\in L^2(\bbR^3) $ and $ V_2\in L^\infty(\bbR^3) $. Then $ -\Delta +V $ is essentially self-adjoint on $ C_0^\infty(\bbR^3) $ and self-adjoint on $ H^2(\bbR^3) $.
\end{lemma}
Clearly then our Hamiltonian $ H $ is a self-adjoint operator with domain $ H^2(\bbR^3) $.  Recall that $\calL_q$   is the half-line Schr\"{o}dinger operator given by conjugating $ H_{0,q} $ by $ r $. Thus, it is automatically self-adjoint on the domain $ r H^2_{\textrm{rad}}(\bbR^3) $ (where we regard functions in $ H^2_{\textrm{rad}}(\bbR^3) $ as functions on $ \bbR_+ $). In particular, for $ g\in \calD(\calL) $ the function $ \frac{g(r)}{r} $ is continuous at $ r=0 $. 

To compute the resolvent of $ \calL_q=\calL$, first observe that a fundamental system of solutions to $ \calL f=-z^2f $ for $ \Re z^2>0 $ is given by the Whittaker functions \cite[13.14]{NIST}
\begin{align*}
	M_{-\frac{q}{ 2 z}, \frac{1}{2}} ( 2  z r),\,\,W_{-\frac{q}{ 2 z}, \frac{1}{2}} ( 2   zr).
\end{align*}
These are solutions of Whittaker's equation
\begin{align*}
	W''(\omega)+\left(-\frac{1}{4}+\frac{\kappa}{\omega}+\frac{\frac{1}{4}-\mu^2}{\omega^2}\right)W(\omega)=0\,,
\end{align*}
which is related to $ \calL f=-z^2f $ via $r = \frac{\omega}{2 z }$ for $ \kappa=-\frac{q}{2z} $ and $ \mu=\frac{1}{2} $. By \cite[(13.14.6)]{NIST}, we have 
 \begin{align}\label{WM}
 \frac{M_{-\frac{q}{ 2 z}, \frac{1}{2}} ( 2  z r)}{2z} =r e^{-rz} \Big[1+\sum_{s=1}^\infty \frac{(q+2z)(q/2+2z)\cdot\ldots\cdot (q/s+2z)}{s!} r^s \Big]
 \end{align}
and thus  $ \phi (z,r) :=  (2 z )^{-1} M_{-\frac{q}{ 2 z}, \frac{1}{2}} ( 2  z r) $ is the unique solution satisfying the boundary condition for $ \calD(\calL) $ which is normalized so that $ \phi'(0,z)=1 $. It is real analytic for $ z^2\leq 0 $. Moreover, $M_{-\frac{q}{2z},\frac{1}{2}}(2zr)$ is (up to scaling) the unique solution that is $L^2$ at infinity (indeed it decays exponentially by \cite[(13.14.21)]{NIST}. Now, by \cite[(13.14.26)]{NIST} we compute the Wronskian as
 \[
	 W[ \phi(z,r ), W_{-\frac{q}{ 2 z},\frac{1}{2}} ( 2zr )] =  W[M_{-\frac{q}{ 2 z}, \frac{1}{2}} ( \cdot ), W_{-\frac{q}{ 2 z},\frac{1}{2}} ( \cdot )]  
  = - \frac{1}{\Gamma(1+\frac{q}{2z})} = -\frac{2z}{q\Gamma(q/(2z))}.
  \]
  so we set $ \psi (z,r) := -\frac{q}{2z} \Gamma( q/(2z)) W_{-\frac{q}{ 2 z},\frac{1}{2}} ( 2  z r)$ to ensure that $ W[\phi,\psi]=1 $. This gives us the following representation of the resolvent of $ \calL $:
  \begin{pr}\label{dftPr}
	For $\Re z>0$, the resolvent kernel of $ \calL $ is given by
	\begin{align*}
		(\calL+z^2)^{-1}(r,s)=\begin{cases}\phi(z,r)\psi(z,s)\,,&0<r\leq s\\
			\psi(z,r)\phi(z,s)\,,&0<s\leq r
		\end{cases}\,.
	\end{align*}
\end{pr}
\begin{proof}
This follows from the form of the resolvent of a Sturm-Liouville operator and the fact that $ \phi $ is the unique solution satisfying the boundary condition of $ \calD(\calL) $.
\end{proof}
In particular, if we let $-z^2 = \sigma^2 \pm i\varepsilon$ for $\sigma^2 \geq0$, then we obtain
 \begin{align*}
	 (\mathcal{L}-(\sigma^2 \pm i0))^{-1}(r,s) =  \begin{cases}\phi( \pm i \sigma ,r) \psi (\pm i\sigma,s)\,, & 0<r\leq s \\
	 \psi(\pm i\sigma, r) \phi(\pm i\sigma,s )\,, & 0<s\leq r \end{cases} \,.
\end{align*}
These considerations suggest that, as in the classical theory, $ \phi $ should give the distorted Fourier basis whereas $\psi$ should play the role of the Weyl solution.
  
To proceed, we need to define the analog of the function $\theta$ in order to determine the measure $\rho$ associated to $\phi$. Due to the strong singularity of $V$ at zero,  we will not be able to pick the $\theta$ function as in \eqref{phitheta}. In \cite{GZ}, Gesztesy and Zinchenko proved that if $V$ is real valued and $H$ is in the limit point case at both end points then the Weyl-m function exists provided that $Hf = zf$ has a solution $\tilde{\phi}(z,x)$ in $\mathcal{O}$, an open neighborhood of $\R$, that is (a) analytic for $x \in(0,\infty)$ and $ z \in\mathcal{O}$ (b) real valued for $x \in(0,\infty)$ and  $z \in \R$, and (c) in $L^2$ around $x=0$  for $z \in \mathbb{C} \backslash \mathbb{R}$ with sufficiently small $|\Im z|$, see Hypothesis~3.1 of \cite{GZ}. They further showed that if the singularity point of $V$ and the end point of $H$ agree, then the Weyl m-function is a scalar function of $z$. We note, that $ \phi(z,r)$ obeys (a),(b) and (c). However, in contrast to the Coulomb potential the potential examined in \cite{GZ} is singular enough to lead to an essentially self-adjoint $H$.

   Despite the Coulomb potential exhibiting a milder singularity, we apply a similar argument as utilized in \cite{GZ}. Consequently, we determine the fundamental system of solutions to $\mathcal{L} f = -z^2f$ at a reference point $x_0=1$. We take the first solution as $ \phi(z,r)$ and pick a $\theta(z,r)$  such that $ W(\phi(z,r), \theta(z,r))=1$, which we are free to do by Picard iteration from this point. Then, we must have 
  \begin{align}
  \psi(z,r) = \theta(z,r) + m(z) \phi(z,r), \,\,\,\,\ m(z) = W(\theta(z,r), \psi(z,r)).
  \end{align}
  
  As in the proof of Proposition \ref{dftPr}, we use Stone's formula to obtain
  \begin{align} \label{ftkernel}
  [e^{it \mathcal{L} } f](r)= \int\limits_0^{\infty}\int\limits_0^{\infty} e^{it \sigma^2} \phi(i\sigma, r) \phi(i\sigma,s) f(s) \:\rho(d\sigma) \:ds, \,\,\,\,\ \frac{d\rho}{d\sigma}= \frac{2 \sigma}{\pi} \Im ( m(i \sigma)).
  \end{align}

  Here, we have used that $\sigma(\mathcal{L})= \sigma_{\textrm{ac}}(\mathcal{L})=[0,\infty)$ since $\mathcal{L}$ is a positive operator and for any $\sigma \in \R$ we have $\phi(\pm i\sigma,r), \psi(\pm i \sigma, r) \in L_{\delta}^{2} \backslash L^2$ for $\delta>\f12$, where $L^{2}_{\delta}:=\{ (1+r^2)^{-\f {\delta}2} f \in L^2\}$. Note that $L^2_{-\delta}$, being the dual space of $L^2_{\delta}$, is dense in $L^2$.

 Finally, let us determine the density of $ \rho $. Note that  $\theta(z,r)$ has to be real analytic for $ z=i \sigma$, therefore,   we have to have $\theta(i \sigma,r) = \Re(  \psi(i\sigma ,r)) + b(\sigma) \phi(i\sigma,r) $ for some real valued $b(\sigma)$ as
 $$  W[ \theta(i\sigma,\cdot), \phi(i\sigma,\cdot)]= \Re( W[   \psi(i\sigma, \cdot ), \phi(i\sigma,\cdot)])  =1   \,.$$  
 Hence, we compute 
 \begin{align*}
  m(i \sigma ) & =  W(\theta(i\sigma,\cdot), \psi(i\sigma, \cdot) ) \\ 
 &=  2^{-1} W[ \psi(i\sigma,\cdot) +\overline{\psi(i \sigma,\cdot)}, \psi(i\sigma,\cdot)] + b(\sigma) \nn \\
 & =   \Im( \psi (i\sigma,\cdot) \overline{\psi^{\prime}(i\sigma,\cdot)}) + b(\sigma) .\nn 
 \end{align*}
 Moreover, we have as $ r \to 0$,  $$\psi(i \sigma ,r) = -1 + c(i \sigma)r -  r \log r + O(r^{2-})$$  where $\Im(c(i \sigma)) =\sigma -q[\frac{\pi}{2}  + \Im (\psi^{(0)}( 1- iq/(2 \sigma)))]$ and $\psi^0 ( z)$ is the digamma function. Therefore, $\Im( m(i \sigma))=  \Im(c(i \sigma))$
and using $\Im(\psi^0( 1+ i y)) = -(2y)^{-1} + \frac{\pi}{2} \coth (\pi y) $, see \cite[(5.7.5)]{NIST}, we obtain $d \rho(\sigma) =  2q \sigma [e^{\frac{q\pi}{\sigma}} -1]^{-1} $ and this in \eqref{ftkernel} gives 
   \begin{align} \label{finalform}
  [e^{it \mathcal{L} } f](r) &= -\frac{q}{2} \int\limits_0^{\infty}\int\limits_0^{\infty} e^{it \sigma^2} M_{\frac{iq}{2\sigma},\f12}(2 i\sigma  r)  M_{\frac{iq}{2\sigma},\f12}(2i \sigma  s)  f(s)  \sigma^{-1} [e^{\frac{q \pi}{\sigma}} -1]^{-1} \: d\sigma \: ds \\
  & =  -\frac{q}{2}\int\limits_0^{\infty}\int\limits_0^{\infty} e^{itq^2 \sigma^2} M_{\frac{i}{2\sigma},\f12}(2 iq\sigma  r)  M_{\frac{i}{2\sigma},\f12}(2i q\sigma  s)  f(s)  \sigma^{-1} [e^{\frac{\pi}{\sigma}} -1]^{-1} \: d\sigma \: ds  \nn  \,.
  \end{align}
Using the fact that $\mathcal{L}= r H_{0,q}r^{-1}$, we obtain \eqref{evolrep}. 

\bibliographystyle{amsplain}
\bibliography{coulombdispersion}

\begin{bibdiv}
\begin{biblist}

\bib{AS}{book}{
      author={Abramowitz, M.},
      author={Stegun, I.A.},
       title={Handbook of mathematical functions: With formulas, graphs, and
  mathematical tables},
      series={Applied mathematics series},
   publisher={Dover Publications},
        date={1965},
        ISBN={9780486612720},
         url={https://books.google.com/books?id=MtU8uP7XMvoC},
}

\bib{Ag}{book}{
      author={Agmon, S.},
       title={Lectures on exponential decay of solutions of second order
  elliptic equations: Bounds on eigenfunctions of n-body schr{\"o}dinger
  operators},
      series={Mathematical Notes - Princeton University Press},
   publisher={Princeton University Press},
        date={1982},
        ISBN={9780691083186},
         url={https://books.google.com/books?id=EJ5sQgAACAAJ},
}

\bib{BG}{article}{
      author={Beceanu, M.},
      author={Goldberg, M.},
       title={Strichartz estimates and maximal operators for the wave equation
  in {$\mathbb{R}^3$}},
        date={2014},
        ISSN={0022-1236,1096-0783},
     journal={J. Funct. Anal.},
      volume={266},
      number={3},
       pages={1476\ndash 1510},
         url={https://doi.org/10.1016/j.jfa.2013.11.010},
      review={\MR{3146823}},
}

\bib{CW}{incollection}{
      author={Bennewitz, C.},
      author={Everitt, W.N.},
       title={The {T}itchmarsh-{W}eyl eigenfunction expansion theorem for
  {S}turm-{L}iouville differential equations},
        date={2005},
   booktitle={Sturm-{L}iouville theory},
   publisher={Birkh\"{a}user, Basel},
       pages={137\ndash 171},
         url={https://doi.org/10.1007/3-7643-7359-8_7},
      review={\MR{2145081}},
}

\bib{CL}{book}{
      author={Coddington, A.},
      author={Levinson, N.},
       title={Theory of ordinary differential equations},
      series={International series in pure and applied mathematics},
   publisher={R.E. Krieger},
        date={1984},
        ISBN={9780898747553},
         url={https://books.google.com/books?id=pcPhnAAACAAJ},
}

\bib{Wexp}{article}{
      author={Costin, O.},
      author={Donninger, R.},
      author={Schlag, W.},
      author={Tanveer, S.},
       title={Semiclassical low energy scattering for one-dimensional
  {S}chr\"{o}dinger operators with exponentially decaying potentials},
        date={2012},
        ISSN={1424-0637,1424-0661},
     journal={Ann. Henri Poincar\'{e}},
      volume={13},
      number={6},
       pages={1371\ndash 1426},
         url={https://doi.org/10.1007/s00023-011-0155-7},
      review={\MR{2966466}},
}

\bib{Wmain}{article}{
      author={Costin, O.},
      author={Schlag, W.},
      author={Staubach, W.},
      author={Tanveer, S.},
       title={Semiclassical analysis of low and zero energy scattering for
  one-dimensional {S}chr\"{o}dinger operators with inverse square potentials},
        date={2008},
        ISSN={0022-1236,1096-0783},
     journal={J. Funct. Anal.},
      volume={255},
      number={9},
       pages={2321\ndash 2362},
         url={https://doi.org/10.1016/j.jfa.2008.07.015},
      review={\MR{2473260}},
}

\bib{NIST}{misc}{
       title={{\it NIST Digital Library of Mathematical Functions}},
         how={\url{https://dlmf.nist.gov/}, Release 1.1.10 of 2023-06-15},
         url={https://dlmf.nist.gov/},
        note={F.~W.~J. Olver, A.~B. {Olde Daalhuis}, D.~W. Lozier, B.~I.
  Schneider, R.~F. Boisvert, C.~W. Clark, B.~R. Miller, B.~V. Saunders, H.~S.
  Cohl, and M.~A. McClain, eds.},
}

\bib{DSS}{article}{
      author={Donninger, R.},
      author={Schlag, W.},
      author={Soffer, A.},
       title={On pointwise decay of linear waves on a {S}chwarzschild black
  hole background},
        date={2012},
        ISSN={0010-3616,1432-0916},
     journal={Comm. Math. Phys.},
      volume={309},
      number={1},
       pages={51\ndash 86},
         url={https://doi.org/10.1007/s00220-011-1393-8},
      review={\MR{2864787}},
}

\bib{DS}{book}{
      author={Dunford, N.},
      author={Schwartz, J.T.},
       title={Linear operators, part 2: Spectral theory, self adjoint operators
  in hilbert space},
      series={Wiley Classics Library},
   publisher={Wiley},
        date={1988},
        ISBN={9780471608479},
         url={https://books.google.com/books?id=fpJ0uQEACAAJ},
}

\bib{EK}{book}{
      author={Eastham, M.S.P.},
      author={Kalf, H.},
       title={Schr{\"o}dinger-type operators with continuous spectra},
      series={Chapman \& Hall/CRC research notes in mathematics series},
   publisher={Pitman},
        date={1982},
        ISBN={9780273085263},
         url={https://books.google.com/books?id=bbemAAAAIAAJ},
}

\bib{ESI}{article}{
      author={Erd{\'e}lyi, A.},
      author={Kennedy, M.},
      author={McGregor, J.L.},
      author={Swanson, C.A.},
       title={Asymptotic forms of coulomb wave functions, i},
        date={1955},
}

\bib{ESII}{article}{
      author={Erd{\'e}lyi, A.},
      author={Swanson, C.A.},
       title={Asymptotic forms of coulomb wave functions, ii},
        date={1955},
         url={https://books.google.com/books?id=QtjuAAAAMAAJ},
}

\bib{EG}{article}{
      author={Erdo\u{g}an, M.B.},
      author={Schlag, W.},
       title={Dispersive estimates for {S}chr\"{o}dinger operators in the
  presence of a resonance and/or an eigenvalue at zero energy in dimension
  three. {I}},
        date={2004},
        ISSN={1548-159X,2163-7873},
     journal={Dyn. Partial Differ. Equ.},
      volume={1},
      number={4},
       pages={359\ndash 379},
         url={https://doi.org/10.4310/DPDE.2004.v1.n4.a1},
      review={\MR{2127577}},
}

\bib{Ev}{article}{
      author={Everitt, W.~N.},
       title={A personal history of the {$m$}-coefficient},
        date={2004},
        ISSN={0377-0427,1879-1778},
     journal={J. Comput. Appl. Math.},
      volume={171},
      number={1-2},
       pages={185\ndash 197},
         url={https://doi.org/10.1016/j.cam.2004.01.010},
      review={\MR{2077204}},
}

\bib{FFF}{article}{
      author={Fanelli, L.},
      author={Felli, V.},
      author={Fontelos, M.A.},
      author={Primo, A.},
       title={Time decay of scaling critical electromagnetic {S}chr\"{o}dinger
  flows},
        date={2013},
        ISSN={0010-3616,1432-0916},
     journal={Comm. Math. Phys.},
      volume={324},
      number={3},
       pages={1033\ndash 1067},
         url={https://doi.org/10.1007/s00220-013-1830-y},
      review={\MR{3123544}},
}

\bib{F}{article}{
      author={Fulton, C.},
       title={Titchmarsh-{W}eyl {$m$}-functions for second-order
  {S}turm-{L}iouville problems with two singular endpoints},
        date={2008},
        ISSN={0025-584X,1522-2616},
     journal={Math. Nachr.},
      volume={281},
      number={10},
       pages={1418\ndash 1475},
         url={https://doi.org/10.1002/mana.200410689},
      review={\MR{2454944}},
}

\bib{GZ}{article}{
      author={Gesztesy, F.},
      author={Zinchenko, M.},
       title={On spectral theory for {S}chr\"{o}dinger operators with strongly
  singular potentials},
        date={2006},
        ISSN={0025-584X,1522-2616},
     journal={Math. Nachr.},
      volume={279},
      number={9-10},
       pages={1041\ndash 1082},
         url={https://doi.org/10.1002/mana.200510410},
      review={\MR{2242965}},
}

\bib{GS}{article}{
      author={Goldberg, M.},
      author={Schlag, W.},
       title={Dispersive estimates for {S}chr\"{o}dinger operators in
  dimensions one and three},
        date={2004},
        ISSN={0010-3616,1432-0916},
     journal={Comm. Math. Phys.},
      volume={251},
      number={1},
       pages={157\ndash 178},
         url={https://doi.org/10.1007/s00220-004-1140-5},
      review={\MR{2096737}},
}

\bib{HN}{article}{
      author={Hayashi, N.},
      author={Naumkin, P.I.},
       title={Asymptotics for large time of solutions to the nonlinear
  {S}chr\"{o}dinger and {H}artree equations},
        date={1998},
        ISSN={0002-9327,1080-6377},
     journal={Amer. J. Math.},
      volume={120},
      number={2},
       pages={369\ndash 389},
  url={http://muse.jhu.edu/journals/american_journal_of_mathematics/v120/120.2hayashi.pdf},
      review={\MR{1613646}},
}

\bib{Hi}{book}{
      author={Hille, E.},
       title={Lectures on ordinary differential equations},
      series={Addison-Wesley series in mathematics},
   publisher={Addison-Wesley Publishing Company},
        date={1968},
        ISBN={9780201530834},
         url={https://books.google.com/books?id=RDjvAAAAMAAJ},
}

\bib{HS}{incollection}{
      author={Hinton, D.},
      author={Schneider, A.},
       title={On the spectral representation for singular selfadjoint boundary
  eigenvalue problems},
        date={1998},
   booktitle={Contributions to operator theory in spaces with an indefinite
  metric ({V}ienna, 1995)},
      series={Oper. Theory Adv. Appl.},
      volume={106},
   publisher={Birkh\"{a}user, Basel},
       pages={217\ndash 251},
      review={\MR{1729597}},
}

\bib{JSS}{article}{
      author={Journ\'{e}, J.-L.},
      author={Soffer, A.},
      author={Sogge, C.~D.},
       title={Decay estimates for {S}chr\"{o}dinger operators},
        date={1991},
        ISSN={0010-3640,1097-0312},
     journal={Comm. Pure Appl. Math.},
      volume={44},
      number={5},
       pages={573\ndash 604},
         url={https://doi.org/10.1002/cpa.3160440504},
      review={\MR{1105875}},
}

\bib{KP}{article}{
      author={Kato, J.},
      author={Pusateri, F.},
       title={A new proof of long-range scattering for critical nonlinear
  {S}chr\"{o}dinger equations},
        date={2011},
        ISSN={0893-4983},
     journal={Differential Integral Equations},
      volume={24},
      number={9-10},
       pages={923\ndash 940},
      review={\MR{2850346}},
}

\bib{Ko}{article}{
      author={Kodaira, K.},
       title={The eigenvalue problem for ordinary differential equations of the
  second order and {H}eisenberg's theory of {$S$}-matrices},
        date={1949},
        ISSN={0002-9327,1080-6377},
     journal={Amer. J. Math.},
      volume={71},
       pages={921\ndash 945},
         url={https://doi.org/10.2307/2372377},
      review={\MR{33421}},
}

\bib{KT}{article}{
      author={Kova\v{r}\'{\i}k, H.},
      author={Truc, F.},
       title={Schr\"{o}dinger operators on a half-line with inverse square
  potentials},
        date={2014},
        ISSN={0973-5348,1760-6101},
     journal={Math. Model. Nat. Phenom.},
      volume={9},
      number={5},
       pages={170\ndash 176},
         url={https://doi.org/10.1051/mmnp/20149511},
      review={\MR{3264315}},
}

\bib{Le}{article}{
      author={Levinson, N.},
       title={A simplified proof of the expansion theorem for singular second
  order linear differential equations},
        date={1951},
        ISSN={0012-7094,1547-7398},
     journal={Duke Math. J.},
      volume={18},
       pages={57\ndash 71},
         url={http://projecteuclid.org/euclid.dmj/1077476389},
      review={\MR{41313}},
}

\bib{LS}{book}{
      author={Levitan, B.M.},
      author={Sargs\t{\i a}n, I.S.},
       title={Introduction to spectral theory: Selfadjoint ordinary
  differential operators},
      series={Translations of mathematical monographs},
   publisher={American Mathematical Society},
        date={1975},
        ISBN={9780821886632},
         url={https://books.google.com/books?id=hGLJBn4q744C},
}

\bib{Miz}{article}{
      author={Mizutani, H.},
       title={Strichartz estimates for {S}chr\"{o}dinger equations with slowly
  decaying potentials},
        date={2020},
        ISSN={0022-1236,1096-0783},
     journal={J. Funct. Anal.},
      volume={279},
      number={12},
       pages={108789, 57},
         url={https://doi.org/10.1016/j.jfa.2020.108789},
      review={\MR{4156128}},
}

\bib{Na}{book}{
      author={Naimark, M.A.},
       title={Linear differential operators. part ii: Linear differential
  operators in hilbert space},
   publisher={Frederick Ungar Publishing Company},
        date={1968},
         url={https://books.google.com/books?id=tbWPzQEACAAJ},
}

\bib{Nak}{article}{
      author={Nakamura, S.},
       title={Low energy asymptotics for {S}chr\"{o}dinger operators with
  slowly decreasing potentials},
        date={1994},
        ISSN={0010-3616,1432-0916},
     journal={Comm. Math. Phys.},
      volume={161},
      number={1},
       pages={63\ndash 76},
         url={http://projecteuclid.org/euclid.cmp/1104269792},
      review={\MR{1266070}},
}

\bib{Olver}{book}{
      author={Olver, F.},
       title={Asymptotics and special functions},
   publisher={CRC Press},
        date={1997},
        ISBN={9781439864548},
         url={https://books.google.com/books?id=-UBZDwAAQBAJ},
}

\bib{PSV}{misc}{
      author={Pasqualotto, F.},
      author={Shlapentokh-Rothman, Y.},
      author={Van~de Moortel, M.},
       title={The asymptotics of massive fields on stationary spherically
  symmetric black holes for all angular momenta},
        date={2023},
}

\bib{Pe}{book}{
      author={Pearson, D.B.},
       title={Quantum scattering and spectral theory},
      series={Techniques of physics},
   publisher={Academic Press},
        date={1988},
        ISBN={9780125482608},
         url={https://books.google.com/books?id=t6J8AAAAIAAJ},
}

\bib{RS}{book}{
      author={Reed, M.},
      author={Simon, B.},
       title={Ii: Fourier analysis, self-adjointness},
      series={Methods of Modern Mathematical Physics},
   publisher={Elsevier Science},
        date={1975},
        ISBN={9780125850025},
         url={https://books.google.com/books?id=Kz7s7bgVe8gC},
}

\bib{RSchlag}{article}{
      author={Rodnianski, I.},
      author={Schlag, W.},
       title={Time decay for solutions of {S}chr\"{o}dinger equations with
  rough and time-dependent potentials},
        date={2004},
        ISSN={0020-9910,1432-1297},
     journal={Invent. Math.},
      volume={155},
      number={3},
       pages={451\ndash 513},
         url={https://doi.org/10.1007/s00222-003-0325-4},
      review={\MR{2038194}},
}

\bib{survey}{article}{
      author={Schlag, W.},
       title={Spectral theory and nonlinear partial differential equations: a
  survey},
        date={2006},
        ISSN={1078-0947,1553-5231},
     journal={Discrete Contin. Dyn. Syst.},
      volume={15},
      number={3},
       pages={703\ndash 723},
         url={https://doi.org/10.3934/dcds.2006.15.703},
      review={\MR{2220744}},
}

\bib{survey2}{article}{
      author={Schlag, W.},
       title={On pointwise decay of waves},
        date={2021},
        ISSN={0022-2488,1089-7658},
     journal={J. Math. Phys.},
      volume={62},
      number={6},
       pages={Paper No. 061509, 27},
         url={https://doi.org/10.1063/5.0042767},
      review={\MR{4276283}},
}

\bib{Schroe}{article}{
      author={Schrödinger, E.},
       title={Quantisierung als eigenwertproblem},
        date={1926},
     journal={Annalen der Physik},
      volume={384},
      number={4},
       pages={361\ndash 376},
  eprint={https://onlinelibrary.wiley.com/doi/pdf/10.1002/andp.19263840404},
  url={https://onlinelibrary.wiley.com/doi/abs/10.1002/andp.19263840404},
}

\bib{SW}{article}{
      author={Soffer, A.},
      author={Weinstein, M.~I.},
       title={Multichannel nonlinear scattering for nonintegrable equations.
  {II}. {T}he case of anisotropic potentials and data},
        date={1992},
        ISSN={0022-0396,1090-2732},
     journal={J. Differential Equations},
      volume={98},
      number={2},
       pages={376\ndash 390},
         url={https://doi.org/10.1016/0022-0396(92)90098-8},
      review={\MR{1170476}},
}

\bib{Ti}{book}{
      author={Titchmarsh, E.C.},
       title={Eigenfunction expansions associated with second-order
  differential equations},
      series={Eigenfunction Expansions Associated with Second-order
  Differential Equations},
   publisher={Clarendon Press},
        date={1962},
      number={pt. 1},
         url={https://books.google.com/books?id=jv9QAAAAMAAJ},
}

\bib{Weder}{article}{
      author={Weder, R.},
       title={Center manifold for nonintegrable nonlinear {S}chr\"{o}dinger
  equations on the line},
        date={2000},
        ISSN={0010-3616,1432-0916},
     journal={Comm. Math. Phys.},
      volume={215},
      number={2},
       pages={343\ndash 356},
         url={https://doi.org/10.1007/s002200000298},
      review={\MR{1799850}},
}

\bib{W}{book}{
      author={Weidmann, J.},
       title={Spectral theory of ordinary differential operators},
      series={Lecture Notes in Mathematics},
   publisher={Springer},
        date={1987},
        ISBN={9780387179025},
         url={https://books.google.com/books?id=U8gZAQAAIAAJ},
}

\bib{Yaj}{article}{
      author={Yajima, K.},
       title={The {$W^{k,p}$}-continuity of wave operators for
  {S}chr\"{o}dinger operators},
        date={1995},
        ISSN={0025-5645,1881-1167},
     journal={J. Math. Soc. Japan},
      volume={47},
      number={3},
       pages={551\ndash 581},
         url={https://doi.org/10.2969/jmsj/04730551},
      review={\MR{1331331}},
}

\bib{YWB}{article}{
      author={Yost, F.~L.},
      author={Wheeler, J.A.},
      author={Breit, G.},
       title={Coulomb wave functions in repulsive fields},
        date={1936Jan},
     journal={Phys. Rev.},
      volume={49},
       pages={174\ndash 189},
         url={https://link.aps.org/doi/10.1103/PhysRev.49.174},
}

\end{biblist}
\end{bibdiv}


\begin{thebibliography}{9} 
\bibitem{AS}
Abramowitz, M.  and  Stegun, I.~A. \emph{Handbook of mathematical functions with
formulas, graphs, and mathematical tables.} National Bureau of Standards
Applied Mathematics Series, 55. For sale by the Superintendent of Documents,
U.S. Government Printing Office, Washington, D.C. 1964. 
\bibitem{Ag} Agmon, S.  \emph{Lectures on Exponential Decay of Solutions of Second-Order Elliptic Equations: Bounds on Eigenfunctions of N-Body Schrodinger Operations.} (MN-29). Princeton University Press. (1982).
\bibitem{BRV}
 Barcel\'o, J.A., Ruiz, A., Vega,  L.
\emph{Some dispersive estimates for Schrödinger equations
with repulsive potentials}, J. Funct. Anal. 236, 1--24 (2006).
\bibitem{BG} Beceanu, M. and Goldberg, M. \emph{Strichartz Estimates and Maximal Operators for the Wave Equation in
$\R^3$.} J. Funct. Anal. 266 (2014), no. 3, 1476--1510.
\bibitem{CW}  Bennewitz, C. and  Everitt, W.N. \emph{ The Titchmarsh–Weyl eigenfunction expansion theorem for Sturm–Liouville differential equations, in Sturm–Liouville Theory, Past and Present}. W. O. Amrein, A. M. Hinz, and D. B. Pearson (eds.), Birkh\"{a}user, Basel, 2005, 137--171.

\bibitem{CL} Coddington, E. A.  and  Levinson, N., \emph{Theory of Ordinary Differential Equations,} Krieger, Mal-abar, 1985.
\bibitem{Wmain} Costin, O., Schlag, W., Staubach, W. and Tanveer S. \emph{ Semiclassical Analysis of low and zero energy scattering for one dimensional Schr\"odinger operators with inverse square potentials.} J. Funct. Anal. (2008), no. 9 2321--2362.
\bibitem{Wexp}  Costin, O., Schlag, W., Staubach, W. and Tanveer S. \emph{Semiclassical low energy scattering for one-dimensional Schr\"odinger operators with exponentially decaying potentials.} Ann. Henri Poincar\'{e} (2012), 13(6) 1371--1426
\bibitem{DSS} Donninger, R., Schlag, W. and Soffer, A. \emph{On pointwise decay of linear waves on a Schwarzschild black hole background.} Comm. Math. Phys. 309 (2012), no.1 51--86
\bibitem{DS}  Dunford, N. and  Schwartz, J. T., \emph{Linear Operators Part II: Spectral Theory, }Interscience, New York, 1988.
\bibitem{ESI} Erdélyi, A. and Kennedy, M. and McGregor, J. L. and Swanson, C. A . \emph{ Asymptotic Forms of Coulomb Wave Functions, I} (1955). California Institute of Technology, Pasadena, CA.
\bibitem{ESII} Erdélyi, A. and Swanson, C. A. \emph{ Asymptotic Forms of Coulomb Wave Functions, II} (1955). California Institute of Technology , Pasadena, CA
\bibitem{EG}  Erdo\u{g}an, M. B. and Schlag, W. \emph{Dispersive estimates for Schr\"odinger operators in the presence of a resonance and/or an eigenvalue at zero energy in dimension three: I.} Dynamics of PDE 1 (2004), 359--
379.
\bibitem{EK}  Eastham, M. S. P. and  Kalf, H. \emph{Schr\"odinger-type Operators with Continuous Spectra}, Pitman, Boston, 1982.
\bibitem{Ev} Everitt, W. N., \emph{A personal history of the m-coefficient,} J. Comp. Appl. Math. 171, 185--197 (2004).
\bibitem{FFF} Fanelli, L., Felli, V., Fontelos, M.A. and Primo, A. \emph{ Time Decay of Scaling Critical Electromagnetic Schr\"odinger Flows.}Commun. Math. Phys. 324, 1033–1067 (2013).
\bibitem{F}Fulton C. \emph{Titchmarsh–Weyl m-functions for second order Sturm–Liouville problems with two singular endpoints} Math. Nachr. 281, No. 10, 1418 – 1475 (2008) 
\bibitem{GPR}
Germain, P.,  Pusateri, F. and  Rousset, F. 
\emph{The nonlinear Schr\"odinger equation with a potential},
Ann. Inst. Henri Poincare (C) Anal. Non Lineaire 35, 1477--1530 (2018).
\bibitem{GZ} Gesztesy, F. and Zinchenko, M.  \emph{On spectral theory for Schrödinger operators with strongly singular potentials.} Mathematische Nachrichten 279 (2005), 1041--1082.
\bibitem{GS}
Goldberg, M. and Schlag, W. \emph{Dispersive Estimates for Schrödinger Operators in Dimensions One and Three}, Commun. Math. Phys. 251, 157--178 (2004).
\bibitem{HPTV}
Hani, Z., Pausader, B.,  Tzvetkov, N., and  Visciglia, N.
\emph{Modified scattering for the cubic Schrödinger equation on product spaces and applications}, Forum Math. Pi, E4 (2015).
\bibitem{Miz} Mizutani, H. \emph{Strichartz estimates for {S}chr\"{o}dinger equations with
              slowly decaying potentials} J. Funct. Anal. (2020), no.12, 108789
\bibitem{HN} Hayashi, N. and Naumkin, P.I.
\emph{Asymptotics for Large Time of Solutions to the Nonlinear Schrödinger and Hartree Equations},
Amer. J. Math. 120, 369--389 (1998).
\bibitem{Hi}  Hille, E., \emph{ Lectures on Ordinary Differential Equations, Addison–Wesley,} Reading, MA, 1969.
\bibitem{JSS}
 Journ\'e, J.-L., Soffer, A. and Sogge, C.
\emph{Decay estimates for Schrödinger operators},
Comm. Pure and App. Math 44, 573--604 (1991).
\bibitem{KP}
 Kato, J. and Pusateri, F. \emph{A new proof of long-range scattering for critical nonlinear Schrödinger equations}, Differ. Integral Equ.  24, 923--940 (2011).
\bibitem{Ko} Kodaira, K. \emph{The eigenvalue problem for ordinary differential equations of the second order and Heisenberg’s theory of the S-matrices,} Amer. J. Math. 71, 921--945 (1949).
\bibitem{KT }Kovařík, H., Truc, F. (2014).\emph{ Schr\"odinger Operators on a Half-Line with Inverse Square Potentials.} Mathematical Modelling of Natural Phenomena, 9(5), 170-176. 
\bibitem{Le} Levinson, N.  \emph{A simplified proof of the expansion theorem for singular second order linear differential equations,} Duke Math. J. 18, 719--722 (1951).
\bibitem{LS}  Levitan B. M. and  Sargsjan, I. S., \emph{Introduction to Spectral Theory: Selfadjoint Ordinary
Differential Operators,}  Amer. Math. Soc., Providence, RI, 1975
\bibitem{Na}  Naimark, M. A. \emph{ Linear Differential Operators, Part II,} Ungar, New York, 1968.
\bibitem{Nak} Nakamura, S. \emph{Low energy asymptotics for Schr\"odinger operators with slowly decreasing potentials.} Commun.Math. Phys. 161, 63–76 (1994)
\bibitem{Olver} Olver, F. Asymptotics and special functions. CRC Press, 1997.
\bibitem{Ozawa}
Ozawa, T. \emph{Long range scattering for nonlinear Schrödinger equations in one space dimension},
Comm. Math. Phys. 139, 479--493 (1991).
\bibitem{PSV} Pasqualotto, F., Shlapentokh-Rothman, Y. andVan de Moortel, M. \emph{ The asymptotics of massive fields on stationary spherically symmetric black holes for all angular momenta} arXiv:2303.17767 
\bibitem{Pe}  Pearson, D. B., \emph{Quantum Scattering and Spectral Theory, } Academic Press, London, 1988.
\bibitem{Ti} Titchmarsh, E. C.  \emph{Eigenfunction Expansions, Part I, 2nd ed.,} Clarendon Press, Oxford, 1962.
\bibitem{RS} Simon B. and Reed M. \emph{ Methods of modern mathematical physics. Vol. 2, Fourier analysis, self-adjointness} New York: Academic Press, 1975.
\bibitem{RSchlag}
 Rodnianski I. and Schlag, W. \emph{Time decay for solutions of Schrödinger equations
with rough and time-dependent potentials},
Invent. Math. 155, 451--513 (2004).
\bibitem{survey}Schlag, W.  \emph{Spectral theory and nonlinear partial differential equations: a survey.} Discrete Contin. Dyn. Syst. 15 (2006), no. 3, 703--723\
\bibitem{survey2}Schlag, W.  \emph{On pointwise decay of waves} J. Math. Phys. 62, 061509 (2021) 
\bibitem{Schroe} Schr\"{o}dinger, E. \emph{Quantisierung als eigenwertproblem.} Annalen der physik 385.13  437-490 (1926).
\bibitem{SW} Soffer, A. and Weinstein, M. I. \emph{Multichannel nonlinear scattering for nonintegrable equations. II. The case of anisotropic potentials and data.} J. Differ. Eqs. 98 (1992), 376--390
\bibitem{W}  Weidmann, J. \emph{ Spectral Theory of Ordinary Differential Operators, Lecture Notes in Math.}
1258, Springer, Berlin, 1987.
\bibitem{Weder} Weder, R.\emph{ Center manifold for nonintegrable nonlinear Schr\"odinger equations on the line.} Commun. Math. Phys. 215 (2000), 343–356.
\bibitem{Yaj} Yajima, K. \emph{ The $W^{k,p}$-continuity of wave operators for Schr\"odinger operators}. J. Math. Soc. Jpn. 47, 551--581 (1995)
\bibitem{YWB}  Yost, F. L.,  Wheeler, J. A. and  Breit G.
\emph{Coulomb Wave Functions in Repulsive Fields},
Phys. Rev. 49, 174--189 (1936).
\bibitem{HS} Hinton, D. and Schneider, A. \emph{On the spectral representation for singular selfadjoint boundary eigenvalue problems}, in Contributions to Operator Theory in Spaces with an Indefinite Metric, A. Dijksma, I. Gohberg, M. A. Kaashoek, R. Mennicken (eds.), Operator Theory: Advances and Applications, Vol. 106 (1998), pp. 217–251.
\end{thebibliography}

\end{document}